\documentclass[a4paper,12pt,twoside]{article}
\usepackage{amsmath,amssymb,ifthen,bm,mathtools}
\usepackage[amsmath,thmmarks]{ntheorem}
\usepackage{graphicx}
\usepackage{paper}
\usepackage[all]{xy}
\usepackage{mathrsfs}
\usepackage{url}
\usepackage{autobreak}
\usepackage{hhline}
\DeclareMathOperator{\rec}{rec}
\DeclareMathOperator{\Ind}{Ind}
\DeclareMathOperator{\Nr}{Nr}
\DeclareMathOperator{\Art}{Art}

\DeclareMathOperator{\Trd}{Trd}
\DeclareMathOperator{\Nrd}{Nrd}

\DeclareMathOperator{\Lie}{Lie}
\DeclareMathOperator{\cInd}{c-Ind}

\DeclareMathOperator{\JL}{JL}
\DeclareMathOperator{\Nil}{Nil}
\DeclareMathOperator{\Spf}{Spf}
\DeclareMathOperator{\Spa}{Spa}
\DeclareMathOperator{\sgn}{sgn}

\DeclareMathOperator{\ab}{ab}
\DeclareMathOperator{\st}{st}
\DeclareMathOperator{\perf}{perf}

\DeclareMathOperator{\qlog}{qlog}
\DeclareMathOperator{\diag}{diag}

\DeclareMathOperator{\LT}{LT}
\DeclareMathOperator{\Cont}{Cont}
\newcommand{\G}{\mathbb{G}}

\newcommand{\xxi}{\bm{\xi}}

\newcommand{\PA}{\mathfrak{P}_{\mathfrak{A}}}
\newcommand{\PD}{\mathfrak{P}_D}
\newcommand{\UA}{U_{\mathfrak{A}}}
\newcommand{\p}{\mathfrak{p}}
\newcommand{\OO}{\mathcal{O}}
\newcommand{\LL}{\mathcal{L}}
\newcommand{\QQ}{\mathcal{Q}}
\newcommand{\BB}{\mathcal{B}}
\newcommand{\HH}{\mathcal{H}}
\newcommand{\ZZ}{\mathcal{Z}}
\newcommand{\JJ}{\mathcal{J}}
\newcommand{\NN}{\mathcal{N}}

\newcommand{\MM}{\mathcal{M}}
\newcommand{\GG}{\mathcal{G}}
\newcommand{\J}{\,\mathfrak{J}}
\newcommand{\Nb}{\mathbf{N}}
\newcommand{\Ub}{\mathbf{U}}
\newcommand{\Tb}{\bm{T}}
\newcommand{\Xb}{\bm{X}}

\newcommand{\Yf}{\mathfrak{Y}}

\newcommand{\jac}[2]{(\frac{#1}{#2})}
\newcommand{\fw}{\flat,\wedge}

\SelectTips{cm}{11}
\title{Geometric realization of the local Langlands and local Jacquet-Langlands correspondences for $\GL_n$ in a partially ramified case}
\author{Haruto Onishi}

\newcommand{\headertitle}{Geometric realization for $\GL_n$ in a partially ramified case}

\begin{document}

\maketitle

\begin{firstfootnote}
 \ Graduate School of Mathematical Sciences, The University of Tokyo, 8--1 Komaba 3--chome, Meguro-ku, Tokyo, 153--8914, Japan

 E-mail address: \texttt{haruto@ms.u-tokyo.ac.jp}

\end{firstfootnote}

\begin{abstract}
We construct a new family of algebraic varieties over a finite field,
realizing the local Langlands and local Jacquet-Langlands correspondences for $\mathrm{GL}_n$ in their compactly supported \'etale cohomology.
The resulting representations have $L$-parameters which are induced from characters of a partially ramified field extension.
\end{abstract}


\section{Introduction}
Let $K$ be a non-archimedean local field with residue field $\mathbb F_q$ of characteristic
$p$, and let $W_K$ be its Weil group. Let $n\geq 1$ be an integer, and let $D$ be the
central division algebra over $K$ of invariant $1/n$. The local Langlands correspondence for
$\GL_n(K)$ gives a bijection between irreducible supercuspidal representations of
$\GL_n(K)$ and irreducible $n$-dimensional representations of $W_K$, while the local
Jacquet-Langlands correspondence gives a bijection between irreducible essentially
square-integrable representations of $\GL_n(K)$ and irreducible smooth
representations of $D^\times$. In this paper, we study the essentially tame case, in which
the $L$-parameters under consideration are induced from characters of the Weil groups of
tamely ramified extensions $L/K$ of degree $n$. We distinguish the unramified, totally
ramified, and partially ramified cases according to the ramification of $L/K$.

For any prime $l\neq p$, it is known that the $l$-adic cohomology of the Lubin-Tate tower realizes the local
Langlands and local Jacquet-Langlands correspondences simultaneously \cite{MR1719811,MR1876802}.
However, this realization relies on global methods and does not by itself provide an
explicit description of the correspondences. Beginning with the work of Yoshida
\cite{MR2676163}, it has become clear that such explicit descriptions can be obtained by
studying the reductions of suitable formal models of Lubin-Tate spaces. In particular,
Yoshida showed that a Deligne-Lusztig variety appears in the reduction of a semistable
model at the first level, reflecting the realization of depth-zero supercuspidal
representations.

Passing to infinite level, the Lubin-Tate tower admits the structure of a perfectoid space,
whose description is considerably simpler than that of its individual finite levels
\cite{MR3529120}. Boyarchenko and Weinstein \cite{MR3402698} constructed a family of special
affinoids in this space and explicitly determined the defining equations of the algebraic
varieties over finite fields obtained as the reductions of their formal models. They proved
that the compactly supported $l$-adic cohomology of these varieties realizes both
correspondences in the unramified case.
Subsequently, Imai-Tsushima \cite{MR4216689}, and
more generally Tokimoto \cite{MR4126896}, obtained analogous geometric realizations in the
totally ramified case.

The purpose of this paper is to generalize these results simultaneously to the partially ramified case.
We construct explicit affinoids in the infinite-level Lubin-Tate space,
determine the reductions of their formal models as the perfections of algebraic varieties
over finite fields, and show that the compactly supported $l$-adic cohomology of these reductions realizes both correspondences.
More precisely, let $K_f/K$ be the unramified extension of degree $f$, and let $E_e/K$ be a totally tamely ramified extension of degree $e$,
with a uniformizer $\varphi_L$ such that $\varpi=\varphi_L^e$ is a uniformizer of $K$.
We put $L=E_eK_f$ and $n=ef$.
We concentrate on the case of $n>1$.
For a character $\theta$ of $L^\times$ satisfying an appropriate regularity condition,
we construct an irreducible supercuspidal representation of $\GL_n(K)$ and
an irreducible representation of $D^\times$ from the compactly supported
$l$-adic cohomology of an explicitly defined algebraic variety over $\mathbb F_q$.
Under the local Langlands correspondence, their $L$-parameters are induced from $\theta$,
with the normalization specified below.
In contrast to the representation-theoretic approach via Yu's construction,
our geometric method constructs the representations directly from cohomology,
without requiring the construction of the associated Heisenberg-Weil representations.
We return to the relationship with Yu's construction in Section~\ref{sec:rep}.

\paragraph{Outline of the paper and main results.}
In Section~\ref{sec:group}, we define filtrations of compact open subgroups
$\J_{\GL}^m\subset\GL_n(K)$ and
$\J_D^m\subset D^\times$, and use them to construct a subgroup
$J_{m,s}^1$ of $\GL_n(K)\times D^\times\times W_K$.
In Section~\ref{sec:variety}, we define an algebraic variety $X_m$ over
$\overline{\F}_q$ by explicit equations and describe the action
of $J_{m,s}^1$ on its coordinates.
The defining equations combine those appearing in the work of Tokimoto~\cite{MR4126896}
with equations of Boyarchenko-Weinstein type~\cite{boyarchenko2013geometricrealizationspecialcases}.
A notable feature is that the latter contributions occur as a sum with shifted indices.
We denote by $X_m^0$ the connected component of $X_m$
containing the origin.
As we will show in Section~\ref{sec:affinoid}, this component is stable under
$J_{m,s}^1$, and its perfection is equivariantly isomorphic
to the reduction of a formal model of a suitable affinoid
in the infinite-level Lubin-Tate space (Theorem~\ref{thm:red}).

In Section~4, we state our main theorem, which describes
the representations arising from the compactly supported
$l$-adic cohomology of $X_m^0$ and establishes their
compatibility with the local Langlands and local
Jacquet-Langlands correspondences.
To state the theorem, let $m\geq 1$ be coprime to $e$,
and put $s=\lfloor m/2\rfloor$.
Fix a nontrivial additive character $\psi:\mathbb F_q\to\overline{\mathbb Q}_l^\times$. For $\lambda\in\mathbb
F_{q^f}$, put
\[
  \psi_\lambda(x)=\psi\bigl(\operatorname{Tr}_{\mathbb F_{q^f}/\mathbb F_q}(\lambda x)\bigr)
  \quad(x\in\mathbb F_{q^f}).
\]
We say that $\lambda$ satisfies the regularity condition (Condition~\ref{cond}) if $\lambda\neq 0$ and
\[
  \Big(\frac{\lambda^{q^{j-1}}}{\lambda}\Big)^e\neq 1
  \quad(2\leq j\leq f).
\]
Here, the condition involving $j$ is empty when $f=1$.
We regard a character of
$L^\times$ as a character of $W_L$ via local class field theory.
We also fix an additive character
$\Psi:K\to\overline{\Q}_l^\times$
which is trivial on $\p_K$ and whose restriction to
$\OO_K$ induces $\psi$ on $\F_q$.
For each integer $m\ge 1$ and $\mu\in\p_L^{-m}$,
we define a character $\psi_\mu^L$ of $U_L^m$ by
\[
\psi_\mu^L(u)
=\Psi(\Tr_{L/K}(\mu(u-1))),
\quad u\in U_L^m.
\]

Put $U_L^m=1+\p_L^m$ and $U_L^{m,1}=U_L^m\cap\Ker N_{L/K}$.
As explained in Remark~\ref{rem:action}, for a character $\theta$ of
$L^\times$, the space
\[
V_\theta:=
\operatorname{Hom}_{U_L^{m,1}}(\theta,H_c^{n-1}(X_m^0,\overline{\mathbb{Q}}_l))
\]
naturally carries an action of
$L^\times \J_{\GL}^m\times L^\times \J_D^m$.
\begin{thm}[Theorem \ref{thm:main}]
Let $m\geq 1$ be an integer coprime to $e$.
\begin{enumerate}
\item Let $\theta$ be a character of $L^\times$ such that
$\theta|_{U_L^m}
=\psi_{\widetilde{\lambda}\varphi_L^{-m}}^L$,
where $\widetilde{\lambda}\in\OO_{K_f}$ is the
Teichm\"uller lift of $\lambda\in\F_{q^f}$.
Then $V_\theta\neq 0$ if and only if $\lambda$ satisfies Condition~\ref{cond}.
If this condition holds, there exist irreducible representations $\overline\pi_{\GL,\theta}$
of $L^\times\J_{\GL}^m$ and $\overline\pi_{D,\theta}$ of $L^\times\J_D^m$ such that
\[
  V_\theta\cong \overline\pi_{\GL,\theta}^{\vee}\boxtimes\overline\pi_{D,\theta}
\]
as representations of $L^\times\J_{\GL}^m\times L^\times\J_D^m$.
\item If the condition of (i) is satisfied, then the induced representations
    \[\pi_{\GL,\theta}:=\cInd_{L^\times \J_{\GL}^m}^{\GL_n(K)} \overline{\pi}_{\GL,\theta}, \quad \pi_{D,\theta}:=\Ind_{L^\times \J_D^m}^{D^\times} \overline{\pi}_{D,\theta}\]
    are irreducible, and we have the correspondences
    \[\JL(\pi_{\GL,\theta})\cong \pi_{D,\theta},\quad \rec_K(\pi_{\GL,\theta})\cong \Ind_{W_L}^{W_K}\theta\Big(\frac{n-1}{2}\Big).\]
Here $\JL$ and $\rec_K$ denote the local Jacquet-Langlands and
local Langlands correspondences, respectively.
\end{enumerate}
\end{thm}

In Sections~\ref{sec:cohomology-m2} and~\ref{sec:cohomology-m1}, we compute the compactly supported $l$-adic cohomology of $X_m$ from
its defining equations, treating $m>1$ and $m=1$ separately. For $m>1$, the computation is
based on Takeuchi's results on quadratic $l$-adic sheaves. In particular, for $f>1$, the
regularity condition on $\lambda$ is equivalent to the condition that the quadratic sheaf
associated with the Boyarchenko-Weinstein-type part of the defining equation is isogeneous.
This allows us to determine the dimensions of the isotypic components for the action of the
central subgroup $\mathbb F_{q^f}$. In the case $m=1$, higher-degree terms in the defining
equation prevent a direct application of the same quadratic-sheaf argument. Instead, we
eliminate the variables successively to compute the cohomology, obtaining the same dimension
and trace formulas as in the case $m>1$. We prove that
\[
  V_\theta\cong H_c^{n-1}(X_m,\overline{\mathbb Q}_l)[\psi_\lambda^{-1}],
\]
which allows us to determine $\dim V_\theta$. The computation of the top-degree cohomology
also determines the connected components of $X_m$: the variety is geometrically connected
when $e>1$, whereas it has exactly $q$ geometrically connected components when $e=1$.

In Section~\ref{sec:rep}, we study the representation theory of the finite groups $S_{\GL}^m$
and $S_D^m$ arising as quotients of $\J_{\GL}^m$ and $\J_D^m$. We prove the
uniqueness of their irreducible representations with prescribed regular central character.
Comparing their dimensions with that of $V_\theta$, we deduce that $V_\theta$ is irreducible
as a representation of $S_{\GL}^m\times S_D^m$, and hence of $\J_{\GL}^m\times\J_D^m$.
Moreover, these finite groups are closely related to the
Heisenberg groups occurring in Yu's construction: up to character twists,
the representations obtained from our cohomology are induced from the corresponding Heisenberg representations.
Together with this comparison, the trace computations in Sections~\ref{sec:cohomology-m2} and~\ref{sec:cohomology-m1}
allow us to identify their extensions to $L^\times$, thereby characterizing
$\overline\pi_{\GL,\theta}$ and $\overline\pi_{D,\theta}$ explicitly. This
identification with the representations arising from Yu's construction also implies that
$\pi_{\GL,\theta}$ and $\pi_{D,\theta}$ are irreducible and that
$\pi_{\GL,\theta}$ is supercuspidal.

In Section~\ref{sec:affinoid}, we construct explicitly a $J_{m,s}^1$-stable affinoid in the infinite-level
Lubin-Tate space and compute the reduction of a suitable formal model. We prove that this
reduction is $J_{m,s}^1$-equivariantly isomorphic to $(X_m^0)^{\mathrm{perf}}$. Following
the approach of Boyarchenko and Weinstein, we choose a CM point and work in the universal
cover of the associated formal $\mathcal O_K$-module, where the Lubin-Tate space is
described by a determinant-type equation. The partially ramified setting makes the required
approximation of this determinant more complicated than in the unramified case. To overcome
this difficulty, we reorganize the determinant expansion so as to separate the contribution
from the totally ramified direction, to which we can apply Tokimoto's results over $K_f$.
The remaining terms are handled by a more delicate calculation, which produces the shifted
Boyarchenko-Weinstein-type terms in the equation of $X_m$. We refer to Section~\ref{sec:affinoid} for the
details.

Finally, the fact that $X_m^0$ occurs as the reduction of an affinoid in the Lubin-Tate
space does not by itself imply that its cohomology contributes to the cohomology
$H_{\mathrm{LT}}^{n-1}$ of the entire Lubin-Tate tower. In Section~\ref{sec:prf}, we use a result of
Mieda \cite{mieda2016geometricapproachexplicitlocal}, together with the geometric construction of Section~\ref{sec:affinoid} and the
cohomological computations of Section~\ref{sec:cohomology-m2} and~\ref{sec:cohomology-m1}, to establish this connection.
More precisely, we show that
\[
  \pi_{\GL,\theta}^{\vee}\boxtimes\pi_{D,\theta}
  \boxtimes\Ind_{W_L}^{W_K}\theta
\]
occurs in $H_{\mathrm{LT}}^{n-1}$.
By non-abelian Lubin-Tate theory, this proves that our
construction realizes the local Langlands and local Jacquet-Langlands correspondences in
the partially ramified case.

\medbreak
\noindent{\bfseries Acknowledgment.}\quad
The author would like to thank Yoichi Mieda for his support and helpful comments.
This research was supported by FoPM, WINGS Program, the University of Tokyo.

\medbreak
\noindent{\bfseries Notation.}\quad
Let $p$ be any prime number (including $2$) and $K$ be a non-archimedean local field with the ring of integers $\mathcal{O}_K$.
Let $\p_K$ denote the maximal ideal of $\mathcal{O}_K$ and $q$ denote the cardinality of the residue field $\mathcal{O}_K/\p_K$
of characteristic $p$.
We fix an algebraic closure $\overline{K}$ of $K$
and regard all algebraic extensions of $K$ as subfields of $\overline{K}$.
We also fix a prime $l$ distinct from $p$ and an isomorphism $\overline{\Q}_l\cong \C$.
For every finite extension $F/K$, let
$\Art_F\colon F^\times\xrightarrow{\sim}W_F^{\mathrm{ab}}$
denote the local reciprocity map, normalized so that
a uniformizer maps to the geometric Frobenius element.
For $\sigma\in W_F$, we write $\Art_F^{-1}(\sigma)$
for the image of $\sigma$ in $W_F^{\mathrm{ab}}$
under the inverse of this isomorphism.

\section{Notation on groups}\label{sec:group}
\subsection{Subgroups of $\GL_n(K)$}\label{subsec:GL}
Let $E_e$ (resp.\ $K_f$) be a totally tamely ramified (resp.\ unramified) extension of $K$ of degree $e$ (resp. $f$) contained in $\overline{K}$
and let $L$ be the composite $E_e\cdot K_f$.
This is an extension of $K$ of degree $n=ef$.
Throughout this paper, we assume $n>1$.

We write $\Phi$ for the arithmetic Frobenius of unramified extensions of $K$,
and we also write $\Phi$ for the element of $\Gal(L/E_e)$ corresponding to the arithmetic Frobenius of $\Gal(K_f/K)$.
We write $\OO_L$ for the ring of integers of $L$ and $\p_L$ for the maximal ideal of $\OO_L$.
We fix a uniformizer $\varphi_L$ of $E_e$ such that $\varpi=\varphi_L^e$ is a uniformizer of $K$.
The subscript $L$ emphasizes that $\varphi_L$ is also a uniformizer of $L$.
We choose an ordered basis $\alpha_1,\dots , \alpha_f$ of $\OO_{K_f}$ over $\OO_K$ and
take the corresponding $K$-embedding $K_f\hookrightarrow M_f(K)$.
We also take the $K$-embedding $E_e\hookrightarrow M_e(K)$
corresponding to the ordered basis $\varphi_L^{e-1},\dots, 1$.
We fix an identification $M_e(K)\otimes_K M_f(K)\xrightarrow{\sim} M_n(K)$ as
\begin{align*}
    A\otimes B\mapsto \begin{pmatrix}
        Ab_{11} & \cdots & Ab_{1f}\\
        \vdots & \ddots & \vdots \\
        Ab_{f1} & \cdots & Ab_{ff}
    \end{pmatrix},
\end{align*}
so that the two embeddings above give the $K$-embedding $i^G\colon L=E_e\otimes_K K_f\hookrightarrow M_n(K)$.
We write $C_1$ for the orthogonal complement of $L$ with respect to the trace pairing $(x,y)\mapsto \Tr(xy)$ on $M_n(K)$,
and let $\pr_1$ denote the projection to the $L$-part of $M_n(K)$.

The following facts are easily verified.
\begin{itemize}
    \item If we write $t_i\in M_e(K)$ ($i=1,\dots ,e$) for the matrix unit with $1$ in the $(i,i)$-entry and $0$ elsewhere,
    then $t_1,\dots ,t_e$ form an $E_e$-basis of $M_e(K)$.
    \item If we write $s\in M_f(K)$ for the representing matrix of the Frobenius of $K_f/K$,
    then $1,\dots, s^{f-1}$ form a $K_f$-basis of $M_f(K)$.
\end{itemize}
Therefore, if we write $s_i^j=t_i\otimes s^j\in M_n(K)$,
then we have the $L$-basis $\{s_i^j\mid 1\leq i\leq e, 0\leq j\leq f-1\}$ of $M_n(K)$.
We define an $\OO_K$-order $\mathfrak{A}$ in $M_n(K)$ as
\[\mathfrak{A}=\bigoplus_{\substack{1\leq i\leq e \\ 0\leq j\leq f-1}} \OO_L s_i^j.\]
We easily see that this is normalized by $L^\times$.
We write $\PA=\varphi_L\mathfrak{A}$ for the radical of $\mathfrak{A}$.
For each integer $m\geq 1$, we set $\UA^m=1+\PA^m$.
We put
\[\J_{\GL}^m=1+\p_L^m+(C_1\cap \PA^{\lceil m/2\rceil})\subset \GL_n(K).\]
This becomes a compact open subgroup of $\GL_n(K)$ by the following lemma.
\begin{lem}
    For each integer $m\geq 1$, we have the following.
    \begin{enumerate}
        \item For an element $a\in M_n(K)$, if we write
        \[a=\sum_{i=1}^e\sum_{j=0}^{f-1} \gamma_{ij}s_i^j,\]
        where $\gamma_{ij}\in L$, then we have $\pr_1(a)=\frac{1}{e}\sum_{i=1}^e \gamma_{i0}$.
        \item $\pr_1(\PA^m)=\p_L^m$.
        \item $\Tr_{L/K}(\p_L^m)=\p_K^{\lceil m/e \rceil}$.
    \end{enumerate}
\end{lem}
\begin{prf}
    (i) is easy to verify, and (ii) is a direct consequence of (i).
    (iii) is a general fact for tamely ramified extensions.
\end{prf}

\begin{lem}\label{lem:GL-norm}
    For $m\geq 1$, we have $\det({\J_{\GL}^m}) = N_{L/K}(U_L^m) = U_K^{\lceil m/e \rceil}$.
\end{lem}
\begin{prf}
    The identity $N_{L/K}(U_L^m)=U_K^{\lceil m/e\rceil}$ is a general fact for tamely ramified extensions.
    The inclusion $N_{L/K}(U_L^m)\subset \det \J_{\GL}^m$ follows from $U_L^m\subset \J_{\GL}^m$,
    so it remains to prove the opposite inclusion.
    We set $h=\lceil m/2 \rceil$ and $r=\lceil m/e \rceil$.
    Let $g=u+x\in \J_{\GL}^m$ where $u\in U_L^m$ and $x\in C_1\cap \PA^h$.
    Since $\det(g)=\det(u)\det(1+u^{-1}x)$,
    it suffices to check that $\det(1+y)\in U_K^r$ for $y\in C_1\cap \PA^h$.
    We write $\lambda_1,\ldots,\lambda_n\in \overline{K}$ for the eigenvalues of $y$.
    Since $\varpi^{-h}y^e\in \mathfrak{A}$,
    its eigenvalues $\varpi^{-h}\lambda_i^e$ are integral over $\OO_K$,
    so that $v_K(\lambda_i)\geq h/e$ for all $i$.
    If we write $c_j(y)$ for the elementary symmetric polynomial of degree $j$ for $1\leq j\leq n$,
    we have
    \begin{align*}
        \det(1+y)=1+\Tr(y)+\sum_{j=2}^n c_j(y).
    \end{align*}
    We have $\Tr(y)=0$ because $y\in C_1$.
    The above evaluation $v_K(\lambda_i)\geq h/e$ implies that $c_j(y)\in \p_K^{\lceil jh/e\rceil}$,
    so that for $j\geq 2$, we obtain $c_j(y)\in \p_K^{\lceil 2h/e\rceil}\subset \p_K^r$.
    This completes the proof.
\end{prf}

\begin{defn}
    \begin{enumerate}
        \item Let $m\geq 1$ be an odd integer which is coprime to $e$.
        Then we just set $S_{\GL}^m=\F_{q^f}$.
        \item For each even integer $m\ge 2$ coprime to $e$,
        we define a finite group $S_{\GL}^m$ as follows.
        Its underlying set is
        \[\left\{(a_0,(a_{ij}))\in \F_{q^f}\times \prod_{i=1}^e \prod_{j=0}^{f-1}\F_{q^f} \,\middle|\,  \sum_{i=1}^e a_{i0}=0\right\},\]
        and the group structure is defined as follows:
        \begin{align*}
            (a_0,(a_{ij}))\cdot (b_0,(b_{ij}))
            =\left(a_0+b_0+\sum_{i=1}^e\sum_{j=0}^{f-1}a_{ij}b_{i+m/2,-j}^{q^j},\, (a_{ij}+b_{ij})\right).
        \end{align*}
        Here, the subscripts are understood modulo $e$ and $f$, respectively. 
        This is a Heisenberg group whose center is $\F_{q^f}$.
    \end{enumerate}
    
\end{defn}

\begin{lem}
    For each integer $m\geq 1$ which is coprime to $e$, there is a canonical isomorphism
    $j_{\GL}\colon \J_{\GL}^m/\J_{\GL}^{m+1}\xrightarrow{\sim} S_{\GL}^m$. More precisely:
    \begin{enumerate}
        \item Suppose $m$ is odd.
        Then the homomorphism $\J_{\GL}^m\to \F_{q^f}$ defined by
        \[1+g\mapsto \overline{e\pr_1(g)\varphi_L^{-m}}\]
        induces an isomorphism $\J_{\GL}^m/\J_{\GL}^{m+1}\xrightarrow{\sim} \F_{q^f}$.
        \item Suppose $m$ is even.
        Then the homomorphism $\J_{\GL}^m\to S_{\GL}^m$ defined by
        \[1+g=1+\sum_{i=1}^e\sum_{j=0}^{f-1} \gamma_{ij}\varphi_L^{m/2}s_i^j\mapsto (\overline{e\pr_1(g)\varphi_L^{-m}}, (\overline{\gamma}_{ij}))\]
        induces an isomorphism $\J_{\GL}^m/\J_{\GL}^{m+1}\xrightarrow{\sim} S_{\GL}^m$.
    \end{enumerate}
\end{lem}

\begin{prf}
    (i) This is clear considering the fact that for an odd $m$ $\lceil m/2\rceil =\lceil (m+1)/2\rceil$.\\
    (ii) We consider a section of the form
    \[(a_0,(a_{ij}))\mapsto 1+e^{-1}\widetilde{a}_0\varphi_L^m+\sum_{i=1}^e\sum_{j=0}^{f-1} \widetilde{a}_{ij}\varphi_L^{m/2}s_i^j,\]
    where $\widetilde{a}_{ij},\widetilde{b}_{kl}\in \OO_{K_f}$.
    (We can choose a section of this form such that $\sum_i \widetilde{a}_{i0}=0$.)
    Then the multiplication of $\J_{\GL}^m/\J_{\GL}^{m+1}$ can be computed as follows.
    \begin{align*}
        &\left(1+e^{-1}\widetilde{a}_0\varphi_L^m+\sum_{i=1}^e\sum_{j=0}^{f-1} \widetilde{a}_{ij}\varphi_L^{m/2}s_i^j\right)
        \left(1+e^{-1}\widetilde{b}_0\varphi_L^m+\sum_{k=1}^e\sum_{l=0}^{f-1} \widetilde{b}_{kl}\varphi_L^{m/2}s_k^l\right)\\
        &\equiv 1+e^{-1}(\widetilde{a}_0+\widetilde{b}_0)\varphi_L^m
        +\sum_{i=1}^e\sum_{j=0}^{f-1} \widetilde{a}_{ij}\varphi_L^{m/2}s_i^j
        +\sum_{k=1}^e\sum_{l=0}^{f-1} \widetilde{b}_{kl}\varphi_L^{m/2}s_k^l\\
        &\quad +\sum_{i=1}^e\sum_{j=0}^{f-1} \sum_{k=1}^e\sum_{l=0}^{f-1} \widetilde{a}_{ij}\varphi_L^{m/2}s_i^j\widetilde{b}_{kl}\varphi_L^{m/2}s_k^l.
    \end{align*}
    Now we have
    \[s_i^j\widetilde{b}_{kl}=\Phi^j(\widetilde{b}_{kl})s_i^j, \quad s_i^j\varphi_L=\varphi_L s_{i+1}^j,\]
    so that the last sum turns to 
    \begin{align*}
        \sum_{i=1}^e\sum_{j=0}^{f-1} \sum_{k=1}^e\sum_{l=0}^{f-1} \widetilde{a}_{ij}\Phi^j(\widetilde{b}_{kl})\varphi_L^m s_{i+m/2}^j s_k^l
        \equiv\sum_{i=1}^e\sum_{j=0}^{f-1} \widetilde{a}_{ij}\Phi^j(\widetilde{b}_{i+m/2,-j})\varphi_L^m s_{i+m/2}
    \end{align*}
    modulo $C_1\cap \PA^{m/2+1}$, which implies the lemma.
\end{prf}

\subsection{Subgroups of $D^\times$}\label{subsec:D}
Let $D$ be the central division algebra over $K$ with invariant $1/n$.
It is known that there exists a $K$-embedding $K_n\hookrightarrow D$
and an element $\varphi_D$ such that
$D$ is generated by $\varphi_D$ over $K_n$ and
the relations 
\[\varphi_D^n=\varpi, \quad\varphi_D\zeta\varphi_D^{-1}=\zeta^q\text{ for $\zeta\in \mu_{q^n-1}(K_n)$}\]
hold.
We write $\OO_D$ for the maximal order of $D$ and
$\mathfrak{P}_D$ for the radical of $\OO_D$,
and we put $U_D^m=1+\mathfrak{P}_D^m$ for $m\geq 1$.
We take the $K_f$-embedding $i^D\colon L\hookrightarrow D$ satisfying $i^D(\varphi_L)=\varphi_D^f$ and regard $L$ as a $K_f$-subalgebra of $D$.
We write $C_2$ for the orthogonal complement of $L$
with respect to the reduced trace pairing $(x,y)\mapsto \Trd(xy)$ on $D$,
and let $\pr_2$ denote the projection to the $L$-part of $D$.
For each $m\geq 1$, we put
\begin{itemize}
    \item if $m$ is odd, 
    \[
    \J_D^m=1+\mathfrak{p}_L^m+(C_2\cap \mathfrak{P}_D^{f(m-1)/2+1})\subset D^\times,
    \]
    \item if $m$ is even, 
    \[
    \J_D^m=1+\mathfrak{p}_L^m+(C_2\cap \mathfrak{P}_D^{fm/2})\subset D^\times.
    \]
\end{itemize}
These are compact open subgroups of $D^\times$ 
and normalized by $L^\times$.

\begin{lem}
    For each integer $m\geq 1$, we have the following.
    \begin{enumerate}
        \item For an element $d\in D$, if we write $d$ as
        \[d=\sum_{i=0}^{n-1}\beta_i \varphi_D^i,\]
        where $\beta_i\in K_n$, then we have
        \[\pr_2(d)=\frac{1}{e}\sum_{i=0}^{e-1}\Tr_{K_n/K_f}(\beta_{fi})\varphi_L^i.\]
        \item $\pr_2(\PD^m)=\p_L^{\lceil m/f\rceil}$.
        \item $\Trd(\PD^m)=\p_K^{\lceil m/n\rceil}$.
    \end{enumerate}
\end{lem}
\begin{prf}
    We put $\varphi=\begin{pmatrix}
         & I_{n-1} \\
        \varpi & 
    \end{pmatrix}\in M_n(K_n)$.
    As noted in \cite[Remark~3.28]{MR4126896},
    there is a $K$-algebra isomorphism from $D$ to the subring
    \begin{align*}
        \{a\in M_n(K_n)\mid \Phi(a)=\varphi a\varphi^{-1}\},
    \end{align*}
    mapping $x\in K_n$ to $\diag(x,\Phi(x),\ldots,\Phi^{n-1}(x))$ and $\varphi_D$ to $\varphi$.
    The statement (i) is easily deduced from this realization.
    (ii) is a direct consequence of (i), and (iii) can be seen by $\Trd=\Tr_{L/K}\circ\pr_2$.
\end{prf}

\begin{lem}\label{lem:D-norm}
    For $m\geq 1$, we have $\Nrd({\J_D^m}) = N_{L/K}(U_L^m) = U_K^{\lceil m/e \rceil}$.
\end{lem}
\begin{prf}
    The argument is exactly the same as Lemma~\ref{lem:GL-norm}.
\end{prf}

\begin{defn}\label{def:S_D^m}
    \begin{enumerate}
        \item Let $m\geq 1$ be an odd integer which is coprime to $e$.
        We define a finite group $S_D^m$.
        The underlying set is $\F_{q^f} \times \F_{q^n}^{f-1}$
        and the group structure is as follows:
        \begin{itemize}
            \item If $m=1$, then
            \begin{align*}
            (a_0,(a_i))\cdot (b_0,(b_i))=\left(a_0+b_0+\sum_{i=1}^{f-1}\Tr_{\F_{q^n}/\F_{q^f}}(a_i b_{f-i}^{q^i}), (a_i+b_i+\sum_{j=1}^{i-1} a_j b_{i-j}^{q^j})\right).
            \end{align*}
            \item If $m\geq 3$, then
            \begin{align*}
            (a_0,(a_i))\cdot (b_0,(b_i))=\left(a_0+b_0+\sum_{i=1}^{f-1}\Tr_{\F_{q^n}/\F_{q^f}}(a_i b_{f-i}^{q^{f(m-1)/2+i}}), (a_i+b_i)\right).
            \end{align*}
        \end{itemize}
        If $m\geq 3$ and $f\geq 2$, this is a Heisenberg group whose center is $\F_{q^f}$.
        \item Let $m\geq 2$ be an even integer which is coprime to $e$.
        The underlying set of $S_D^m$ is 
        \[\{(a_0,a_1)\in \F_{q^f}\times \F_{q^n}\mid \Tr_{\F_{q^n}/\F_{q^f}}(a_1)=0\},\]
        and the group structure is defined as follows:
        \[(a_0,a_1)\cdot (b_0,b_1)=\left(a_0+b_0+\Tr_{\F_{q^n}/\F_{q^f}}(a_1b_1^{q^{fm/2}}), a_1+b_1\right).\]
        This is a Heisenberg group whose center is $\F_{\! q^f}$.
    \end{enumerate}
\end{defn}


\begin{lem}
    For each integer $m\geq 1$ which is coprime to $e$, there is a canonical isomorphism
    $j_D\colon \J_D^m/\J_D^{m+1}\xrightarrow{\sim} S_D^m$. More precisely:
    \begin{enumerate}
        \item Suppose $m$ is odd.
        Then the homomorphism $\J_D^m\to S_D^m$ defined by
        \[1+d=1+\sum_{i=1}^n \beta_i\varphi_D^{f(m-1)/2+i}\mapsto (\overline{e\pr_2(d)\varphi_L^{-m}}, (\overline{\beta}_i)_{1\leq i\leq f-1}), \quad \text{where $\beta_i\in K_n$,}\]
        induces an isomorphism $\J_D^m/\J_D^{m+1}\xrightarrow{\sim} S_D^m$.
        \item Suppose $m$ is even.
        Then the homomorphism $\J_D^m\to S_D^m$ defined by
        \[1+d=1+\sum_{i=1}^n \beta_i\varphi_D^{fm/2+i-1}\mapsto (\overline{e\pr_2(d)\varphi_L^{-m}}, \overline{\beta}_1), \quad \text{where $\beta_i\in K_n$,}\]
        induces an isomorphism $\J_D^m/\J_D^{m+1}\xrightarrow{\sim} S_D^m$.
    \end{enumerate}
\end{lem}
\begin{prf}
    (i) Consider a section of the form
    \[(a_0,(a_i))\mapsto 1+e^{-1}\widetilde{a}_0\varphi_L^m+\sum_{i=1}^{f-1} \widetilde{a}_i\varphi_D^{f(m-1)/2+i}.\]
    Then the statement can be easily verified.
    The statement (ii) can also be verified by considering a section of the form
    \[(a_0,a_1)\mapsto 1+e^{-1}\widetilde{a}_0\varphi_L^m+\widetilde{a}_1\varphi_D^{fm/2}.\]
    If $j_D(1+d)=0$, then $d\in \PD^{fm/2+1}$ and $\pr_2(d)\in \p_L^{m+1}$.
    Therefore we get $d-\pr_2(d)\in C_2\cap \PD^{fm/2+1}$,
    which implies $\Ker(\J_D^m\to S_D^m)=\J_D^{m+1}$.
\end{prf}

\subsection{Subgroups of $\GL_n(K)\times D^\times \times W_K$}
We write $\Delta_L\colon L\hookrightarrow M_n(K)\times D\times W_K$ for the embedding $a\mapsto (a,a,1)$.
We put $J_m=L^\times \J_{\GL}^m\times L^\times \J_D^m\times W_L$.
Let $\GG^1$ be the kernel of the map $\GL_n(K) \times D^\times \times W_K\to K^\times; (g,d,\sigma)\mapsto \det(g)^{-1}\Nrd(d)\Art_K^{-1}(\sigma)$,
and for a subgroup $J$ of $\GL_n(K) \times D^\times \times W_K$
we write $J^1$ for the intersection $J\cap \GG^1$.
Let $J_{m,s}^1$ be the subgroup of $J_m^1$ generated by the following three types of elements:
\begin{itemize}
    \item $\Delta_L(a)=(a,a,1)$ for $a\in L^\times$,
    \item $(g,d,1)$ for $g\in \J_{\GL}^m$ and $d\in \J_D^m$ with $\det g=\Nrd d$,
    \item $(1,\Art_L^{-1}(\sigma)^{-1},\sigma)$ for $\sigma\in W_L$.
\end{itemize}

\begin{lem}\label{lem:surj}
    The composites
    \begin{align*}
        &J_{m,s}^1\hookrightarrow J_m\twoheadrightarrow L^\times \J_D^m \times W_L,\\ 
        &J_{m,s}^1\hookrightarrow J_m\twoheadrightarrow L^\times \J_{\GL}^m \times W_L,\\ 
        &J_{m,s}^1\hookrightarrow J_m\twoheadrightarrow L^\times \J_{\GL}^m \times L^\times \J_D^m
    \end{align*}
    are surjective.
\end{lem}
\begin{prf}
We call these projections $p_1,p_2$ and $p_3$ in this order.
By Lemmas~\ref{lem:GL-norm} and~\ref{lem:D-norm}, we have
$(d,1)\in \J_D^m\times W_L$ and $(g,1)\in \J_{\GL}^m\times W_L$
are in the images of $p_1$ and $p_2$, respectively.
Since $\Delta_L(L^\times)\subset J_{m,s}^1$,
the subgroup $L^\times\times \{1\}$ of $L^\times \J_D^m \times W_L$ (resp.\ $L^\times \J_{\GL}^m \times W_L$)
is in the image of $p_1$ (resp.\ $p_2$).
For $\sigma\in W_L$, put
$a_\sigma=\operatorname{Art}_L^{-1}(\sigma)$.
We have $p_2(1,a_\sigma^{-1},\sigma)=(1,\sigma)$,
so $p_2$ is surjective.
Moreover, the product
\[
(a_\sigma,a_\sigma,1)
(1,a_\sigma^{-1},\sigma)
=(a_\sigma,1,\sigma)
\]
maps to $(1,\sigma)$ under $p_1$,
so that $p_1$ is also surjective.

We now consider the third projection.
For every $a\in L^\times$, choose
$\sigma\in W_L$ such that $a_\sigma=a$.
The images of $(a,a,1)$ and
$(1,a^{-1},\sigma)$ under $p_3$
are $(a,a)$ and $(1,a^{-1})$, respectively.
Hence the image contains both $(a,1)$
and $(1,a)$ for every $a\in L^\times$.

Next, let $g\in \J^m_{\GL}$.
By Lemma~\ref{lem:GL-norm}, we can choose $a\in U_L^m$
such that $N_{L/K}(a)=\det g$.
Then $(g,a,1)\in J^1_{m,s}$, so that
$(g,a)$ belongs to the image of $p_3$.
Since $(1,a^{-1})$ also belongs
to the image, so does $(g,1)$.
Similarly, for every $d\in \J^m_D$,
Lemma~\ref{lem:D-norm} gives an element $a\in U_L^m$
such that $N_{L/K}(a)=\operatorname{Nrd}d$.
Thus $(a,d,1)\in J^1_{m,s}$, and
multiplying its image by $(a^{-1},1)$
shows that $(1,d)$ belongs to the image.
This shows the surjectivity of $p_3$.
\end{prf}

\section{Algebraic varieties over $\overline{\F}_q$}\label{sec:variety}
\begin{defn}\label{defn:var}
    Let $m\geq 1$ be an integer.
    We set $s=\lfloor m/2\rfloor$.
    We define the closed subvariety $X_m$ of $\A^{n+1}=\Spec \overline{\F}_q[z,y_{ij}\mid 1\leq i\leq e, 1\leq j\leq f]$
    by 
    \[y_{11}+\cdots +y_{e1}=0,\quad z^{q^f}-z=f_m(y_{11},\dots ,y_{e1})+g_m((y_{ij})_{j\neq 1}),\]
    where $f_m$ and $g_m$ are defined as follows
    (the subscripts $i,j$ are understood modulo $e,f$, respectively):
    \begin{enumerate}
        \item We first define $f_m$.
        For each integer $\mu$, we define a set $T(\mu)$ by
        \begin{align*}
            T(\mu)=
            \begin{cases}
                \{(i,j)\in \Z^2 \mid 1\leq i<j\leq e,\, \overline{\mu}< j-i < e-\overline{\mu}\} & \text{if $\overline{\mu}<e/2$},\\
                \{(i,j)\in \Z^2 \mid 1\leq i<j\leq e,\, e-\overline{\mu}\leq j-i \leq \overline{\mu}\} & \text{if $\overline{\mu}\geq e/2$},
            \end{cases}
        \end{align*}
        where $\overline{\mu}$ is the representative of $\mu+e\Z$ in $\{1,\dots, e\}$.
        \begin{enumerate}
            \item If $m$ is odd, then
            \begin{align*}
                f_m=(-1)^{\lfloor m/e\rfloor +1}\sum_{(i,j)\in T(s)}y_{i1}y_{j1}.
            \end{align*}
            \item If $m$ is even, then
            \begin{itemize}
                \item if $\lfloor m/e\rfloor$ is even, then
                \begin{align*}
                \begin{autobreak}
                f_m=\sum_{(i,j)\in T(s)}(y_{i1}^{q^f}-y_{i1})(y_{j1}^{q^f}-y_{j1})
                + \sum_{i=1}^e (y_{i1}^{q^f}-y_{i1})y_{\overline{i+s},1}^{q^f}
                \end{autobreak}
                \end{align*}
                \item if $\lfloor m/e\rfloor$ is odd, then
                \begin{align*}
                \begin{autobreak}
                f_m=-\sum_{(i,j)\in T(-s)}(y_{i1}^{q^f}-y_{i1})(y_{j1}^{q^f}-y_{j1})
                + \sum_{i=1}^e (y_{i1}^{q^f}-y_{i1})y_{\overline{i+s},1}
                \end{autobreak}
                \end{align*}
            \end{itemize}
        \end{enumerate}
        \item We next define $g_m$.
        If $f=1$, we set $g_m=0$.
        In the following, suppose $f\ge 2$.
        For a square matrix $M=(m_{ij})$ of degree~$f$,
        let $|M|^*$ denote the sum $\sum_{\sigma\in \mathfrak{S}_f\setminus \{1\}} \sgn \sigma \prod_{i=1}^f m_{i\sigma(i)}$.
        \begin{itemize}
            \item If $m=1$, then
            \[
            g_1=-\sum_{i=1}^e
            \begingroup
            \renewcommand{\arraystretch}{1.4}
                \begin{vmatrix}
                    1 & y_{i,2}^q & y_{i,3}^{q^2} & \cdots & y_{i,f}^{q^{f-1}}\\
                    y_{i,f}^{q^f}-y_{\overline{i+1},f} & 1 & y_{i,2}^{q^2} & \cdots & y_{i,f-1}^{q^{f-1}}\\
                    y_{i,f-1}^{q^f}-y_{\overline{i+1},f-1} & 0 & 1 &  \cdots & y_{i,f-2}^{q^{f-1}}\\
                    \vdots & \vdots & \vdots & \ddots & \vdots \\
                    y_{i,2}^{q^f}-y_{\overline{i+1},2} & 0 & 0 & \cdots & 1
                \end{vmatrix}^*
            \endgroup.
            \]
            \item If $m\geq 3$ is odd, then
            \[
            g_m=-\sum_{i=1}^e
            \begingroup
            \renewcommand{\arraystretch}{1.4}
                \begin{vmatrix}
                    1 & y_{\overline{i-s},2}^q & y_{\overline{i-s},3}^{q^2} & \cdots & y_{\overline{i-s},f}^{q^{f-1}}\\
                    y_{i,f}^{q^f}-y_{\overline{i+1},f} & 1 & 0 & \cdots & 0\\
                    y_{i,f-1}^{q^f}-y_{\overline{i+1},f-1} & 0 & 1 & \cdots & 0\\
                    \vdots & \vdots & \vdots & \ddots & \vdots \\
                    y_{i,2}^{q^f}-y_{\overline{i+1},2} & 0 & 0 & \cdots & 1
                \end{vmatrix}^*
            \endgroup.
            \]
            \item If $m$ is even, then
            \[
            g_m=-\sum_{i=1}^e
            \begingroup
            \renewcommand{\arraystretch}{1.4}
                \begin{vmatrix}
                    1 & y_{\overline{i+s},2}^{q^{f-1}} & y_{\overline{i+s},3}^{q^{f-2}} & \cdots & y_{\overline{i+s},f}^{q}\\
                    y_{i,f}^{q^f}-y_{i,f} & 1 & 0 & \cdots & 0\\
                    y_{i,f-1}^{q^f}-y_{i,f-1} & 0 & 1 & \cdots & 0\\
                    \vdots & \vdots & \vdots & \ddots & \vdots \\
                    y_{i,2}^{q^f}-y_{i,2} & 0 & 0 & \cdots & 1
                \end{vmatrix}^*
            \endgroup.
            \]
        \end{itemize}
    \end{enumerate}
\end{defn}

\begin{defn}
    We use the notation of Definition~\ref{defn:var}.
    For each $m\geq 1$ coprime to $e$,
    we define an action of $S_{\GL}^m\times S_D^m$ on $X_m$ as follows.
    \begin{enumerate}
        \item We first describe the right action of $S_{\GL}^m$.
        \begin{itemize}
            \item If $m$ is odd, then the action of $a_0 \in S_{\GL}^m=\F_{q^f}$ is
            \[y_{ij}\mapsto y_{ij},\quad z\mapsto z+a_0.\]
            \item If $m$ is even, then the action of $(a_0,(a_{ij}))\in S_{\GL}^m$ is
            \[y_{ij}\mapsto y_{ij}+a_{i,j-1},\quad z\mapsto z+\sum_{i=1}^e \sum_{j=0}^{f-1}a_{ij}^{q^{-j}}y_{i-s,1-j}+a_0.\]
        \end{itemize}
        \item We next describe the left action of $S_D^m$.
        \begin{itemize}
            \item If $m=1$, then the action of $(a_0,(a_i))\in S_D^1$ is
            \begin{align*}
                &y_{i1}\mapsto y_{i1},\quad y_{ij}\mapsto y_{ij}+\sum_{1<l<j}a_{j-l}^{q^{f(i-1)+1-j}}y_{il}+a_{j-1}^{q^{f(i-1)+1-j}}\,(j\geq 2),\\
                &z\mapsto z+\sum_{i=1}^e \sum_{j=2}^f a_{f+1-j}^{q^{f(i-2)}}y_{ij}+a_0.
            \end{align*}    
            \item If $m\geq 3$ is odd, then the action of $(a_0,(a_i))\in S_D^m$ is
            \begin{align*}
                &y_{i1}\mapsto y_{i1},\quad y_{ij}\mapsto y_{ij}+a_{j-1}^{q^{f(i-s-1)+1-j}}\,(j\geq 2),\\
                &z\mapsto z+\sum_{i=1}^e \sum_{j=2}^f a_{f+1-j}^{q^{f(i-m-1)}}y_{ij}+a_0.
            \end{align*}
            \item If $m$ is even, then the action of $(a_0,a_1)\in S_D^m$ is
            \begin{align*}
                &y_{i1}\mapsto y_{i1}+a_1^{q^{f(i-s-1)}}, \quad y_{ij}\mapsto y_{ij}\, (j\geq 2),\\
                &z\mapsto z+\sum_{i=1}^e a_1^{q^{f(i-m-1)}}y_{i1}+a_0.
            \end{align*}
        \end{itemize}
    \end{enumerate}
    In addition, we write $\Frob_{q^f}$ for the automorphism of $X_m$ over $\F_{q^f}$ (not $\overline{\F}_q$)
    defined by $y_{ij}\mapsto y_{ij},\, z\mapsto z$ and $a\mapsto a^{q^{-f}}$ for $a\in \overline{\F}_q$.
\end{defn}

\begin{defn}\label{defn:action-on-reduction}
    For each $m\geq 1$ coprime to $e$,
    we define an action of $J_{m,s}^1$ on $X_m$ as follows.
    \begin{itemize}
        \item $(\varphi_L,\varphi_L,1)$ acts by $y_{ij}\mapsto y_{\overline{i-1},j},\, z\mapsto z$.
        \item $(u,u,1)$ with $u\in U_L$ acts by
        \begin{itemize}
            \item if $m$ is odd, $y_{ij}\mapsto \overline{u}^{1-q^{1-j}}y_{ij}, \, z\mapsto z$.
            \item if $m$ is even, $y_{ij}\mapsto \overline{u}^{q^{j-1}-1}y_{ij}, \, z\mapsto z$.
        \end{itemize}
        \item $(g,d,1)$ with $\det g=\Nrd d$ acts by $(j_{\GL}(g),j_D(d^{-1}))\in S_{\GL}^m\times S_D^m$.
        \item For $\sigma\in W_L$, we set $a_\sigma=\Art_L^{-1}(\sigma)$, $n_\sigma=v_L(a_\sigma)$ and $u_\sigma=a_\sigma \varphi_L^{-n_\sigma}$.
        Let $\chi_2\colon \F_{q^f}^\times\to\{\pm1\}$
        denote the quadratic character if $p\neq 2$,
        and the trivial character if $p=2$.
        The element $(1,a_\sigma^{-1},\sigma)$ acts by the composite of $\Frob_{q^f}^{n_\sigma}$ and
        \begin{itemize}
            \item if $m$ is odd, $y_{i1}\mapsto \chi_2(\overline{u}_\sigma)y_{i1},\, y_{ij}\mapsto \overline{u}_\sigma^{q^{1-j}-1}y_{ij}\, (j\geq 2),\, z\mapsto z$.
            \item if $m$ is even, $y_{ij}\mapsto y_{ij},\, z\mapsto z$.
        \end{itemize} 
    \end{itemize}
\end{defn}

\begin{rem}
    Even though $\J_{\GL}^m\times \J_D^m$ is not contained in $J_{m,s}^1$,
    the action of $\J_{\GL}^m\times \J_D^m$ on $X_m$ can be defined as the action of $(j_{\GL}(g),j_D(d^{-1}))$.
    We have to be careful about the direction of these actions.
    Since the action of $S_{\GL}^m$ is a right action and that of $S_D^m$ is a left action,
    the action of $\J_{\GL}^m\times \J_D^m$ defined in this way is a right action on both factors,
    so that the action on the cohomology is a left action.
\end{rem}

\begin{rem}\label{rem:conn}
    Let $X_m^0$ be the connected component of $X_m$ containing the origin.
    We will prove in Proposition~\ref{prop:connected} that $X_m=X_m^0$ if $e>1$,
    while $X_m$ has $q$ connected components if $e=1$.
    Theorem~\ref{thm:red} and its proof imply that $X_m^0$ is $J_{m,s}^1$-stable.
\end{rem}

\begin{rem}\label{rem:action}
    By Remark~\ref{rem:conn},
    $H_c^{n-1}(X_m^0,\overline{\Q}_l)$ becomes a representation of $J_{m,s}^1$.
    By Lemma~\ref{lem:surj}, for a character $\theta$ of $L^\times$,
    $\Hom_{U_L^{m,1}}(\theta, H_c^{n-1}(X_m^0,\overline{\Q}_l))$ becomes a representation of $L^\times \J_{\GL}^m\times L^\times \J_D^m$.
    Here, $U_L^{m,1}=\{u\in U_L^m \mid N_{L/K}(u)=1\}$ and
    its action on $H^{n-1}_c(X_m^0,\overline{\Q}_l)$ is taken as that of $(1,1,\Art_L(U_L^{m,1}))\subset J_{m,s}^1$. 
    The action of $(g,d)\in L^\times \J_{\GL}^m\times L^\times \J_D^m$ on $f\in \Hom_{U_L^{m,1}}(\theta, H_c^{n-1}(X_m^0,\overline{\Q}_l))$ is given by
    \[((g,d)f)(v)=(g,d,\sigma)(f(\sigma^{-1} v)) \quad (v\in \theta),\]
    where $\sigma\in W_L$ is an element such that $(g,d,\sigma)\in J_{m,s}^1$.
    This action is independent of the choice of $\sigma$ 
    because the ambiguity of $\sigma$ is exactly $U_L^{m,1}$, on which $f$ is equivariant.
    We note that the action of $(1,1,\Art_L(u))$ for $u\in U_L^{m,1}$ on $X_m$ is the same as
    that of $(u^{-1},1,1)$ and $(1,u,1)$, since we have the decompositions
    \begin{align*}
        (1,1,\Art_L(u))&=(1,u,1)(1,u^{-1},\Art_L(u))\\
        &=(u^{-1},1,1)(u,u,1)(1,u^{-1},\Art_L(u)),
    \end{align*}
    and $(u,u,1)$ and $(1,u^{-1},\Art_L(u))$ act on $X_m$ trivially.
\end{rem}

\section{Main Theorem}\label{sec:main}
We fix a nontrivial character $\psi$ of $\F_q$.
For each $\lambda\in \F_{q^f}$, we define a character $\psi_\lambda$ of $\F_{q^f}$
by $\psi_\lambda(x)=\psi(\Tr_{\F_{q^f}/\F_{q}}(\lambda x))$.
We consider the following condition for $\lambda$.
\begin{cond}\label{cond}
    We require that $\lambda\neq 0$ and
    \[\bigg(\frac{\lambda^{q^{j-1}}}{\lambda}\bigg)^e\neq 1\]
    for all $j=2,\dots,f$.
\end{cond}
We also fix an additive character
$\Psi:K\to\overline{\Q}_l^\times$
which is trivial on $\p_K$ and whose restriction to
$\OO_K$ induces $\psi$ on $\F_q$.
For each integer $m\ge 1$ and $\mu\in\p_L^{-m}$,
we define a character $\psi_\mu^L$ of $U_L^m$ by
\[
\psi_\mu^L(u)
=\Psi(\Tr_{L/K}(\mu(u-1))),
\quad u\in U_L^m.
\]
The following lemma is easy to verify.
\begin{lem}
    Let $\theta$ be a character of $L^\times$ such that
    \[
    \theta|_{U_L^m}
    =\psi_{\widetilde{\lambda}\varphi_L^{-m}}^L,
    \]
    where $\widetilde{\lambda}\in\OO_{K_f}$ is the
    Teichm\"uller lift of $\lambda\in\F_{q^f}$.
    Suppose $\lambda\neq 0$.
    Then Condition~\ref{cond} for $\lambda$ implies the irreducibility of the induced representation $\Ind_{W_L}^{W_K}\theta$ of $W_K$.
    Here we view $\theta$ as a character of $W_L$ via $\Art_L$.
\end{lem}

The following is our main theorem.
\begin{thm}\label{thm:main}
    Let $m\geq 1$ be an integer coprime to $e$.
    \begin{enumerate}
    \item Let $\theta$ be a character of $L^\times$ such that
    $\theta|_{U_L^m}
    =\psi_{\widetilde{\lambda}\varphi_L^{-m}}^L$,
    where $\widetilde{\lambda}\in\OO_{K_f}$ is the
    Teichm\"uller lift of $\lambda\in\F_{q^f}$.
    Then
    \[\Hom_{U_L^{m,1}}(\theta, H^{n-1}_c(X_m^0,\overline{\Q}_l))\neq 0\]
    if and only if $\lambda$ satisfies Condition~\ref{cond}, and in that case,
    there are irreducible representations $\overline{\pi}_{\GL,\theta}$ of $L^\times \J_{\GL}^m$
    and $\overline{\pi}_{D,\theta}$ of $L^\times \J_D^m$ such that
    \[\Hom_{U_L^{m,1}}(\theta, H^{n-1}_c(X_m^0,\overline{\Q}_l))\cong \overline{\pi}_{\GL,\theta}^\vee \boxtimes \overline{\pi}_{D,\theta}\]
    as $L^\times \J_{\GL}^m\times L^\times \J_D^m$-representations.
    \item If the condition of (i) is satisfied, then the induced representations
    \[\pi_{\GL,\theta}:=\cInd_{L^\times \J_{\GL}^m}^{\GL_n(K)} \overline{\pi}_{\GL,\theta}, \quad \pi_{D,\theta}:=\Ind_{L^\times \J_D^m}^{D^\times} \overline{\pi}_{D,\theta}\]
    are irreducible, and we have the correspondences
    \[\JL(\pi_{\GL,\theta})\cong \pi_{D,\theta},\quad \rec_K(\pi_{\GL,\theta})\cong \Ind_{W_L}^{W_K}\theta\Big(\frac{n-1}{2}\Big).\]
    Here $\JL$ and $\rec_K$ denote the local Jacquet-Langlands and
    local Langlands correspondences, respectively.
    \end{enumerate}
\end{thm}

\begin{rem}
    The geometric interpretation of $X_m^0$ will be established in Section~\ref{sec:affinoid}.
    More precisely, we will show in Theorem~\ref{thm:red} that the perfection of $X_m^0$ arises
    as the reduction of a suitable equivariant formal model of an affinoid in the infinite-level Lubin-Tate space.
    This geometric realization will play a crucial role in the proof of the correspondence assertions in Theorem~\ref{thm:main}.
\end{rem}

This geometric construction of $\overline{\pi}_{\GL,\theta}$ and $\overline{\pi}_{D,\theta}$
gives us the explicit characterizations described in the following theorem.
Here, $(\frac{x}{q^f})$ is the Jacobi symbol,
and for $p\neq 2$ and $\lambda\neq 0$,
$\mathfrak{g}(\psi_\lambda)$ denotes the quadratic Gauss sum
$\sum_{x\in \F_{q^f}} \psi_\lambda(x^2)=\sum_{x\in \F_{q^f}^\times} (\frac{x}{q^f})\psi_\lambda (x)$.
When $m$ is odd, we set
\begin{align*}
    \mathfrak{g}_\lambda=
    \begin{cases}
        \jac{q^f}{e}q^{f(e-1)/2} & \text{if $e$ is odd and $p=2$,}\\
        \jac{e}{q^f}\mathfrak{g}(\psi_\lambda)^{e-1} & \text{if $e$ is odd and $p\neq 2$,}\\
        -\jac{(-1)^s 2e}{q^f}\mathfrak{g}(\psi_\lambda)^{e-1} & \text{if $e$ is even.}
    \end{cases}
\end{align*}
Also, for $u\in U_L$ we define a subset $S_u$ of $\{2,\ldots, f\}$ as
\[S_u:=\{j\in \{2,\ldots, f\}\mid \overline{u}\in \F_{q^{j-1}}\}.\]
If $f=1$, we set $S_u=\emptyset$.
\begin{rem}
    We note that $S_u$ is closed under the involution $j\mapsto j^*=f-j+2$.
    Indeed, there exists a minimal integer $a\mid f$ such that $\overline{u}\in \F_{q^a}$.
    Then $\overline{u}\in \F_{q^{j-1}}$ if and only if $a\mid j-1$.
    However, since $(j-1)+(j^*-1)=f$, $a\mid j-1$ is equivalent to $a\mid j^*-1$,
    which implies the claim.
\end{rem}
\begin{thm}\label{thm:characterize}
    Let $m\geq 1$ be an integer which is coprime to $e$.
    Let $\theta$ be a character of $L^\times$ such that
    $\theta|_{U_L^m}
    =\psi_{\widetilde{\lambda}\varphi_L^{-m}}^L$,
    where $\widetilde{\lambda}\in\OO_{K_f}$ is the
    Teichm\"uller lift of $\lambda\in\F_{q^f}$.
    Suppose that $\lambda$ satisfies Condition~\ref{cond}.
    Let $\overline{\pi}_{\GL,\theta}$ and $\overline{\pi}_{D,\theta}$
    be the irreducible representations of $L^\times \J_{\GL}^m$ and $L^\times \J_D^m$ respectively,
    such that 
    \[\Hom_{U_L^{m,1}}(\theta, H^{n-1}_c(X_m^0,\overline{\Q}_l))\cong \overline{\pi}_{\GL,\theta}^\vee \boxtimes \overline{\pi}_{D,\theta}\]
    as $L^\times \J_{\GL}^m\times L^\times \J_D^m$-representations.
    \begin{enumerate}
        \item The irreducible representation $\overline{\pi}_{\GL,\theta}$ of $L^\times \J_{\GL}^m$
        is characterized by the following properties:
        \begin{itemize}
            \item $\overline{\pi}_{\GL,\theta}|_{\J_{\GL}^m}$ is the inflation of the unique irreducible representation
            of $\J_{\GL}^m/\J_{\GL}^{m+1}$ with the central character $\psi_\lambda$.
            \item $U_L^1$ acts by $\theta$.
            \item For the action of $\varphi_L$, we have
            \begin{itemize}
                \item If $m$ is odd, then $\Tr \overline{\pi}_{\GL,\theta}(\varphi_L)=(-1)^{n-1}q^{-ef(f-1)/2}\mathfrak{g}_{-\lambda}^{-1} \theta(\varphi_L)$.
                \item If $m$ is even, then $\Tr \overline{\pi}_{\GL,\theta}(\varphi_L)=(-1)^{f-1}q^{-(e-1)f^2/2} \theta(\varphi_L)$.
            \end{itemize}
            \item If $u\in U_L$, then 
            \begin{itemize}
                \item If $m$ is odd, then $\Tr \overline{\pi}_{\GL,\theta}(u)=\chi_2(\overline{u})^{e-1}\theta(u)$.
                \item If $m$ is even, then $\Tr \overline{\pi}_{\GL,\theta}(u)=(-1)^{e(f-1-|S_u|)}q^{(n|S_u|+(e-1)f)/2} \theta(u)$.
            \end{itemize}
        \end{itemize}
        \item The irreducible representation $\overline{\pi}_{D,\theta}$ of $L^\times \J_D^m$
        is characterized by the following properties:
        \begin{itemize}
            \item $\overline{\pi}_{D,\theta}|_{\J_D^m}$ is the inflation of the unique irreducible representation
            of $\J_D^m/\J_D^{m+1}$ with the central character $\psi_\lambda$.
            \item $U_L^1$ acts by $\theta$.
            \item For the action of $\varphi_L^{-1}$, we have
            \begin{itemize}
                \item If $m$ is odd, then $\Tr \overline{\pi}_{D,\theta}(\varphi_L^{-1})=(-1)^{f-1}q^{(e+1)f(f-1)/2}\mathfrak{g}_{-\lambda} \theta(\varphi_L)^{-1}$.
                \item If $m$ is even, then $\Tr \overline{\pi}_{D,\theta}(\varphi_L^{-1})=(-1)^{e(f-1)}q^{(n-1)f/2}\theta(\varphi_L)^{-1}$.
            \end{itemize}
            \item If $u\in U_L$, then 
            \begin{itemize}
                \item If $m$ is odd, then $\Tr \overline{\pi}_{D,\theta}(u)=(-1)^{e(f-1-|S_u|)}q^{n|S_u|/2}\chi_2(\overline{u})^{e-1}\theta(u)$.
                \item If $m$ is even, then $\Tr \overline{\pi}_{D,\theta}(u)=q^{(e-1)f/2}\theta(u)$.
            \end{itemize}
        \end{itemize}
    \end{enumerate}
\end{thm}

\section{\'Etale cohomology for $m>1$}\label{sec:cohomology-m2}
In this section,
we compute the \'etale cohomology of $X_m$
and the traces of the group actions defined in Section~\ref{sec:variety}.
However, since the $m=1$ case involves a more complicated expression,
we treat this separately in the next section.
The results of this section, Propositions~\ref{prop:coh_dim}, \ref{prop:connected}, \ref{prop:trace} and \ref{prop:forget-supp}
are true even if $m=1$, but the proofs of this case are given in the next section.

We write $\LL_\lambda$ for the Artin-Schreier sheaf on $\G_a$ over $\F_{q^f}$
associated to the character $\psi_\lambda$.
For two polynomials $h,h'\in \F_{q^f}[x_1^{q^{-\infty}},\ldots,x_r^{q^{-\infty}}]$,
we write $h\doteq h'$ if there is an isomorphism $h^*\LL_\lambda\cong h'^*\LL_\lambda$ over $\F_{q^f}$
as local systems over $\G_a^r$.
All cohomology groups are taken after base change to $\overline{\F}_q$.

We first rewrite $g_m$ in a more convenient form.
The polynomial $g_m$ can be viewed as a morphism $U:=\G_a^{e(f-1)}\to \G_a$.
We write $U=V_2\times \cdots \times V_f$, where $V_r=\G_a^e$ is the space of $y_r=(y_{1,r},\ldots ,y_{e,r})$.
In this and the next sections,
we always treat these group schemes and $X_m$ as perfections,
and this treatment does not affect the \'etale topos.
We define two automorphisms $F$ and $T$ on $V_r$ as follows:
\[F y_r=y_r^q=(y_{1,r}^q,\ldots ,y_{e,r}^q),\quad T y_r=(y_{2,r},\ldots ,y_{e,r},y_{1,r}).\]
Clearly they commute with each other.
Let $r^*:=f-r+2$ for $r\in \{2,\ldots,f\}$. 
Then we can rewrite $g_m$ as
\begin{align}\label{eqn:g}
    g_m=
    \begin{dcases}
        \sum_{r=2}^f {}^t\!(T^{-s}F^{r-1}y_r)(F^f-T)y_{r^*} & \text{if $m\geq 3$ is odd},\\
        \sum_{r=2}^f {}^t\!(T^s F^{r^*-1}y_r)(F^f-1) y_{r^*}& \text{if $m$ is even}.
    \end{dcases}
\end{align}

We first work with the case where $m\geq 3$ is odd.
We set $A_r:=T^{-s}F^{r-1}$ and $B:=F^f-T$.
Thus $g_m=\sum_{r=2}^f {}^t\!(A_r y_r)(B y_{r^*})$.
Now, for a given endomorphism $D$ of $\G_a^e$,
we choose an ``adjoint'' endomorphism $D^\dagger$ of $\G_a^e$ satisfying
\begin{align*}
    ({}^t\!(Dx)y)^*\LL_\lambda\cong ({}^t\!x(D^\dagger y))^*\LL_\lambda
\end{align*}
as local systems on $\G_a^{2e}$.
Suppose $\lambda \neq 0$.
Since $\Tr(\lambda x^{q^a}y)=\Tr(\lambda^{q^{-a}}xy^{q^{-a}})$,
we can choose 
\begin{align*}
    (F^a)^\dagger=\alpha_a F^{-a},
\end{align*}
where $\alpha_a:=\lambda^{q^{-a}}/\lambda$.
Obviously $T^\dagger=T^{-1}$, so we can choose
\begin{align*}
    A_r^\dagger=\alpha_{r-1}F^{-r+1}T^s,\quad B^\dagger=F^{-f}-T^{-1}=-T^{-1}F^{-f}B.
\end{align*}
Suppose $f$ is odd and write $f=2\mu+1$.
Then we have
\begin{align*}
    g_m(y)&=\sum_{r=2}^f {}^t\!(A_r y_r)(B y_{r^*})\\
    &\doteq \sum_{r=2}^f {}^t\! y_r (A_r^\dagger B y_{r^*})\\
    &\doteq \sum_{r=2}^{\mu+1} {}^t\! y_r (A_r^\dagger B y_{r^*}) + \sum_{r=\mu+2}^f {}^t\! y_{r^*}  (B^\dagger A_r y_r)\\
    &\doteq \sum_{r=2}^{\mu+1} {}^t\! y_r (A_r^\dagger B+B^\dagger A_{r^*})y_{r^*}.
\end{align*}
We put $C_r:=A_r^\dagger B+B^\dagger A_{r^*}$, so that
$g_m\doteq \sum_{r=2}^{\mu+1} {}^t\! y_r C_r y_{r^*}$.
We note that $C_r^\dagger =C_{r^*}$,
and we have
\begin{align*}
    C_r&=\alpha_{r-1}F^{-r+1}T^s B -T^{-s-1}F^{-r+1}B\\
    &=T^{-s-1}(\alpha_{r-1}T^m-1)F^{-r+1}B.
\end{align*}

Next we suppose $f=2\mu$ is even. Then
\begin{align*}
    g_m(y)&\doteq \sum_{r=2}^f {}^t\!(A_r y_r)(B y_{r^*})\\
    &\doteq \sum_{r=2}^\mu {}^t\! y_r (A_r^\dagger B y_{r^*}) + \sum_{r=\mu+2}^f {}^t\! y_{r^*}  (B^\dagger A_r y_r) + {}^t\!(A_{\mu+1} y_{\mu+1})(B y_{\mu+1})\\
    &\doteq \sum_{r=2}^\mu {}^t\! y_r C_r y_{r^*} + {}^t\!(A_{\mu+1} y_{\mu+1})(B y_{\mu+1}).
\end{align*}

For the case where $m$ is even, we put $A_r',B',C_r'$ in the same way:
\[A_r':=T^s F^{r^*-1},\quad B':=F^f-1,\quad C_r':=(A_r')^\dagger B'+(B')^\dagger A_{r^*}'.\]
We have
\begin{align*}
    C_r'&=\alpha_{r^*-1}F^{-r^*+1}T^{-s}B'-T^s F^{-r^*+1}B'\\
    &=T^s(\alpha_{r^*-1}T^{-m}-1) F^{-r^*+1}B'.
\end{align*}

Before computing the dimension of the compactly supported \'etale cohomology spaces
$H^i_c(\A^{e-1},f_m^*\LL_\lambda)$ and $H^i_c(\A^{e(f-1)},g_m^*\LL_\lambda)$,
we explain the result from \cite{takeuchi2023quadraticelladicsheafheisenberg},
which provides a method for computing the dimensions of the compactly supported \'etale cohomology
of quadratic sheaves.
Let $U$ be a perfect unipotent commutative group scheme over a finite field extension $k/\F_p$
and $\QQ$ be a quadratic $\overline{\Q}_l$-sheaf on $U$.
We write $\BB_\QQ$ for the corresponding bimultiplicative sheaf on $U\times U$, that is,
$\BB_\QQ=m^*\QQ\otimes \pr_1^*\QQ^{-1}\otimes \pr_2^*\QQ^{-1}$,
where $m\colon U\times U\to U$ is the multiplication
and $\pr_i\colon U\times U\to U$ is the $i$-th projection.
Let $\ell_\QQ$ be the induced duality morphism $U\to U^*$.
In \cite{takeuchi2023quadraticelladicsheafheisenberg}, Takeuchi proves the following theorem.

\begin{thm}[{{\cite[Theorem~4.2 and Corollary~4.6]{takeuchi2023quadraticelladicsheafheisenberg}}}]\label{thm:Takeuchi}
    Let $U$, $\QQ$, $\BB_\QQ$ and $\ell_\QQ$ be as above.
    \begin{enumerate}
        \item If $\ell_\QQ$ is finite \'etale of degree $p^{2r}$,
        then $H^i_c(U_{\overline{k}},\QQ)$ is nonzero if and only if $i=\dim U$,
        and the dimension is $p^r$.
        \item Let $K^\circ$ be the neutral connected component of $\Ker \ell_\QQ$, which is perfect unipotent.
        If $\QQ|_{K^\circ}$ is nontrivial, then $H^i_c(U_{\overline{k}},\QQ)=0$ for all $i$.
        If $\QQ|_{K^\circ}$ is trivial, then $H^i_c(U_{\overline{k}},\QQ)$ is nonzero if and only if $i=\dim U+\dim K^\circ$,
        and its dimension equals the square root of the cardinality of $\pi_0(\Ker \ell_\QQ)$.
    \end{enumerate}
\end{thm}
\begin{rem}
    We remark that the original statement of Corollary~4.6 in \cite{takeuchi2023quadraticelladicsheafheisenberg}
    appears to contain a typo:
    in the case where $\QQ|_{K^\circ}$ is trivial, it is claimed that there is a descent $\overline{\QQ}$ of $\QQ$ to $U/K^\circ$,
    and the computation of $R\Gamma_c(U_{\overline{k}},\QQ)$ is reduced to that of $R\Gamma_c((U/K^\circ)_{\overline{k}},\overline{\QQ})$,
    concluding that $H^i_c(U_{\overline{k}},\QQ)$ vanishes unless $i=\dim (U/K^\circ)$.
    However, we should consider the contribution of $R\Gamma_c(K^\circ_{\overline{k}}, \QQ|_{K^\circ})$,
    which causes $2\dim K^\circ$-shift on degree.
    Therefore, the correct degree on which $H^i_c(U_{\overline{k}},\QQ)$ is concentrated turns out to be $\dim (U/K^\circ)+2\dim K^\circ=\dim U+\dim K^\circ$.
    Also the dimension should be
    $\sqrt{\#\pi_0(\Ker\ell_{\mathcal Q})}$ rather than
    $\#\pi_0(\Ker\ell_{\mathcal Q})$.
\end{rem}

\begin{prop} \label{prop:coh_dim}
    Let $\lambda\in \F_{q^f}$.
    \begin{enumerate}
        \item If $\lambda\neq 0$, we have $H^i_c(\A^{e-1},f_m^*\LL_\lambda)\neq 0$ if and only if $i=e-1$, and
        \begin{align*}
            \dim H^{e-1}_c(\A^{e-1},f_m^*\LL_\lambda)=
            \begin{cases}
                1 & \text{if $m$ is odd},\\
                q^{f(e-1)} & \text{if $m$ is even}.
            \end{cases}
        \end{align*}
        If $\lambda=0$, we have $H^i_c(\A^{e-1},f_m^*\LL_0)\neq 0$ if and only if $i=2(e-1)$, and
        \begin{align*}
            \dim H^{2(e-1)}_c(\A^{e-1},f_m^*\LL_0)=1.
        \end{align*}
        \item Suppose $f>1$.
        If $\lambda$ satisfies Condition~\ref{cond}, we have $H^i_c(\A^{e(f-1)},g_m^*\LL_\lambda)\neq 0$ if and only if $i=e(f-1)$, and
        \begin{align*}
            \dim H^{e(f-1)}_c(\A^{e(f-1)},g_m^*\LL_\lambda)=q^{ef(f-1)/2}.
        \end{align*}
        Otherwise, we have $H^{e(f-1)}_c(\A^{e(f-1)},g_m^*\LL_\lambda)=0$.
        Moreover, the condition for 
        \[H^{2e(f-1)}_c(\A^{e(f-1)},g_m^*\LL_\lambda)\neq 0\]
        is $\lambda=0$ when $e\geq 2$ and $\lambda\in \F_q$ when $e=1$.
        The dimension of this space is one.
        \item We have
        \begin{align*}
            H^{n-1}_c(X_m,\overline{\Q}_l)[\psi_\lambda]\cong H^{e-1}_c(\A^{e-1},f_m^*\LL_\lambda)\otimes H^{e(f-1)}_c(\A^{e(f-1)},g_m^*\LL_\lambda),
        \end{align*}
        which is nonzero if and only if $\lambda$ satisfies Condition~\ref{cond}.
    \end{enumerate}
\end{prop}

\begin{prf}
    (i) follows from \cite[Propositions~4.11 and 4.12]{MR4126896}.
    We prove (ii) in the case where $m>1$.
    We apply \cite[Theorem~4.2 and Corollary~4.6]{takeuchi2023quadraticelladicsheafheisenberg}
    taking the perfect unipotent commutative group scheme $U=(\G_a^{\perf})^{e(f-1)}$
    and the quadratic sheaf $\QQ=g_m^*\LL_\lambda$.
    We set
    \[b_m(x,y):=g_m(x+y)-g_m(x)-g_m(y),\]
    so that the corresponding bimultiplicative sheaf on $U\times U$ is $\BB_\QQ=b_m^*\LL_\lambda$.
    Let $\ell\colon U\to U^*$ denote the group morphism induced by $\BB_\QQ$.
    By the above computation, if $m$ is odd, we have
    \begin{align*}
        b_m(x,y)&\doteq \sum_{r=2}^f ({}^t\! x_r A_r^\dagger B y_{r^*} + {}^t\! y_r A_r^\dagger B x_{r^*})\\
        &\doteq \sum_{r=2}^f ({}^t\! x_r A_r^\dagger B y_{r^*} + {}^t\! x_r B^\dagger A_{r^*}  y_{r^*})\\
        &= \sum_{r=2}^f {}^t\! x_r C_r y_{r^*}.
    \end{align*}
    Therefore, if we combine $\ell$ with the identification
    \[U\xrightarrow{\sim} U^*; y\mapsto ({}^t\! xy)^*\LL_\lambda,\]
    the resulting morphism $U\to U$ is written as
    \[y=(y_2,\ldots ,y_f)\mapsto (C_2y_f,\ldots ,C_fy_2).\]
    Now, since $T^e=1$ and $m$ is coprime to $e$, 
    Condition~\ref{cond}, which says that $\alpha_{r-1}$ is not an $e$-th root of unity,
    implies that $\alpha_{r-1}T^m-1$ is invertible.
    Also $T$ and $F$ are automorphisms of $U$,
    so $\Ker C_r=\Ker B$.
    Since
    \[y_r\in \Ker B \Longleftrightarrow \forall i\in \{1,\ldots,e\},\, y_{i,r}^{q^f}=y_{i+1,r},\]
    we have $\Ker B\cong \F_{q^n}$.
    Therefore, we see that $\ell$ is finite \'etale of degree $q^{n(f-1)}=q^{ef(f-1)}$,
    which implies the statement for $\dim H^{e(f-1)}_c$ by \cite[Theorem~4.2]{takeuchi2023quadraticelladicsheafheisenberg}.

    If $m$ is even, the same computation gives $b_m(x,y)\doteq \sum_{r=2}^f {}^t\! x_r C_r' y_{r^*}$,
    and $\ell$ can be expressed as $(y_2,\ldots ,y_f)\mapsto (C_2' y_f,\ldots ,C_f' y_2)$.
    Again $\Ker C_r'=\Ker B'$, and this is isomorphic to $\F_{\!q^f}^{\,e}$.
    Thus $\ell$ is finite \'etale of degree $q^{ef(f-1)}$, the same as when $m$ is odd.

    When $\lambda$ does not satisfy Condition~\ref{cond},
    it is easy to see that $\ell_\QQ$ has the kernel of positive dimension,
    thus Theorem~\ref{thm:Takeuchi} (ii) implies that $H^{e(f-1)}_c(\A^{e(f-1)},g_m^*\LL_\lambda)$ vanishes.
    In addition, the non-vanishing of $H^{2e(f-1)}_c(\A^{e(f-1)},g_m^*\LL_\lambda)$ occurs exactly when the sheaf $g_m^*\LL_\lambda$ is trivial,
    in which case the dimension is obviously one.
    This occurs if and only if $\lambda=0$ when $e\geq 2$ and $\lambda\in \F_q$ when $e=1$.

    For (iii), we have
    \[R\Gamma_c(X_m,\overline{\Q}_l)[\psi_\lambda]\cong R\Gamma_c(\A^{e-1},f_m^*\LL_\lambda)\otimes R\Gamma_c(\A^{e(f-1)},g_m^*\LL_\lambda)\]
    by K\"unneth's formula.
    However, we have $H^i_c(\A^{e-1},f_m^*\LL_\lambda)=0$ if $i<e-1$
    and $H^i_c(\A^{e(f-1)},g_m^*\LL_\lambda)=0$ if $i<e(f-1)$
    by affine vanishing.
    Therefore the statement follows by (i), (ii) and noting that $(e-1)+e(f-1)=n-1$.
\end{prf}

We have two actions of $U_L^m$ on $X_m$: the action of $U_L^m\subset \J_{\GL}^m$ and that of $U_L^m\subset \J_D^m$,
both of which factor through $U_L^m/U_L^{m+1}\cong \F_{q^f}$.
The next proposition is on these actions, but one is the inverse of the other,
so it does not matter which action we choose.
\begin{prop}\label{prop:connected}
    The quotient group $\Coker(U_L^{m,1}\to \F_{q^f})$ acts on the component set of $X_m$ simply transitively.
    In particular, $X_m$ is connected if $e\geq 2$,
    and $X_m$ has $q$ connected components, pairwise isomorphic, if $e=1$.
\end{prop}
\begin{prf}
    By K\"unneth's formula and affine vanishing, we have
    \[H^{2(n-1)}_c(X_m,\overline{\Q}_l)[\psi_\lambda]\cong H^{2(e-1)}_c(\A^{e-1},f_m^*\LL_\lambda)\otimes H^{2e(f-1)}_c(\A^{e(f-1)},g_m^*\LL_\lambda).\]
    If $e\geq 2$, the right hand side remains only when $\lambda=0$ and has dimension one.
    Since $H^{2(n-1)}_c(X_m)\cong \bigoplus_{\lambda\in \F_{q^f}} H^{2(n-1)}_c(X_m)[\psi_\lambda]$,
    this means that $H^{2(n-1)}_c(X_m)$ is one-dimensional, which implies that $X_m$ is geometrically connected.
    In addition, since $e$ does not divide $m$, we see that $N_{L/K}(U_L^m)=N_{L/K}(U_L^{m+1})$,
    which implies that the composition $U_L^{m,1}\hookrightarrow U_L^m\to U_L^m/U_L^{m+1}\cong \F_{q^f}$ is surjective.
    Then the proposition follows in this case.

    If $e=1$, then the right hand side remains when $\lambda\in \F_q$
    (note that in this case $\A^{e-1}$ is one point, so the first term is always trivial one-dimensional space),
    thus $H^{2(n-1)}_c(X_m)$ is $q$-dimensional.
    Therefore $X_m$ has $q$ connected components (cf.\ \cite[Remark~4.5.2]{MR3402698}).
    Moreover, Boyarchenko and Weinstein constructed in \cite{MR3402698} an algebraic group $\Ub_m$ over $\F_q$
    equipped with the norm map $\Nb\colon \Ub_m\to \G_a$
    satisfying the following properties: 
    \begin{itemize}
        \item $\Ub_m$ is an algebraic group over $\F_q$ such that
        \[\Ub_m(R)=\{1+r_1e_1+\cdots +r_n e_n\mid r_i\in R\}\]
        for any $\F_q$-algebra, with a certain group law, and
        $\Ub_m(\F_{q^n})$ is identified with the group $S_G^m$ for even $m$ (or, equivalently, $S_D^m$ for odd $m$).
        The center $Z(\Ub_m)$ of $\Ub_m$ is $\{1+r_n e_n\}\cong \G_a$, and we write $\pr_n\colon \Ub_m\to \G_a$ for the projection to $n$-th component.
        \item $X_m$ is identified with the closed subvariety of $\Ub_{m,\overline{\F}_q}$ defined by $\pr_n\circ L=0$,
        where $L\colon \Ub_m\to \Ub_m$ is $g\mapsto \Fr_{q^n}(g)g^{-1}$.
        The right action of $S_G^m$ for even $m$ (resp.\ left action of $S_D^m$ for odd $m$)
        is identified with the right (resp.\ left) multiplication of $\Ub_m$.
        \item $\Nb(gh)=\Nb(g)+\Nb(h)$ for $g\in \Ub_m$ and $h\in \Ub(\F_{q^n})$.
        \item $\Nb$ restricts to the trace map $\Tr_{\F_{q^n}/\F_q}$ on the center $Z(\Ub_m(\F_{q^n}))\cong \F_{q^n}$. 
        \item $\pr_n(L(g))=\Nb(g)^q-\Nb(g)$ for any $g\in \Ub_m$.
    \end{itemize}
    Then the last property shows that $X_m$ is the disjoint union of closed subvarieties $X_{m,b}$ of the form $\Nb(g)=b$ where $b$ runs over $\F_q$.
    The argument above on the dimension of $H^{2(n-1)}_c$ implies that each $X_{m,b}$ is a connected component of $X_m$.
    For $z\in Z(\Ub_m(\F_{q^n}))$ and $x\in X_{m,b}$, we have
    $\Nb(xz)=\Nb(x)+\Nb(z)=b+\Nb(z)$, thus $z$ maps $X_{m,b}$ to $X_{m,b+\Nb(z)}$.
    Finally, since $\Nb(z)=\Tr(\pr_n(z))$, the proposition follows
    by the fact that the image of $U_L^{m,1}\hookrightarrow U_L^m\to U_L^m/U_L^{m+1}\cong \F_{q^f}$ is the subspace of elements with trace zero.
\end{prf}

\begin{lem}\label{lem:trisum}
    Suppose that $\lambda\in \F_{q^f}$ satisfies Condition~\ref{cond}. 
    Let $S(g_m)_\lambda$ denote the triangular sum
    \[\sum_{y\in U(\F_{q^f})} \psi_\lambda(g_m(y)).\]
    Then we have 
    \begin{align*}
        S(g_m)_\lambda =
        \begin{cases}
            (-1)^{(e-1)(f-1)}q^{(e+1)f(f-1)/2} & \text{ if $m$ is odd,}\\
            q^{ef(f-1)} & \text{ if $m$ is even}.
        \end{cases}
    \end{align*}
\end{lem}

\begin{prf}
    We first treat the case where $mf$ is odd.
    Then we have
    \begin{align*}
        S(g_m)_\lambda
        &=\sum_{y\in U(\F_{\!q^f})}\prod_{r=2}^{\mu+1} \psi_\lambda({}^t\! y_r C_r y_{r^*})\\
        &=\prod_{r=2}^{\mu+1} \sum_{y_{r^*}\in \F_{\!q^f}^{\,e}} \sum_{y_{r}\in \F_{\!q^f}^{\,e}} \psi_\lambda ({}^t\! y_r C_r y_{r^*}).
    \end{align*}
    For a fixed $y_{r^*}\in \F_{q^f}^e$, if $C_r y_{r^*}=(z_1,\ldots ,z_e)\neq 0$, there exists an index $i_0$ such that $z_{i_0}\neq 0$, and
    \begin{align*}
        \sum_{y_{r}\in \F_{q^f}^e} \psi_\lambda ({}^t\! y_r C_r y_{r^*})
        &=\prod_{i=1}^e \sum_{y_{r,i}\in \F_{q^f}}\psi_\lambda (y_{r,i} z_i)
    \end{align*}
    is zero, since
    \begin{align*}
        \sum_{y_{r,{i_0}}\in \F_{q^f}}\psi_\lambda (y_{r,{i_0}} z_{i_0})
        =\sum_{u\in \F_{q^f}}\psi_\lambda (u)=0.
    \end{align*}
    Therefore, we get
    \begin{align*}
        S(g_m)_\lambda
        &= \prod_{r=2}^{\mu+1} \sum_{y_{r^*}\in (\Ker C_r)(\F_{q^f})} \sum_{y_{r}\in \F_{q^f}^e} \psi_\lambda (0)\\
        &=(q^f\cdot q^n)^{\mu}
        =q^{(e+1)f(f-1)/2}.
    \end{align*}
    
    Next we consider the case where $m$ is odd and $f$ is even.
    We have 
    \begin{align*}
        S(g_m)_\lambda
        &= \sum_{y\in U(\F_{q^f})}\left(\prod_{r=2}^{\mu} \psi_\lambda({}^t\! y_r C_r y_{r^*})\cdot \psi_\lambda ({}^t\!(A_{\mu+1} y_{\mu+1})(B y_{\mu+1}))\right)\\
        &= \left(\prod_{r=2}^{\mu} \sum_{y_{r^*}\in \F_{q^f}^e} \sum_{y_{r}\in \F_{q^f}^e} \psi_\lambda ({}^t\! y_r C_r y_{r^*})\right)
        \left(\sum_{u\in \F_{q^f}^e} \psi_\lambda({}^t\!(A_{\mu+1} u)(B u))\right),
    \end{align*}
    and the first bracket can be computed in the same way as the previous case:
    \begin{align*}
        \prod_{r=2}^{\mu} \sum_{y_{r^*}\in \F_{q^f}^e} \sum_{y_{r}\in \F_{q^f}^e} \psi_\lambda ({}^t\! y_r C_r y_{r^*})
        =(q^f\cdot q^n)^{\mu-1}.
    \end{align*}
    For the computation of the second bracket, we set $Q=q^\mu$ so that $q^f=Q^2$,
    and for $u\in \F_{Q^2}$, we write $F^\mu u=\overline{u}$, the conjugate over $\F_Q$.
    Then we have
    \begin{align*}
        {}^t\!(A_{\mu+1} u)(B u)
        ={}^t\!(T^{-s} \overline{u})(1-T)u
        ={}^t\! \overline{u} (T^s(1-T))u.
    \end{align*}
    Thus, if we write $M_\lambda=\lambda T^s(1-T)$
    and $H_\lambda=M_\lambda + {}^t\! \overline{M}_\lambda$,
    the sum in the second bracket turns to 
    \begin{align*}
        \sum_{u\in \F_{Q^2}^e} \psi(\Tr_{\F_{Q^2}/\F_q}({}^t\! \overline{u} M_\lambda u))
        = \sum_{u\in \F_{Q^2}^e} \psi^Q({}^t\! \overline{u} H_\lambda u),
    \end{align*}
    where $\psi^Q=\psi\circ \Tr_{\F_Q/\F_q}$.
    Since $H_\lambda$ is Hermitian, it can be diagonalized by some invertible matrix,
    so a new coordinate $v$ allows us to rewrite the sum as
    \begin{align*}
        \sum_{v\in \F_{Q^2}^e} \psi^Q(c_1 v_1^{1+Q}+\cdots +c_r v_r^{1+Q})=(Q^2)^{e-r}\prod_{i=1}^r \sum_{v_i\in \F_{Q^2}}\psi^Q(c_i v_i^{1+Q}),
    \end{align*}
    where $c_i\in \F_Q^\times$ and $r$ is the rank of the matrix $H_\lambda$.
    We have
    \begin{align*}
        H_\lambda &= \lambda T^s(1-T)+\overline{\lambda} T^{-s}(1-T^{-1})\\
        &=\lambda T^{-s-1} (\lambda^{-1}\overline{\lambda}-T^m)(T-1),
    \end{align*}
    so that $\Ker H_\lambda=\Ker (T-1)$ is $1$-dimensional, which implies $r=e-1$.
    Moreover, we have
    \begin{align*}
        \sum_{v_i\in \F_{Q^2}}\psi^Q(c_i v_i^{1+Q})
        &=\sum_{v_i\in \F_{Q^2}}\psi^Q(c_i \Nr_{\F_{Q^2}/\F_Q}(v_i))\\
        &=1+(Q+1)\sum_{w\in \F_Q^\times}\psi^Q(w)\\
        &=1-(Q+1)
        =-Q.
    \end{align*}
    Therefore, we conclude that
    \begin{align*}
        S(g_m)_\lambda
        &=(q^f\cdot q^n)^{\mu-1}\cdot Q^2(-Q)^{e-1}\\
        &=(-1)^{e-1}q^{(e+1)f(f-1)/2}.
    \end{align*}

    Finally, if $m$ is even, $(F^f-1)y_{r^*}=0$ for $y_{r^*}\in \F_{q^f}^e$,
    so that $\psi_\lambda(g_m(y))=1$ for every $y\in U(\F_{q^f})$.
    Thus the statement is trivial.
\end{prf}

\begin{prop}\label{prop:trace}
    Suppose $\lambda\in \F_{q^f}$ satisfies Condition~\ref{cond}.
    The traces of the following elements of $J_{m,s}^1$ on the cohomology
    are summarized in the tables below, where $u\in U_L$,
    and $\mathfrak{g}_\lambda$ and $S_u$ are as defined in Theorem~\ref{thm:characterize}.
    If $m$ is odd, we have the following table.
    \begin{center}
        {
        \renewcommand{\arraystretch}{1.4}
        \begin{tabular}{|c||c|c|}\hline
             & $H^{e-1}_c(\A^{e-1},f_m^*\LL_\lambda)$ & $H^{e(f-1)}_c(\A^{e(f-1)},g_m^*\LL_\lambda)$\\ \hhline{|=#=|=|}
            \begin{tabular}{c}
                $(1,\varphi_L^{-1},\Art_L(\varphi_L))$ \\ $=\Frob_{q^f}$
            \end{tabular}
            & $\mathfrak{g}_\lambda$ & $(-1)^{f-1}q^{(e+1)f(f-1)/2}$\\ \hline
            $(\varphi_L^{\pm 1},\varphi_L^{\pm 1},1)$ & $(-1)^{e-1}$ & $(-1)^{(e-1)(f-1)}q^{f(f-1)/2}$\\ \hline
            $(u,u,1)$ & $1$ & $(-1)^{e(f-1-|S_u|)}q^{n|S_u|/2}$\\ \hline
            $(1,u^{-1},\Art_L(u))$ & $\chi_2(\overline{u})^{e-1}$ & $(-1)^{e(f-1-|S_u|)}q^{n|S_u|/2}$\\\hline
        \end{tabular}
        }
    \end{center}
    
    If $m$ is even, we have the following table.
    \begin{center}
        {
        \renewcommand{\arraystretch}{1.4}
        \begin{tabular}{|c||c|c|}\hline
             & $H^{e-1}_c(\A^{e-1},f_m^*\LL_\lambda)$ & $H^{e(f-1)}_c(\A^{e(f-1)},g_m^*\LL_\lambda)$\\ \hhline{|=#=|=|}
            \begin{tabular}{c}
                $(1,\varphi_L^{-1},\Art_L(\varphi_L))$ \\ $=\Frob_{q^f}$
            \end{tabular} 
            & $q^{f(e-1)}$ & $(-1)^{e(f-1)}q^{ef(f-1)}$\\ \hline
            $(\varphi_L^{\pm 1},\varphi_L^{\pm 1},1)$ & $1$ & $(-1)^{(e-1)(f-1)}q^{f(f-1)/2}$\\ \hline
            $(u,u,1)$ & $q^{f(e-1)}$ & $(-1)^{e(f-1-|S_u|)}q^{n|S_u|/2}$\\ \hline
            $(1,u^{-1},\Art_L(u))$ & $q^{f(e-1)}$ & $q^{ef(f-1)/2}$\\\hline
        \end{tabular}
        }
    \end{center}
\end{prop}
\begin{prf}
    The statement for $H^{e-1}_c(\A^{e-1},f_m^*\LL_\lambda)$
    follows by \cite[Propositions~4.11 and 4.12]{MR4126896}.
    We consider the case of $H^{e(f-1)}_c(X_{g_m})_\lambda:=H^{e(f-1)}_c(\A^{e(f-1)},g_m^*\LL_\lambda)$.
    By Proposition~\ref{prop:coh_dim}, we have
    \begin{align*}
        \Tr(\Frob_{q^f}; H^{e(f-1)}_c(X_{g_m})_\lambda)&=(-1)^{e(f-1)}\sum_i (-1)^i\Tr(\Frob_{q^f}; H^i_c(X_{g_m})_\lambda)\\
        &=(-1)^{e(f-1)}S(g_m)_\lambda,
    \end{align*}
    where the last equality is by the Grothendieck-Lefschetz trace formula.
    Thus we get the desired trace by Lemma~\ref{lem:trisum}.

    Next let $\delta$ denote one of the automorphisms
    $(\varphi_L,\varphi_L,1)$, $(u,u,1)$ or $(1,u^{-1},\Art_L(u))$
    for $u\in U_L$,
    and we compute their traces.
    To this end, we use the fixed point formula from \cite[Theorem 3.2]{MR393266}:
    \begin{quote}
        Let $X$ be a scheme, separated and of finite type over $\overline{\F}_q$
        and let $\sigma\colon X\to X$ be the automorphism of finite order of $X$.
        We decompose $\sigma$ as $\sigma=s\cdot u$ where $s$ and $u$ are powers of $\sigma$ respectively
        of order prime to $p$ and a power of $p$.
        Then we have 
        \[\sum_i (-1)^i \Tr(\sigma; H_c^i(X,\overline{\Q}_l))=\sum_i (-1)^i \Tr(u; H_c^i(X^s,\overline{\Q}_l)),\]
        where $X^s$ denotes the fixed point variety with respect to the action of $s$.
    \end{quote}
    We put
    \[e_{\psi_\lambda}=\frac{1}{q^f}\sum_{x\in \F_{q^f}}\psi_\lambda (x)^{-1}x,\]
    which is the idempotent element of the group algebra $\overline{\Q}_l[\F_{q^f}]$
    corresponding to $\psi_\lambda$.
    Then we have
    \begin{align*}
        \sum_i (-1)^i\Tr(\delta; H_c^i(X_{g_m})_\lambda)
        &=\sum_i (-1)^i\Tr(\delta e_{\psi_\lambda}; H_c^i(X_{g_m}))\\
        &=\frac{1}{q^f}\sum_{x\in \F_{q^f}}\psi_\lambda(x)^{-1}\sum_i (-1)^i \Tr(\delta x; H_c^i(X_{g_m}))\\
        &=\frac{1}{q^f}\sum_{x\in \F_{q^f}}\psi_\lambda(x)^{-1}\sum_i (-1)^i \Tr(x; H_c^i(X_{g_m}^{\delta}))\\
        &=\sum_i (-1)^i \Tr(e_{\psi_\lambda};H_c^i(X_{g_m}^{\delta}))\\
        &=\sum_i (-1)^i \dim H_c^i(X_{g_m}^{\delta})_\lambda.
    \end{align*}
    If $\delta=(\varphi_L^{\pm 1},\varphi_L^{\pm 1},1)$, its action on $X_{g_m}$ is $T^{\mp 1}$,
    so the fixed point variety $X_{g_m}^\delta$ is $z^{q^f}-z=g_m^\delta$,
    where $g_m^\delta$ is a polynomial function on $\G_a^{f-1}=\{(y_2,\ldots ,y_f)\}$
    obtained by replacing $T$ by $1$ in the expression (\ref{eqn:g}),
    and $C_r$ (resp.\ $C_r'$) is replaced by $(\alpha_{r-1}-1)F^{-r+1}(F^f-1)$ (resp.\ $(\alpha_{r^*-1}-1)F^{-r^*+1}(F^f-1)$).
    Then the same argument in the proof of Proposition~\ref{prop:coh_dim} using \cite[Theorem~4.2]{takeuchi2023quadraticelladicsheafheisenberg}
    concludes that
    \begin{align*}
        \sum_i (-1)^i \dim H_c^i(X_{g_m}^{\delta})_\lambda
        &=(-1)^{f-1}q^{f(f-1)/2}.
    \end{align*}

    If $\delta=(u,u,1)$, its action on $X_{g_m}$ is $(y_{ij},z)\mapsto (\overline{u}^{1-q^{1-j}}y_{ij},z)$
    or $(y_{ij},z)\mapsto (\overline{u}^{q^{j-1}-1}y_{ij},z)$,
    depending on the parity of $m$.
    Since $\overline{u}^{1-q^{1-j}}=1$ if and only if $\overline{u}\in \F_{q^{j-1}}$,
    the fixed point variety $X_{g_m}^\delta$ is the intersection of $X_{g_m}$ and
    $y_{ij}=0$ for $j\not\in S_u$.
    Noting that $S_u$ is closed under $j\mapsto j^*$,
    the same argument as above gives
    \begin{align*}
        \sum_i (-1)^i \dim H_c^i(X_{g_m}^{\delta})_\lambda
        =(-1)^{e|S_u|}q^{n|S_u|/2}.
    \end{align*}

    Finally, if $\delta=(1,u^{-1},\Art_L(u))$,
    its action on $X_{g_m}$ is the same as $(u^{-1},u^{-1},1)$ if $m$ is odd and trivial if $m$ is even,
    so the trace is as desired.
\end{prf}

\begin{prop}\label{prop:forget-supp}
    The natural morphism $H^{n-1}_c(X_m,\overline{\Q}_l)\to H^{n-1}(X_m,\overline{\Q}_l)$ is an isomorphism.
\end{prop}
\begin{prf}
    Since we have the decomposition
    \[H^{n-1}_c(X_m,\overline{\Q}_l)\cong \bigoplus_{\lambda\in \F_{q^f}} H^{n-1}_c(X_m,\overline{\Q}_l)[\psi_\lambda]\]
    and the similar one for $H^{n-1}(X_m,\overline{\Q}_l)$,
    it suffices to prove that $H^{n-1}_c(X_m,\overline{\Q}_l)[\psi_\lambda]\to H^{n-1}(X_m,\overline{\Q}_l)[\psi_\lambda]$
    is an isomorphism for each $\lambda \in \F_{q^f}$.
    By the K\"unneth formula and affine vanishing, we have
    \[H^{n-1}_c(X_m,\overline{\Q}_l)[\psi_\lambda]\cong H^{e-1}_c(\A^{e-1},f_m^*\LL_\lambda)\otimes H^{e(f-1)}_c(\A^{e(f-1)},g_m^*\LL_\lambda),\]
    and the same for $H^{n-1}(X_m,\overline{\Q}_l)[\psi_\lambda]$.

    By Proposition~\ref{prop:coh_dim} (i), we have $H^{e-1}_c(\A^{e-1},f_m^*\LL_\lambda)=0$ if $e>1$ and $\lambda =0$,
    which implies the vanishing $H^{e-1}(\A^{e-1},f_m^*\LL_\lambda)=0$ by the Poincar\'e duality.
    Similarly, by Proposition~\ref{prop:coh_dim} (ii) and Poincar\'e duality,
    we have 
    \[H^{e(f-1)}_c(\A^{e(f-1)},g_m^*\LL_\lambda)=H^{e(f-1)}(\A^{e(f-1)},g_m^*\LL_\lambda)=0\]
    for $\lambda$ not satisfying Condition~\ref{cond}.
    Therefore, we only need to prove that
    \begin{itemize}
        \item $H^{e-1}_c(\A^{e-1},f_m^*\LL_\lambda)\to H^{e-1}(\A^{e-1},f_m^*\LL_\lambda)$ for $\lambda\neq 0$, and
        \item $H^{e(f-1)}_c(\A^{e(f-1)},g_m^*\LL_\lambda) \to H^{e(f-1)}(\A^{e(f-1)},g_m^*\LL_\lambda)$ for $\lambda$ satisfying Condition~\ref{cond}
    \end{itemize}
    are isomorphisms.
    We show these by applying \cite[Proposition~5.13]{takeuchi2023quadraticelladicsheafheisenberg}.
    To this end, it suffices to show that $f_m^*\LL_\lambda$ and $g_m^*\LL_\lambda$ are isogeneous quadratic sheaves,
    that is, each corresponding $\ell_\QQ$ is finite \'etale.
    This is already shown for $g_m^*\LL_\lambda$ under Condition~\ref{cond},
    so we consider $\QQ=f_m^*\LL_\lambda$.
    Since $f_m(x+y)-f_m(x)-f_m(y)$ is biadditive, $f_m^*\LL_\lambda$ is a quadratic sheaf.
    We know that $H^i_c(\A^{e-1},f_m^*\LL_\lambda)$ is concentrated on the middle degree,
    so Theorem~\ref{thm:Takeuchi} implies that $\Ker \ell_\QQ$ is $0$-dimensional,
    thus $f_m^*\LL_\lambda$ is isogeneous.
\end{prf}

\section{\'Etale cohomology for $m=1$}\label{sec:cohomology-m1}

We prove the assertions of Propositions~\ref{prop:coh_dim}, \ref{prop:connected}, \ref{prop:trace}, and~\ref{prop:forget-supp} for $m=1$.
The argument for $m>1$ does not apply directly, because $g_1$
contains terms of length greater than two in the determinant
expansion. We first show that these extra terms become
Artin-Schreier coboundaries after the preceding variables have
been eliminated. We then use the same relative Fourier
calculation as in the quadratic case.

The propositions we are going to prove are obvious when $\lambda=0$,
so we suppose $\lambda\neq 0$ unless otherwise stated.

\subsection{The determinant expansion and adjoint operators}
\label{subsec:m1-adjoint}

Set $k=\F_{q^f}$, $N=e(f-1)$ and $r_0=\lceil f/2 \rceil$,
and we freely use the notations from the previous section.
For $m=1$, the definition of $g_1$ in Definition~\ref{defn:var} (ii)
yields
\begin{equation}\label{eq:m1-compositions}
 g_1=\sum_{i=1}^e
 \sum_{\substack{L\geq2,\ d_1,\ldots,d_t\geq1\\
                  d_1+\cdots+d_t=f}}
 (-1)^t\left(\prod_{\nu=1}^{L-1}
 y_{i,d_\nu+1}^{q^{s_\nu}}\right)(By_{d_t+1})_i,
 \quad s_\nu=\sum_{j=1}^\nu d_j.
\end{equation}
Indeed, every nonidentity permutation contributing a nonzero
term to the defining determinant consists of the single cycle
\[
 1\longmapsto1+s_1\longmapsto\cdots
 \longmapsto1+s_{t-1}\longmapsto1
\]
and fixed points. Its sign is $(-1)^{t-1}$, and the minus sign
in the definition of $g_1$ gives $(-1)^t$.
Conversely, these cycles correspond bijectively to ordered
compositions of $f$.
In particular, the terms with $t=2$
give
\begin{equation}\label{eq:m1-paired-part}
 \sum_{r=2}^f{}^t(F^{r-1}y_r)By_{r^*},
\end{equation}
but the terms with $t\geq 3$ must be retained.

For $2\leq r\leq f$, define
\begin{equation}\label{eq:m1-Cr}
 C_r=(F^{r^*-1})^\dagger B+B^\dagger F^{r-1}.
\end{equation}
As $r^*-1=f-r+1$, we have
\begin{equation}\label{eq:m1-Cr-factor}
 C_r=
 (\alpha_{r^*-1}I-T^{-1})F^{-(r^*-1)}B.
\end{equation}
Since $p\nmid e$, the matrix
$\alpha_{r^*-1}I-T^{-1}$ is invertible precisely when
$\alpha_{r^*-1}\notin\mu_e$, equivalently
$\alpha_{r-1}\notin\mu_e$.
If it is not invertible, its kernel has dimension one.

Define $h_1=g_1$ and, for $1\leq r\leq r_0$,
\begin{equation}\label{eq:m1-hr}
 h_r(y_2,\ldots,y_{r^*-1})
 =g_1(y_2,\ldots,y_{r^*-1},0,\ldots,0).
\end{equation}

\begin{lem}\label{lem:m1-step}
Let $2\leq r\leq r_0$. On any subscheme on which
$By_j=0$ for $2\leq j<r$, one has
\begin{align}
 h_{r-1}={}&h_r+
 {}^t(F^{r-1}y_r)By_{r^*}
 +{}^t(F^{r^*-1}y_{r^*})By_r
 +{}^tP_rBy_{r^*},\label{eq:m1-step}\\
 (P_r)_i={}&
 \sum_{\substack{t\geq2,\ d_1,\ldots,d_t\geq1\\
                  d_1+\cdots+d_t=r-1}}
 (-1)^{t+1}
 \prod_{\nu=1}^t y_{i,d_\nu+1}^{q^{d_1+\cdots+d_\nu}}.
 \label{eq:m1-Pr}
\end{align}
In particular $BP_r=0$ on this subscheme and
\begin{equation}\label{eq:m1-step-adjoint}
 h_{r-1}\equiv_\lambda
 h_r+{}^ty_{r^*}C_ry_r.
\end{equation}
\end{lem}
\begin{prf}
Give $y_j$ weight $j-1$. Every monomial in
\eqref{eq:m1-compositions} has total weight $f$, counting
the final factor $By_{d_L+1}$ by weight $d_L$.
After setting $y_{r^*+1}=\cdots=y_f=0$, terms involving
$y_{r^*}$ arise in exactly two ways.
If the final factor is $By_{r^*}$, the preceding weights
sum to $r-1$; the two-factor contribution is
${}^t(F^{r-1}y_r)By_{r^*}$, and the remaining
contributions are ${}^tP_rBy_{r^*}$.
If $y_{r^*}$ occurs earlier, the last factor must
have weight at least $r-1$, so the total-weight
condition forces the two-factor term
${}^t(F^{r^*-1}y_{r^*})By_r$.

Because $P_r$ is a polynomial over $\mathbb F_q$ in
$y_2,\ldots,y_{r-1}$, the identities $F^fy_j=Ty_j$
imply $F^fP_r=TP_r$, hence $BP_r=0$.
Thus ${}^t\!P_rBy_{r^*}\equiv_\lambda {}^t\!y_{r^*}B^\dagger P_r =0$.
\end{prf}

\subsection{The schemes $K_r,S_r,Z_r,Z_r'$}
\label{subsec:m1-schemes}

Put
\[
 K_r=\Ker C_r,\quad
 S_1=\operatorname{Spec}k,\quad
 S_r=K_2\times\cdots\times K_r,
\]
and denote the projection $S_r\to S_{r-1}$ by $\pi_r$.
For $2\leq r\leq r_0+1$ let
\[
 \begin{aligned}
 Z_r&=S_{r-1}\times V_r\times\cdots\times V_{r^*},\\
 Z_r'&=S_{r-1}\times V_r\times\cdots\times V_{r^*-1},
 \end{aligned}
\]
and let $p_r:Z_r\to Z_r'$ drop $V_{r^*}$.
Thus $Z_2=U$.

The additive morphism $B=F^f-T$ is finite \'etale
and surjective of degree $q^{ef}$.
It follows from \eqref{eq:m1-Cr-factor} that
\begin{equation}\label{eq:m1-Kr}
 \begin{cases}
 K_r=\Ker B\text{ is finite \'etale of degree }q^{ef},
 &\alpha_{r-1}\notin\mu_e,\\
 K_r\text{ is smooth of pure dimension }1
 \text{ with finitely many components},
 &\alpha_{r-1}\in\mu_e.
 \end{cases}
\end{equation}

\begin{lem}\label{lem:m1-relative}
Let $\lambda\ne0$, and assume
$K_j=\Ker B$ for $2\leq j<r$. Put
\[
 Q_r={}^ty_{r^*}C_ry_r:Z_r\longrightarrow\mathbb G_a.
\]
Write $i_r:W_r\hookrightarrow Z_r'$ for the zero
scheme of $C_ry_r$. Then
\begin{equation}\label{eq:m1-rp}
 Rp_{r!}Q_r^*\mathcal L_\lambda
 \cong i_{r*}\overline{\Q}_l(-e)[-2e],
\end{equation}
and consequently
\begin{equation}\label{eq:m1-red}
 R\Gamma_c(Z_r,h_{r-1}^*\mathcal L_\lambda)
 \cong R\Gamma_c(W_r,h_r^*\mathcal L_\lambda)
 (-e)[-2e].
\end{equation}
If moreover $K_r=\Ker B$, then
$W_r=Z_{r+1}$.
These isomorphisms preserve the Weil structures.
\end{lem}
\begin{prf}
Using \cite[Lemma~9.4]{MR2734333},
the isomorphism \eqref{eq:m1-rp} is verified at the stalk level.
The equations $By_j=0$ hold identically on
$S_{r-1}$, so Lemma~\ref{lem:m1-step} applies over
the whole base.
The projection formula and
\eqref{eq:m1-step-adjoint} give \eqref{eq:m1-red}.
\end{prf}

\subsection{The cohomology of $g_1^*\mathcal L_\lambda$}
\label{subsec:m1-main-cohomology}

\begin{lem}\label{lem:m1-middle}
Let $f>1$ and $\lambda\in k$.
If $\lambda$ satisfies Condition~\ref{cond}, then
\begin{equation}\label{eq:m1-middle}
 H_c^i(U_{\bar k},g_1^*\mathcal L_\lambda)=0
 \quad(i\ne N),\quad
 \dim H_c^N(U_{\bar k},g_1^*\mathcal L_\lambda)
 =q^{ef(f-1)/2}.
\end{equation}
Otherwise $H_c^N(U_{\bar k},g_1^*\mathcal L_\lambda)=0$.
\end{lem}
\begin{prf}
Assume first that Condition~\ref{cond} holds.
Then $K_r=\Ker B$ for all $2\leq r\leq r_0$.
Apply Lemma~\ref{lem:m1-relative} successively.
If $f=2\mu+1$, then $r_0=\mu+1$ and the
terminal scheme is $S_{\mu+1}=(\Ker B)^\mu$.
Each term of $h_{\mu+1}$ has final factor $By_j$
for some $j\leq\mu+1$, so $h_{\mu+1}=0$
on this scheme. Thus
\begin{equation}\label{eqn:f-odd}
 R\Gamma_c(U_{\bar k},g_1^*\mathcal L_\lambda)
 \cong
 H^0(S_{\mu+1,\bar k},\overline{\Q}_l)
 (-e\mu)[-2e\mu].
\end{equation}
It is concentrated in degree $2e\mu=N$,
of dimension $q^{ef\mu}$.

If $f=2\mu$, the paired eliminations stop
at $r=\mu$.
On $S_\mu\times V_{\mu+1}$,
we have
\[
 h_\mu\equiv_\lambda Q_0(y),\quad
 Q_0(y)={}^t(F^\mu y)(F^{2\mu}-T)y.
\]
Indeed, all terms whose final factor is
$By_j$ with $j<\mu+1$ vanish; the terms
with final factor $By_{\mu+1}$ consist of
$Q_0$ and a higher-order term
${}^tP_{\mu+1}By_{\mu+1}$ with
$BP_{\mu+1}=0$.
Writing $\mathcal Q=Q_0^*\mathcal L_\lambda$,
we obtain
\begin{equation}\label{eq:m1-even-terminal}
 R\Gamma_c(U_{\bar k},g_1^*\mathcal L_\lambda)
 \cong H^0(S_{\mu,\bar k},\overline{\Q}_l)
 \otimes R\Gamma_c(V_{\mu+1,\bar k},\mathcal Q)
 (-e(\mu-1))[-2e(\mu-1)].
\end{equation}
The polarization map is
\[
 \ell_{\mathcal Q}=(F^\mu)^\dagger B+B^\dagger F^\mu
   =C_{\mu+1}.
\]
Under Condition~\ref{cond},
it is finite \'etale of degree $q^{ef}$,
so \cite[Theorem~4.2]{takeuchi2023quadraticelladicsheafheisenberg} implies
\[
 H_c^i(V_{\mu+1,\bar k},\mathcal Q)=0\ (i\ne e),
 \quad\dim H_c^e(V_{\mu+1,\bar k},\mathcal Q)
 =q^{ef/2}.
\]
Since $\#S_\mu(\bar k)=q^{ef(\mu-1)}$,
\eqref{eq:m1-even-terminal} is concentrated
in degree $e+2e(\mu-1)=N$, with the
claimed dimension.

Now suppose $\lambda\ne0$ violates Condition~\ref{cond}.
Define
\[
 r_1=\max\left\{1\leq r\leq r_0:
     \alpha_{j-1}\notin\mu_e
     \text{ for all }2\leq j\leq r\right\}.
\]
If $R=r_1+1\leq r_0$,
Lemma~\ref{lem:m1-relative} gives
\[
 R\Gamma_c(U_{\bar k},g_1^*\mathcal L_\lambda)
 \cong R\Gamma_c(W_R,h_R^*\mathcal L_\lambda)
 (-e(R-1))[-2e(R-1)],
\]
where
\[
 W_R=S_{R-1}\times K_R\times
 V_{R+1}\times\cdots\times V_{R^*-1}.
\]
By \eqref{eq:m1-Kr}, $W_R$ is a smooth affine
scheme of dimension
$d=1+e(f-2R+1)$. Hence
\[
 H_c^N(U_{\bar k},g_1^*\mathcal L_\lambda)
 \cong H_c^{d-1}(W_R,h_R^*\mathcal L_\lambda)
 (-e(R-1))=0,
\]
by affine vanishing and Poincar\'e duality.

The remaining case is $f=2\mu$, where
all paired indices are nondegenerate
but $\alpha_\mu\in\mu_e$.
The reduction \eqref{eq:m1-even-terminal}
still holds, while
$\Ker \ell_{\mathcal Q}$ has
one-dimensional identity component.
Then the middle degree vanishes
by Theorem~\ref{thm:Takeuchi}.
\end{prf}

\subsection{Proof of Proposition~\ref{prop:coh_dim} and \ref{prop:connected}}
\label{subsec:m1-proof-53}
It remains to prove (ii), since the other parts of Proposition~\ref{prop:coh_dim} are the same as $m>1$,
and Proposition~\ref{prop:connected} follows by Proposition~\ref{prop:coh_dim} in the same way.
Suppose $f>1$.
The middle-degree assertion was established in
Lemma~\ref{lem:m1-middle}.
To determine the highest degree,
recall that, for a rank-one local system $\mathcal F$ on a
geometrically connected smooth variety of dimension $N$,
$H_c^{2N}(U_{\bar k},\mathcal F)$ is the dual of
$H^0(U_{\bar k},\mathcal F^\vee)(-N)$, hence is nonzero
if and only if $\mathcal F$ is geometrically trivial,
in which case its dimension is one.

The statement for $e=1$ is established in \cite{MR3402698},
so we only treat the case $e>1$.
Restrict to the line
$y_{1,2}=t$ and $y_{i,j}=0$ for $(i,j)\ne(1,2)$.
The only composition contributing a nonzero term has all parts equal to one, and
\begin{equation}\label{eq:m1-line-e>1}
 g_1\big|_{\mathbb A^1}
 =(-1)^f t^{q+q^2+\cdots+q^f}.
\end{equation}
After extracting the initial $q$-power, its degree is
$1+q+\cdots+q^{f-1}$, which is prime to $p$ and exceeds one.
Consequently $\lambda g_1|_{\mathbb A^1}$ is not an
Artin-Schreier coboundary. Thus $g_1^*\mathcal L_\lambda$
is not geometrically trivial.

\subsection{Proof of Proposition~\ref{prop:trace}}
\label{subsec:m1-proof-56}
By the proof in the previous section,
we only need to check that
\[
\Tr(\Frob_{q^f}; H_c^N(X_{g_1})_\lambda),
\quad \dim H^i_c(X_{g_1}^\delta)_\lambda,\quad \delta=T,\,(u,u,1)
\]
have the same expressions as those of $m>1$.

\emph{Frobenius.}
If $f=2\mu+1$, by \eqref{eqn:f-odd},
we have
\begin{align*}
     \Tr(\Frob_{q^f}; H_c^N(X_{g_1})_\lambda)
     =(q^f)^{e\mu}|S_{\mu+1}(k)|
     =q^{(e+1)f\mu}
\end{align*}
as desired.

If $f=2\mu$, 
\eqref{eq:m1-even-terminal} gives
\begin{align*}
     \Tr(\Frob_{q^f}; H_c^N(X_{g_1})_\lambda)
     =(-1)^e(q^f)^{e(\mu-1)}|S_{\mu}(k)|\sum_{y\in k^e}\psi_\lambda(Q_0(y)).
\end{align*}
The last sum was computed in the previous section to be $(-1)^{e-1}Q^{e+1}$,
so that the trace is $-q^{(e+1)f(f-1)/2}$ as desired.
Both cases give the expected value.

\emph{Fixed locus of $\delta$.}
As in the case of $m>1$, restricting to the fixed locus of $T$
is the same as reducing to $e=1$,
and the computation is completely the same as the previous section.
The case of $\delta=(u,u,1)$ is also the same.

\subsection{Proof of Proposition~\ref{prop:forget-supp}}
\label{subsec:m1-proof-57}
All we need to alter from the case $m>1$
is the proof of the isomorphism
\begin{align*}
     H^{e(f-1)}_c(U,g_1^*\LL_\lambda)\xrightarrow{\sim} H^{e(f-1)}(U,g_1^*\LL_\lambda)
\end{align*}
for $\lambda$ satisfying Condition~\ref{cond}.
For $2\le r\le r_0$, the preceding elimination calculation gives
an isomorphism of local systems on $Z_r$,
\begin{equation}\label{eq:m1-57-elimination}
 (h_{r-1}|_{Z_r})^*\mathcal L_\lambda
 \cong
 \bigl({}^ty_{r^*}C_ry_r\bigr)^*\mathcal L_\lambda
 \otimes p_r^*(h_r|_{Z_r'})^*\mathcal L_\lambda.
\end{equation}
By Condition~\ref{cond}, the additive morphism
$C_r\colon V_r\to V_r$ is finite \'etale and surjective, and
$K_r=\Ker C_r=\Ker(F^f-T)$ is finite \'etale.

Put $\mathcal P_r=\bigl({}^ty_{r^*}C_ry_r\bigr)^*\mathcal L_\lambda$.
By the comparison isomorphism of Fourier-Deligne transform (cf. \cite[Theorem~5.12]{takeuchi2023quadraticelladicsheafheisenberg}),
we obtain
\begin{equation}\label{eq:m1-57-relative}
 Rp_{r!}\mathcal P_r\xrightarrow{\ \sim\ }Rp_{r*}\mathcal P_r.
\end{equation}
Therefore the elimination steps of $R\Gamma_c$ also work for $R\Gamma$
in a compatible way with the comparison maps.
More precisely, we obtain the following
commutative diagram in the derived category:
\begin{equation}\label{eq:m1-57-comparison-diagram}
\begin{array}{ccc}
 R\Gamma_c\bigl(Z_r,(h_{r-1}|_{Z_r})^*\mathcal L_\lambda\bigr)
 &\longrightarrow&
 R\Gamma\bigl(Z_r,(h_{r-1}|_{Z_r})^*\mathcal L_\lambda\bigr)
 \\
 \big\downarrow\wr && \big\downarrow\wr
 \\
 R\Gamma_c\bigl(Z_{r+1},(h_r|_{Z_{r+1}})^*\mathcal L_\lambda\bigr)
 (-e)[-2e]
 &\longrightarrow&
 R\Gamma\bigl(Z_{r+1},(h_r|_{Z_{r+1}})^*\mathcal L_\lambda\bigr)
 (-e)[-2e].
\end{array}
\end{equation}
Thus the comparison map
\begin{equation}\label{eq:m1-57-goal}
 R\Gamma_c(U,g_1^*\mathcal L_\lambda)
 \longrightarrow R\Gamma(U,g_1^*\mathcal L_\lambda)
\end{equation}
is identified, up to the same shift and Tate twist on both
sides, with the comparison map on the terminal space $Z_{r_0+1}$.

If $f=2\mu+1$, then $r_0=\mu+1$ and the terminal space is
\[
 S_{\mu+1}=K_2\times\cdots\times K_{\mu+1}.
\]
On this space $h_{\mu+1}=0$. Since $S_{\mu+1}$ is finite
\'etale over $k$, it is proper; hence
\[
 R\Gamma_c(S_{\mu+1},\overline{\Q}_l)
 \xrightarrow{\ \sim\ }
 R\Gamma(S_{\mu+1},\overline{\Q}_l).
\]

If $f=2\mu$, the terminal space is
$S_\mu\times V_{\mu+1}$.
The restriction of $h_\mu^*\LL_\lambda$
is the pullback of an isogeneous quadratic sheaf,
hence, by \cite[Proposition~5.13]{takeuchi2023quadraticelladicsheafheisenberg}, we have
\[
 H_c^e(S_\mu\times V_{\mu+1},h_{\mu}^*\LL_\lambda)
 \xrightarrow{\sim}
 H^e(S_\mu\times V_{\mu+1},h_{\mu}^*\LL_\lambda).
\]
This completes the proof.

\section{Representation theory}\label{sec:rep}

\subsection{Comparison with Yu's construction}\label{subsec:comparison}
Kaletha constructed regular supercuspidal representations from tame
elliptic regular pairs using Yu's construction, under suitable
assumptions on the residual characteristic \cite{MR4013740,MR1824988}.
We recall the relevant construction for $G=\GL_n(K)$ or $D^\times$.
In odd residual characteristic we use the usual Kaletha-Yu
construction, while in residual characteristic $2$ we use Takaya's
extension of Yu's construction
\cite{takaya2026characteristicfreeapproachesyusconstruction}.

Let $\theta$ be the character of $L^\times$ appearing in
Theorem~\ref{thm:main}. Thus $\theta$ is trivial on $U_L^{m+1}$,
its restriction to $U_L^m$ corresponds to $\psi_\lambda$, and
$\lambda$ satisfies Condition~\ref{cond}. The Howe factorization
can be chosen with a single positive-depth stage. In the notation
of Kaletha's construction, the corresponding Yu datum is $(\vec{G},y,\vec{r},\rho,\vec{\phi})$,
where
\begin{align*}
    \vec{G}=(L^\times, G), \quad \vec{\phi}=(\theta,1), \quad \rho=1,\quad \vec{r}=(m/e,m/e)
\end{align*}
and $y$ is the point of the reduced Bruhat-Tits building fixed by
the elliptic torus $L^\times$. Put $r_0=m/e$ and $s_0=r_0/2$.

We first verify the genericity conditions used in Yu's construction
and in Takaya's extension. In the terminology of
\cite[Definition~3.1]{takaya2026characteristicfreeapproachesyusconstruction},
these are the depth condition (GE0), the root-theoretic condition
(GE1), and the additional condition (GE2) required for
supergenericity. For the root datum of $\GL_n$ and its inner forms,
(GE2) follows from (GE1); see \cite[Lemma~8.1]{MR1824988}.

\begin{lem}\label{lem:takaya-genericity}
Suppose that $\lambda$ satisfies Condition~\ref{cond}.
Then $\theta$ is $(G,L^\times)$-supergeneric
of depth $r_0=m/e$ in the sense of
\cite[Definition~3.1]{takaya2026characteristicfreeapproachesyusconstruction},
that is, (GE0), (GE1), and (GE2) hold.
Furthermore, when $p=2$, the absence of symmetric ramified roots
required in
\cite[Assumption~2.14]{takaya2026characteristicfreeapproachesyusconstruction}
is automatic.
\end{lem}

\begin{prf}
Recall from Section~\ref{sec:main} that $\Psi\colon K\to\overline{\mathbb Q}_l^\times$
is an additive character
which is trivial on $\mathfrak p_K$ and whose restriction to
$\mathcal O_K$ induces $\psi$ on $\mathbb F_q$. Identify the dual
of $\Lie(L^\times)=L$ with $L$ by the trace pairing.
Let $\widetilde\lambda\in\mathcal O_{K_f}^\times$ be the
Teichm\"uller lift of $\lambda$.
Then the depth $-r_0$ element
\[
    X=\widetilde\lambda\varphi_L^{-m}\in L\cong \Lie^*(L^\times)
\]
represents $\theta|_{U_L^m}$,
which verifies (GE0).

Choose a tamely ramified Galois extension $E/K$ which splits
$L^\times$ and contains a primitive $e$-th root of unity $\zeta_e$.
The $n=ef$ embeddings
of $L$ into $E$ can be indexed by pairs $(i,j)$ with
$0\leq i<e$ and $0\leq j<f$, so that
\[
    \sigma_{i,j}(\varphi_L)=\zeta_e^i\varphi_L,
    \quad
    \sigma_{i,j}(\widetilde\lambda)
    =\widetilde\lambda^{q^j}.
\]
Under an identification of $G_E=G\otimes_K E$ with $\GL_{n,E}$, the element
$X$ is diagonal with entries
\[
    x_{ij}= \widetilde\lambda^{q^j} \zeta_e^{-im}\varphi_L^{-m}.
\]
For any root $\alpha=\alpha_{(i,j),(i',j')}$ of $G_E$
relative to $L_E^\times$, we have
\[
    X(H_\alpha)=x_{ij}-x_{i'j'}.
\]
Thus, to verify (GE1), it suffices to show that
\[
    v_L(x_{ij}-x_{i'j'})=-m
\]
for all distinct pairs $(i,j)$ and $(i',j')$,
where $v_L$ is extended to $E$.
Let $\overline\zeta_e$ denote the reduction of $\zeta_e$.
It is enough to show that the coefficients
$c_{ij}=\lambda^{q^j}\overline\zeta_e^{-im}$
are pairwise distinct.

Suppose that $c_{ij}=c_{i'j'}$.
If $j=j'$, then $\bar\zeta_e^{(i-i')m}=1$.
Since $\gcd(m,e)=1$ and $\overline\zeta_e$ has order $e$,
this implies $i=i'$.
If $j\ne j'$, taking the $e$-th power of the equality gives
$\lambda^{eq^j}=\lambda^{eq^{j'}}$,
which contradicts Condition~\ref{cond}.
Therefore the $c_{ij}$ are pairwise distinct, and hence
$v_L(x_{ij}-x_{i'j'})=-m$ for every root $\alpha$.
This proves (GE1).

Finally, if $p=2$, tameness implies that $e$ is odd.
Then \cite[Assumption~2.14]{takaya2026characteristicfreeapproachesyusconstruction}
is automatic as noted there.
\end{prf}

We now recall the construction of an irreducible supercuspidal
representation from this datum.
\begin{enumerate}
    \item
    Define the Moy-Prasad subgroups
    \begin{align*}
        \JJ=(L^\times,G)_{y,(r_0,s_0)},\quad
        \JJ_+=(L^\times,G)_{y,(r_0,s_0^+)},\quad
        \NN=(L^\times,G)_{y,(r_0^+,s_0^+)}.
    \end{align*}
    Let $\widetilde\psi_\lambda$ be the character of
    $\JJ_+$ whose restriction to $U_L^m$ corresponds to
    $\psi_\lambda$ and which is trivial on the root subgroup
    factors. This character is trivial on $\NN$.

    \item
    The quotient $\HH=\JJ/\NN$ fits into a central extension
    \[
       1\to \ZZ=\JJ_+/\NN
       \to \HH
       \to W=\JJ/\JJ_+
       \to 1.
    \]
    The commutator pairing
    \[
       W\times W\to\overline{\mathbb Q}_l^\times,
       \quad
       (\overline x,\overline y)\mapsto
       \widetilde\psi_\lambda([x,y])
    \]
    is nondegenerate by (GE1); this is proved in
    \cite[Lemma~11.1]{MR1824988} for odd residue characteristic,
    and the same argument is valid for $p=2$ by
    \cite[Lemma~3.2]{takaya2026characteristicfreeapproachesyusconstruction}.
    In particular, $\ZZ$ is the center of $\HH$, and
    $\HH/\Ker\widetilde\psi_\lambda$ is a finite
    Heisenberg group.
    By the Stone-von Neumann theorem, $\HH$ has, up to
    isomorphism, a unique irreducible representation
    $\tau_\lambda'$ with central character
    $\widetilde\psi_\lambda$. We use the same notation
    for its inflation to $\JJ$.

    \item
    We extend $\tau_\lambda'$ to a
    representation $\kappa_\theta$ of $L^\times\JJ$
    by combining $\theta$ with a Heisenberg-Weil extension.
    If $p\ne2$, we use the usual extension from the
    Kaletha-Yu construction \cite{MR4013740,MR1824988}.
    If $p=2$, we use Takaya's construction
    \cite[Theorem~2.27 and Lemma~3.7]
          {takaya2026characteristicfreeapproachesyusconstruction}.
    In either case, $\kappa_\theta|_\JJ\simeq
    \tau_\lambda'$.

    \item
    Set
    \[
        \pi_{G,\theta}'
        =\cInd_{L^\times\JJ}^{G}\kappa_\theta.
    \]
    When $p\ne 2$, this representation is irreducible and
    supercuspidal by 
    \cite[Theorem~15.1]{MR1824988}; see also the corrected proof
    \cite[Theorem~3.1]{MR4357723}.
    When $p=2$, the same conclusion follows from
    \cite[Theorem~3.8]{takaya2026characteristicfreeapproachesyusconstruction},
    using the genericity verified above.
\end{enumerate}

We are going to investigate the representations of $S_{\GL}^m$ and $S_D^m$
in comparison to the Heisenberg representation of $\JJ/\NN$.
We write $S^m$ for $S_{\GL}^m$ or $S_D^m$ and $Z$ for its center.
We write $\J$ for $\J_{\GL}^m$ or $\J_D^m$, $\J_+$ for $U_L^m \J_{\GL}^{m+1}$ or $U_L^m \J_D^{m+1}$
and $\mathfrak N$ for $\J_{\GL}^{m+1}$ or $\J_D^{m+1}$, respectively.
Then we have the following sequence of subgroups:
\[\J_+\subset \JJ_+\subset \JJ\subset \J.\]
These four subgroups are of the form $1+\p_L^m+(C_1\cap \PA^i)$ or $1+\p_L^m+(C_2\cap \PD^i)$,
where the integer $i$ is as in the following tables.
\begin{center}
    {
        \renewcommand{\arraystretch}{1.2}
    \begin{tabular}[t]{|c||c|c|}\hline
        $\GL$ & $m$: odd & $m$: even\\
         & \vphantom{$f$: odd} & \\ \hhline{|=#=|=|}
        $\J$ & $\frac{m+1}{2}$ & $\frac{m}{2}$\\
        $\JJ$  & $\frac{m+1}{2}$ & $\frac{m}{2}$\\
        $\JJ_+$ & $\frac{m+1}{2}$ & $\frac{m}{2}+1$\\
        $\J_+$ & $\frac{m+1}{2}$ & $\frac{m}{2}+1$\\ \hline
    \end{tabular}
    \quad
    \begin{tabular}[t]{|c||cc|c|}\hline
        $D$ & \multicolumn{2}{c|}{$m$: odd} & $m$: even\\
         & $f$: odd & $f$: even & \\ \hhline{|=#==|=|}
        $\J$ & \multicolumn{2}{c|}{$\frac{f(m-1)}{2}+1$} & $\frac{fm}{2}$\\
        $\JJ$ & $\frac{fm+1}{2}$ & $\frac{fm}{2}$ & $\frac{fm}{2}$\\
        $\JJ_+$ & $\frac{fm+1}{2}$ & $\frac{fm}{2}+1$ & $\frac{fm}{2}+1$\\
        $\J_+$ & \multicolumn{2}{c|}{$\frac{f(m+1)}{2}$} & $\frac{fm}{2}+1$\\ \hline
    \end{tabular}
    }
\end{center}
Let $\langle , \rangle_\lambda$ denote the pairing on both $\JJ/\JJ_+$ and $\J/\J_+$
induced by the commutator and $\psi_\lambda$.

\begin{prop}\label{prop:nondeg}
    If $m\geq 2$, then $S^m$ is a Heisenberg group or an abelian group,
    and for $\lambda\in \F_{q^f}$ satisfying Condition~\ref{cond}, the pairing $\langle , \rangle_\lambda$
    on $\JJ/\JJ_+\cong S^m/Z$ is non-degenerate.
    Therefore, there exists a unique irreducible representation of $S^m$ with the central character $\psi_\lambda$.
\end{prop}
\begin{prf}
    The first statement is obvious from the multiplication rule of $S^m$,
    so we only need to check that the pairing $\langle , \rangle_\lambda$
    on $\JJ/\JJ_+\cong S^m/Z$ is non-degenerate.
    We work with the case where $m$ is even and $S^m/Z=S_{\GL}^m/Z(S_{\GL}^m)=\big\{(a_{ij})\in \prod_{i=1}^e\prod_{j=0}^{f-1}\F_{q^f}\mid \sum_{i=1}^e a_{i0}=0\big\}$.
    Then the pairing is as follows:
    \begin{align*}
        \langle a,b\rangle_\lambda 
        &=\psi\Big(\Tr\Big(\lambda \Big(\sum_{i,j} a_{ij}b_{i+s,-j}^{q^j}-b_{ij}a_{i+s,-j}^{q^j}\Big)\Big)\Big)\\
        &=\psi \Big(\Tr \Big( \sum_{i,j}\Big( \lambda^{q^j}a_{i-s,-j}^{q^j}-\lambda a_{i+s,-j}^{q^j} \Big)b_{ij} \Big)\Big).
    \end{align*}
    Therefore, the map $b\mapsto \langle a,b\rangle_\lambda$ is zero
    if and only if
    \begin{itemize}
        \item $\lambda^{q^j}a_{i-s,-j}^{q^j}-\lambda a_{i+s,-j}^{q^j}=0$ for all $i$ and $j\neq 0$, and
        \item $a_{i-s,0}- a_{i+s,0}$ is independent of $i$.
    \end{itemize}
    For $j\neq 0$, the condition is equivalent to
    \begin{align*}
        a_{i,j}=\frac{\lambda^{q^j}}{\lambda}a_{i+m,j}
    \end{align*}
    for all $i$ and $j\neq 0$, and Condition~\ref{cond} and the fact that $m$ and $e$ are coprime
    implies that $a_{i,j}=0$ for all $i$ and $j\neq 0$.
    For $j=0$, the condition forces $a_{i,0}$ to be constant,
    but the sum $\sum_i a_{i,0}$ must be zero,
    so the fact that $e$ is coprime to $p$ implies that $a_{i,0}=0$ for all $i$.
    Thus the pairing $\langle ,\rangle_\lambda$ is non-degenerate.
    We can prove the other cases similarly, so we omit them.
\end{prf}

\begin{rem}
    If $m=1$, $S^m$ need not be a Heisenberg group.
    However, as we will see below, the uniqueness of the irreducible representation of $S^m$
    with the central character $\psi_\lambda$
    also follows even if $m=1$.
\end{rem}

\begin{defn}
    Suppose $\lambda\in \F_{q^f}$ satisfies Condition~\ref{cond}.
    We define the representation of $\J$ as
    $\tau_\lambda=\Ind_{\JJ}^{\J}\tau_\lambda'$.
    This factors through a representation of $S^m$,
    which is also denoted as $\tau_\lambda$.
\end{defn}
Let $\overline{\JJ}$ and $\overline{\JJ}_+$ denote the image of $\JJ$ and $\JJ_+$ in $S^m=\J/\mathfrak N$, respectively.
Then $\tau_\lambda$ is also described as $\Ind_{\overline \JJ}^{S^m}\tau_\lambda'$.

\begin{prop}\label{prop:Heisenberg}
    Suppose $\lambda\in \F_{q^f}$ satisfies Condition~\ref{cond}.
    For any integer $m\geq 1$, the representation $\tau_\lambda$
    is the unique irreducible representation of $S^m$ with the central character $\psi_\lambda$.
    Its dimension is equal to $\sqrt{|S^m/Z|}$.
    That is, we have the following table.
    \[
    \begin{tabular}{|c|cc|}\hline
         & $\dim \tau_{\GL,\lambda}$ & $\dim \tau_{D,\lambda}$\\ \hline
        $m$: odd & $1$ & $q^{n(f-1)/2}$\\
        $m$: even & $q^{f(n-1)/2}$ & $q^{(e-1)f/2}$\\ \hline
    \end{tabular}
    \]
\end{prop}
\begin{prf}
    The statement for dimension follows by the equation
    $\dim \tau_\lambda=|S^m/\overline \JJ|\dim \tau_\lambda'=\sqrt{|\J/\J_+|}=\sqrt{|S^m/Z|}$.

    First we see that it suffices to prove that $\tau_\lambda$ is irreducible.
    Indeed, we have $\Ind_{Z}^{S^m}\psi_\lambda\cong \bigoplus_\tau \tau^{\dim \tau}$,
    where $\tau$ runs through all irreducible representations with central character $\psi_\lambda$.
    Thus we have $|S^m/Z|=\sum_\tau (\dim\tau)^2$,
    but $(\dim \tau_\lambda)^2$ is already $|S^m/Z|$,
    so the uniqueness follows.

    When $m>1$, or $m=1$ and $G=\GL_n(K)$,
    since $S^m$ is also a Heisenberg group or an abelian group,
    the irreducibility of $\tau_\lambda$ is deduced by the dimension counting.

    When $m=1$ and $G=D^\times$, 
    if further $f=1$, $\J=\J_+$ and $S^1$ is abelian,
    so the claim is obvious.
    Suppose $f>1$.
    We have to check that if $x\in S^1=\F_{q^f}\times \F_{q^n}^{f-1}$ satisfies $\widetilde{\psi}_\lambda^x\cong \widetilde{\psi}_\lambda$,
    then $x\in \overline{\JJ}$.
    If $f=2$, $\overline{\JJ}=S^1$,
    so the claim is again obvious.
    Hence we only consider the case of $f>2$.
    We suppose $x\not\in \overline{\JJ}$ and prove that there exists $y\in \overline{\JJ}_+$ such that 
    $\widetilde{\psi}_\lambda(xyx^{-1})\neq \widetilde{\psi}_\lambda(y)$,
    that is, $\widetilde{\psi}_\lambda([x,y])\neq 1$.
    We write $x=(a_0,a_1,\ldots,a_{f-1})$.
    There exists an integer $1\leq j< \lceil f/2\rceil$ such that $a_j\neq 0$,
    and take $j$ as the minimum.
    For any $b\in \F_{q^n}$, we put $y_b$ for the element of $\overline{\JJ}_+$
    whose $(f-j)$-th coefficient is $b$ and others are zero.
    Then simple computation shows that $[x,y_b]\in Z$
    and $[x,y_b]=\Tr_{\F_{q^n}/\F_{q^f}}(a_j b^{q^j}-b a_j^{q^{f-j}})$.
    In a similar way to the previous computation,
    we see that Condition~\ref{cond} implies that the pairing $(a,b)\mapsto \psi_\lambda(\Tr(ab^{q^j}-ba^{q^{f-j}}))$ is non-degenerate,
    so there exists $b$ such that $\widetilde{\psi}_\lambda([x,y_b])\neq 1$.
\end{prf}

\subsection{Proof of Theorem~\ref{thm:main} (i) and Theorem~\ref{thm:characterize}}
This subsection is devoted to proving Theorem~\ref{thm:main} (i) and Theorem~\ref{thm:characterize}.
Thus we fix $\lambda\in \F_{q^f}$ and a character $\theta$ of $L^\times$
such that
$\theta|_{U_L^m}
=\psi_{\widetilde{\lambda}\varphi_L^{-m}}^L$,
where $\widetilde{\lambda}\in\OO_{K_f}$ is the
Teichm\"uller lift of $\lambda\in\F_{q^f}$.

\begin{lem}\label{lem:central_character}
    The subgroup $\J_{\GL}^{m+1}\times \J_D^{m+1}$ acts on $\Hom_{U_L^{m,1}}(\theta,H^{n-1}_c(X_m^0,\overline{\Q}_l))$ trivially,
    and $U_L^1 \times U_L^1$ acts on it by $\theta^{-1}\boxtimes \theta$.
\end{lem}
\begin{prf}
    Let $d\in \J_D^{m+1}$ and consider the action of $(1,d)\in \J_{\GL}^{m+1}\times \J_D^{m+1}$ on $f\in \Hom_{U_L^{m,1}}(\theta,H^{n-1}_c(X_m^0,\overline{\Q}_l))$.
    By the definition of the action (cf.\ Remark~\ref{rem:action}),
    for $\sigma\in W_L$ such that $(1,d,\sigma)\in J_{m,s}^1$, we have
    $((1,d)f)(v)=(1,d,\sigma)(f(\sigma^{-1}v))$.
    Since $\Nrd(\J_D^{m+1})=N_{L/K}(U_L^{m+1})$, we can take $\sigma$ as $\Art_L(u)$ for $u\in U_L^{m+1}$.
    Indeed, if we take $u\in U_L^{m+1}$ such that $\Nrd(d)=N_{L/K}(u)^{-1}$,
    then $(1,d,\Art_L(u))=(1,du,1)(1,u^{-1},\Art_L(u))$ is an element of $J_{m,s}^1$.
    However, the action of $(1,du,1)$ and $(1,u^{-1},\Art_L(u))$ on $X_m$ is trivial,
    and $\sigma^{-1}v=v$ since $\theta$ is trivial on $U_L^{m+1}$,
    thus we get $((1,d)f)(v)=f(v)$, which implies that $\{1\}\times\J_D^{m+1}$ acts trivially.
    For the action of $g\in \J_{\GL}^{m+1}$, we can prove the statement similarly,
    except for that we consider the decomposition $(g,1,\sigma)=(gu^{-1},1,1)(1,u^{-1},\Art_L(u))(u,u,1)\in J_{m,s}^1$.
    Since the action of $(u,u,1)$ on $X_m$ is trivial, the argument works the same.
    
    Next we consider the action of $(1,u)$ where $u\in U_L^1$.
    Since $(1,u,\Art_L(u)^{-1})$ acts on $X_m$ trivially,
    we have $((1,u)f)(v)=(1,u,\Art_L(u)^{-1})(f(\theta(u)v))=\theta(u)f(v)$,
    thus $\{1\}\times U_L^1$ acts by $\theta$.
    In addition, since $(u,u,1)$ acts on $X_m$ trivially,
    the action of $(u,u)$ is also trivial, which implies that $U_L^1\times \{1\}$ acts by $\theta^{-1}$.
\end{prf}

\begin{lem}\label{lem:Hom}
    We have
    \begin{align*}
        \Hom_{U_L^{m,1}}(\theta, H^{n-1}_c(X_m^0,\overline{\Q}_l))&\cong H^{n-1}_c(X_m,\overline{\Q}_l)[\psi_\lambda^{-1}].
    \end{align*}
\end{lem}
\begin{prf}
    If $e\geq 2$, since $e$ does not divide $m$, we see that $N_{L/K}(U_L^m)=N_{L/K}(U_L^{m+1})$,
    which implies that the composition $U_L^{m,1}\hookrightarrow U_L^m\to U_L^m/U_L^{m+1}\cong \F_{q^f}$ is surjective.
    Also, by Proposition~\ref{prop:connected}, we have $X_m^0=X_m$.
    Therefore the lemma follows from Lemma~\ref{lem:central_character}.
    If $e=1$, then the image of $U_L^{m,1}\to \F_{q^f}$ is the trace-zero subspace $\F_{q^f}^0$ of $\F_{q^f}$, so we have
    \begin{align*}
        \Hom_{U_L^{m,1}}(\theta, H^{n-1}_c(X_m^0))
        &= \Hom_{\F_{q^f}^0} (\psi_\lambda|_{\F_{q^f}^0}, H^{n-1}_c(X_m^0))\\
        &\cong \Hom_{\F_{q^f}}(\psi_\lambda, \Ind_{\F_{q^f}^0}^{\F_{q^f}}H^{n-1}_c(X_m^0))\\
        &\cong \Hom_{\F_{q^f}}(\psi_\lambda, H^{n-1}_c(X_m))\\
        &\cong H^{n-1}_c(X_m)[\psi_\lambda^{-1}].
    \end{align*}
    From the second line to the third line, we used the fact that,
    by Proposition~\ref{prop:connected}, 
    $\F_{q^f}/\F_{q^f}^0\cong \F_q$ acts on the component set of $X_m$
    simply transitively.
\end{prf}

Now we are in a position to prove Theorem~\ref{thm:main} (i) and Theorem~\ref{thm:characterize}.

\begin{prf}[of Theorem~\ref{thm:main} (i)]
    By Lemma~\ref{lem:Hom}, we have
    \begin{align*}
        \Hom_{U_L^{m,1}}(\theta, H^{n-1}_c(X_m^0,\overline{\Q}_l))&\cong H^{n-1}_c(X_m,\overline{\Q}_l)[\psi_\lambda^{-1}]\\
        &\cong H^{e-1}_c(\A^{e-1},f_m^*\LL_{-\lambda})\otimes H^{e(f-1)}_c(\A^{e(f-1)},g_m^*\LL_{-\lambda}).
    \end{align*}
    By Proposition~\ref{prop:coh_dim} this is nonzero if and only if $\lambda$ satisfies Condition~\ref{cond},
    and its dimension is
    \[
    \begin{cases}
        q^{ef(f-1)/2} & \text{ if $m$ is odd,}\\
        q^{f(e-1)+ef(f-1)/2} & \text{ if $m$ is even.}
    \end{cases}
    \]
    On the other hand, by Lemma~\ref{lem:central_character},
    the action of $\J_{\GL}^m\times \J_D^m$ on this space factors through $S_{\GL}^m\times S_D^m$,
    and its center acts by $\psi_\lambda^{-1}\boxtimes \psi_\lambda$.
    Therefore it is isomorphic to $(\tau_{\GL,\lambda}^\vee \boxtimes \tau_{D,\lambda})^a$ for some $a\in \Z_{\geq 0}$.
    However, $\dim \tau_{\GL,\lambda}\dim \tau_{D,\lambda}$ is exactly the same as the dimension of $\Hom_{U_L^{m,1}}(\theta, H^{n-1}_c(X_m^0,\overline{\Q}_l))$,
    which forces $a=1$.
    This finishes the proof.
\end{prf}

\begin{prf}[of Theorem~\ref{thm:characterize}]
    We have already proved the properties for the action of $\J_{\GL}^m$ (resp.\ $\J_D^m$) and $U_L^1$.
    Since the canonical map
    \begin{align*}
        \Hom_{U_L^{m,1}}(\theta, H^{n-1}_c(X_m^0,\overline{\Q}_l))\otimes \theta\to H^{n-1}_c(X_m)[\psi_\lambda^{-1}]
    \end{align*}
    is $J_{m,s}^1$-equivariant, for any element $(g,d,\Art_L(x))\in J_{m,s}^1$, we have
    \begin{align*}
        \Tr \overline{\pi}_{\GL,\theta}^\vee(g) \Tr \overline{\pi}_{D,\theta}(d) \theta(x) = \Tr ((g,d,\Art_L(x)); H^{n-1}_c(X_m)[\psi_\lambda^{-1}]).
    \end{align*}
    The value of the right hand side is obtained by Proposition~\ref{prop:trace},
    taking the product of the two columns of the tables there.
    If we take $(1,\varphi_L^{-1},\Art_L(\varphi_L))$,
    we can compute the value of $\Tr \overline{\pi}_{D,\theta}(\varphi_L^{-1})$,
    and then $(\varphi_L^{-1},\varphi_L^{-1},1)$ gives us the value of $\Tr \overline{\pi}_{\GL,\theta}^\vee(\varphi_L^{-1})=\Tr \overline{\pi}_{\GL,\theta}(\varphi_L)$.
    Similarly, $(1,u,\Art_L(u)^{-1})$ and $(u,u,1)$ give us the values of
    $\Tr \overline{\pi}_{\GL,\theta}(u^{-1})$ and $\Tr \overline{\pi}_{D,\theta}(u)$.
\end{prf}

\subsection{Irreducibility}
In this subsection, we continue to treat the cases
$G=\GL_n(K)$ and $G=D^\times$ simultaneously.
We write $\overline{\pi}_\theta$ for
$\overline{\pi}_{\GL,\theta}$ or $\overline{\pi}_{D,\theta}$,
and $\pi_\theta$ for $\pi_{\GL,\theta}$ or $\pi_{D,\theta}$,
according to $G$.
We also write $\pi_\theta'$ for the representation
$\pi_{G,\theta}'$ constructed in Subsection~\ref{subsec:comparison}.
\begin{lem}\label{lem:extension}
    We have
    \begin{align*}
        (\Ind_{L^\times \JJ}^{L^\times \J}\kappa_\theta)|_{\J}\cong \overline{\pi}_\theta|_{\J}.
    \end{align*}
\end{lem}
\begin{prf}
    By Mackey's formula, we have
    \begin{align*}
        (\Ind_{L^\times \JJ}^{L^\times \J}\kappa_\theta)|_{\J}
        &\cong \Ind_{L^\times \JJ\cap \J}^{\J}(\kappa_\theta|_{L^\times \JJ\cap \J})
        = \Ind_{\JJ}^{\J}(\kappa_\theta|_{\JJ})\\
        &\cong \Ind_{\JJ}^{\J}(\tau_\lambda')
        =\tau_\lambda=\overline{\pi}_\theta|_{\J}.
    \end{align*}
\end{prf}

\begin{lem}\label{lem:twist}
    If $\rho_1$ and $\rho_2$ are two representations of $L^\times \J$
    whose restrictions are isomorphic with an irreducible representation $\tau$ on $\J$,
    then there exists a character $\chi$ of $L^\times$, trivial on $L^\times\cap \J=U_L^m$,
    such that $\rho_1=\chi \rho_2$.
    Here, $\chi$ is viewed as a character of $L^\times \J$ by making it to be trivial on $\J$.
\end{lem}
\begin{prf}
    We set $\chi=\rho_1\rho_2^{-1}$.
    Then for any $x\in L^\times \J$ and $j\in \J$, we have
    \begin{align*}
        \chi(x)\tau(j)&=\rho_1(x)\rho_2(x)^{-1}\tau(j)\\
        &=\rho_1(x)\rho_2(x^{-1}jx)\rho_2(x)^{-1}\\
        &=\rho_1(jx)\rho_2(x)^{-1}\\
        &=\tau(j)\chi(x),
    \end{align*}
    so that $\chi(x)\in \overline{\Q}_l^\times$ by Schur's lemma. 
    For $x,y\in L^\times \J$, we have
    \begin{align*}
        \chi(xy)&=\rho_1(x)\rho_1(y)\rho_2(y)^{-1}\rho_2(x)^{-1}\\
        &=\rho_1(x)\chi(y)\rho_2(x)^{-1}\\
        &=\chi(x)\chi(y).
    \end{align*}
    Therefore $\chi$ is a character of $L^\times \J$.
    This is clearly trivial on $\J$.
\end{prf}

\begin{prop}\label{prop:comparison}
    There exists a character $\chi$ of $L^\times$ such that
    $\Ind_{L^\times\JJ}^{L^\times \J}\kappa_{\theta\chi}\cong \overline{\pi}_{\theta}$.
    In particular, we have $\pi_\theta\cong \pi_{\theta\chi}'$,
    which is an irreducible supercuspidal representation of $G$.
\end{prop}
\begin{prf}
    By Lemma~\ref{lem:extension} and Lemma~\ref{lem:twist},
    $\chi=\overline{\pi}_\theta(\Ind_{L^\times \JJ}^{L^\times \J}\kappa_\theta)^{-1}$ is
    a character of $L^\times\J$, trivial on $U_L^m$.
    By the construction of $\kappa_\theta$, we see that $\kappa_{\theta\chi}=\chi\kappa_\theta$, therefore
    \begin{align*}
        \Ind_{L^\times \JJ}^{L^\times \J}\kappa_{\theta\chi}
        =\Ind_{L^\times \JJ}^{L^\times \J}(\chi\kappa_\theta)
        =\chi(\Ind_{L^\times \JJ}^{L^\times \J}\kappa_\theta)
        =\overline{\pi}_\theta.
    \end{align*}
\end{prf}

\section{Affinoids of the Lubin-Tate space}\label{sec:affinoid}
\subsection{Non-abelian Lubin-Tate theory}
We briefly recall the basic notation for the Lubin-Tate theory.
Let $\OO_{\breve{K}}$ be the completion of the ring of integers $\OO_{K^{\mathrm{ur}}}$
of the maximal unramified extension of $K$
and $\breve{K}$ be the fraction of $\OO_{\breve{K}}$.
We fix a one-dimensional formal $\OO_K$-module $H_0$ over $\overline{\F}_q$ of height $n$,
which is unique up to isomorphism.
The $\OO_K$-algebra $\End_{\OO_K}{H_0}$ can be identified with
the maximal order $\OO_D$ of the central division algebra $D$ with invariant $1/n$.
If we take $H_0$ as the reduction of the standard formal $\OO_K$-module $H_{\st}$ of height $n$ over $\OO_{\breve{K}}$
in the sense of \cite[Section~2.3]{MR3402698},
then $\End_{\OO_K}H_0$ contains the elements $\Pi\colon T\mapsto T^q$ and $[\zeta]\colon T\mapsto \overline{\zeta}T$ for $\zeta\in \mu_{q^n-1}(K_n)$,
and the identification $\OO_D \xrightarrow{\sim} \End_{\OO_K}H_0$ can be taken
as $\varphi_D\mapsto \Pi$ and  $\zeta\mapsto [\zeta]$.

Let $\widehat{\mathcal{C}}$ be the category of complete Noetherian local $\OO_{\breve{K}}$-algebras
whose residue fields are isomorphic to $\overline{\F}_q$.
For each integer $m\geq 0$, we consider the functor $\MM_{H_0,m}^{(j)}$ from $\widehat{\mathcal{C}}$ to the category of sets
which sends $A\in \widehat{\mathcal{C}}$ to the set of equivalence classes
of triples $(H,\iota, \eta)$,
where $H$ is a one-dimensional formal $\OO_K$-module over $A$,
$\iota\colon H_0\to H\otimes_A \overline{\F}_q$ is a quasi-isogeny of height $j$,
and $\eta$ is a Drinfeld $m$-level structure of $H$.
We write $\MM_{H_0,m}=\bigsqcup_{j\in \Z}\MM_{H_0,m}^{(j)}$.
The obtained projective system $\{\MM_{H_0,m}\}$ is called the Lubin-Tate tower.
Each $\MM_{H_0,m}$ is equipped with an action of $D^\times$ and a Weil descent datum,
and $\GL_n(K)$ acts on the pro-object $\{\MM_{H_0,m}\}$.
We can write $\MM_{H_0,m}^{(j)}=\Spf A_m^{(j)}$ for $A_m^{(j)}\in \widehat{\mathcal{C}}$.
If $n=1$, it is known that $A_m^{(j)}\cong \OO_{\breve{K}_m}$, where $\breve K_m$ is the $m$-th Lubin-Tate extension of $K^{\mathrm{ur}}$.

We denote by $M_{H_0,m}$, $M_{H_0,m}^{(j)}$ the rigid generic fibers of $\MM_{H_0,m}$, $\MM_{H_0,m}^{(j)}$, respectively.
Then we obtain a projective system $\{M_{H_0,m}\}$ of rigid spaces,
whose transition maps are finite and \'etale.
We put
\[H^i_{\mathrm{LT}}=\varinjlim_m H_c^i(M_{H_0,m}\otimes_{\breve{K}} C, \overline{\Q}_l),\]
where $C$ denotes the completion of $\overline{K}$.
It is equipped with an action of $\GL_n(K)\times D^\times \times W_K$,
where the action of $W_K$ is defined by the Weil descent datum on $\MM_{H_0,m}$.
The action of $\GL_n(K)\times D^\times$ is known to be smooth,
and we consider the discrete topology on $W_K$, so that the action of $W_K$ is also smooth.

For a technical reason, we fix $0\neq \gamma\in \p_K$
and consider the discrete subgroup $\Gamma=\gamma^\Z\subset K^\times$,
which is cocompact in $D^\times$.
We put
\[H^i_{\mathrm{LT},\Gamma}=\varinjlim_m H_c^i((M_{H_0,m}/\Gamma)\otimes_{\breve{K}} C, \overline{\Q}_l),\]
which is a representation of $\GL_n(K)/\Gamma \times D^\times/\Gamma \times W_K$.

The following theorem is known as the non-abelian Lubin-Tate theory.
\begin{thm}[\cite{MR1876802}, \cite{MR1719811}, \cite{MR2511742}, \cite{MR2680204}]\label{thm:NALT}
    \begin{enumerate}
        \item For an irreducible supercuspidal representation $\pi$ of $\GL_n(K)$, we have
        \begin{align*}
            \Hom_{\GL_n(K)}(H_{\mathrm{LT}}^i,\pi)\cong
            \begin{cases}
                \mathrm{JL}(\pi)\boxtimes \rec_K(\pi)(\frac{n-1}{2}) & \text{if $i=n-1$},\\
                0 & \text{otherwise}
            \end{cases}
        \end{align*}
        as representations of $D^\times \times W_K$.
        \item For an irreducible supercuspidal representation $\pi$ of $\GL_n(K)$
        whose central character is trivial on $\Gamma$,
        we write $H_{\mathrm{LT},\Gamma,\pi}^i$
        for the $\pi$-isotypic part of $H_{\mathrm{LT},\Gamma}^i$.
        Then we have
        \begin{align*}
            H_{\mathrm{LT},\Gamma,\pi}^i\cong
            \begin{cases}
                \pi\boxtimes \mathrm{JL}(\pi^\vee)\boxtimes \rec_K(\pi^\vee)(-\frac{n-1}{2}) & \text{if $i=n-1$},\\
                0 & \text{otherwise}
            \end{cases}
        \end{align*}
        as representations of $\GL_n(K)\times D^\times \times W_K$.
    \end{enumerate}
\end{thm}

\subsection{Lubin-Tate space at infinite level}
Next we recall the Lubin-Tate space at infinite level.
We fix $H=H_{\st}$ as the standard formal $\OO_K$-module
and write $H_0$ for its reduction.
We have a one-dimensional formal $\OO_K$-module $\wedge H_0$ of height $1$
whose Dieudonn\'e module $M(\wedge H_0)$ is the top exterior product of the Dieudonn\'e module $M(H_0)$ of $H_0$.
We choose a lift $\wedge H$ of $\wedge H_0$
over $\OO_{\breve K}$ such that
\begin{align*}
    \log_H(T)=\sum_{i\geq 0}\frac{T^{q^{in}}}{\varpi^i}, \quad \log_{\wedge H}(T)=\sum_{i\geq 0}(-1)^{(n-1)i}\frac{T^{q^i}}{\varpi^i}.
\end{align*}

The $K$-vector space object $\widetilde{H}:=\varprojlim_{[\varpi]}H$ is called the universal cover of $H$,
and it depends only on $H_0$.
More precisely, any endomorphism of $H_0$ lifts naturally to an endomorphism of $\widetilde{H}$ (\cite[Lemma~2.5.3 (3)]{MR3402698}),
and we use the same symbol to denote the lift.
The subset of elements $(\ldots, s_2,s_1)$ such that $s_1\in H[\varpi]$
is an $\OO_K$-submodule, called the Tate module $TH$.
The $K$-subspace $VH=TH\otimes_{\OO_K} K\subset \widetilde{H}$ is called the rational Tate module,
and it is comprised of such elements with $s_1\in H[\varpi^\infty]$.

We write $\Nil^{\flat}$ for the functor from $\widehat{\mathcal{C}}$ to the category of sets
which sends $A\in \widehat{\mathcal{C}}$ to $\varprojlim_q \Nil(A)$
where $\Nil(A)$ is the set of topologically nilpotent elements in $A$
and the limit is taken along the $q$-power Frobenius.
This functor $\Nil^\flat$ is obviously represented by the formal scheme $\Spf \mathcal{O}_{\breve{K}}[[T^{q^{-\infty}}]]$.
It is known that there exist natural isomorphisms 
$\lambda\colon \widetilde{H} \xrightarrow{\sim} \Nil^\flat$
and $\lambda'\colon \widetilde{\wedge H}\xrightarrow{\sim} \Nil^{\flat}$.
We write $\lambda_0\colon\widetilde H\to\Nil$
for the composite of $\lambda$ with the projection
onto the first component of $\Nil^\flat$.
Thus, if $\lambda(s)=(T,T^{1/q},T^{1/q^2},\ldots)$,
then $\lambda_0(s)=T$.
The composite maps $\Nil^\flat\xrightarrow{\lambda^{-1}}\widetilde{H}\xrightarrow{\pr_1}H\xrightarrow{\log_H}\G_a$
and $\Nil^\flat\xrightarrow{\lambda'^{-1}}\widetilde{\wedge H}\xrightarrow{\pr_1}\wedge H\xrightarrow{\log_{\wedge H}}\G_a$
are expressed as $\Tb \mapsto h(T)$ and $\Tb \mapsto h'(T)$ respectively, where
\begin{align*}
    h(T)=\sum_{i\in \Z}\frac{T^{q^{in}}}{\varpi^i}, \quad h'(T)=\sum_{i\in \Z}(-1)^{(n-1)i}\frac{T^{q^i}}{\varpi^i}.
\end{align*}
From now on, we use boldface letters, such as $\Tb$, for sections of $\Nil^\flat$,
whereas we use the corresponding Roman letters, such as $T$,
for the first component of the system $\Tb=\{T^{q^{-r}}\}_{r\geq 0}$.

We write $A_\infty^{(j)}$ for the completion of $\varinjlim_m A_m^{(j)}$
and set $\MM_{H_0,\infty}^{(j)}=\Spf A_\infty^{(j)}$ and $\MM_{H_0,\infty}=\bigsqcup_{j\in\Z} \MM_{H_0,\infty}^{(j)}$.
A point $x\in \MM_{H_0,\infty}(A)$ for $A\in \widehat{\mathcal{C}}$
corresponds to a
deformation $(H,\iota)$ of $H_0$ together with a
compatible system of Drinfeld level structures.
In \cite[Section~2.7]{MR3402698}, they defined the quasilogarithm map
$\qlog_H\colon\widetilde{H}(A)\to M(H_0)\otimes (A\otimes K)$.
They chose a convenient basis $e_1,\ldots ,e_n$ of $M(H_0)\otimes K$, called the standard basis,
for which $\qlog_H$ is expressed as
\begin{align}\label{eqn:qlog}
    s\mapsto (h(\lambda_0(s)), h(\lambda_0(s)^q),\ldots, h(\lambda_0(s)^{q^{n-1}})),
\end{align}
where $h(T)$ is as above
(\cite[Lemma~2.7.1]{MR3402698}).
For the basis $e_1\wedge \cdots \wedge e_n$ of the rank one module $M(\wedge H_0)\otimes K$,
$\qlog_{\wedge H}$ is of the form $(\ldots,s_2,s_1)\mapsto \log_{\wedge H}s_1$.
Thus the kernel of $\qlog_{\wedge H}$ is the rational Tate module $V(\wedge H)$.
It is proved in \cite[Theorem~2.10.3]{MR3402698} that
there exists a $K$-alternating map $\delta\colon \widetilde{H}^n\to \widetilde{\wedge H}$
such that the following diagram commutes:
\[
\xymatrix{
\widetilde{H}(A)^n \ar[r]^\delta \ar[d]_{\qlog_H\times \cdots \times \qlog_H}&
\widetilde{\wedge H}(A) \ar[d]^{\qlog_{\wedge H}}\\
M(H_0)^n\otimes (A\otimes K) \ar[r]^{\det}&
M(\wedge H_0)\otimes (A\otimes K).
}
\]
We also have the map $\MM_{H_0,\infty}\to \widetilde{H}^n$,
induced by the level structures.
In this setting Weinstein found the following Cartesian diagram (\cite[Theorem 2.7.3]{MR3529120}):
\[
\xymatrix{
\MM_{H_0,\infty} \ar[r] \ar[d]&
\MM_{\wedge H_0,\infty} \ar[d]\\
\widetilde{H}^n \ar[r]^\delta&
\widetilde{\wedge H}.
}
\]
We also write $\delta$ for the top morphism $\MM_{H_0,\infty}\to \MM_{\wedge H_0,\infty}$. 
Our computation is carried out using the coordinate of $\Nil^\flat$,
so we restate the above Cartesian diagram by $\lambda$ and $\lambda'$ as follows:
\[
\xymatrix{
\MM_{H_0,\infty} \ar[r] \ar[d]&
\MM_{\wedge H_0,\infty} \ar[d]\\
\Nil^{\flat,n} \ar[r]^\Delta&
\Nil^{\flat}.
}
\]


We now pass to the adic generic fibers.
We define
$\MM_{H_0,\infty,\OO_C}
:=
\MM_{H_0,\infty}\widehat{\otimes}_{\OO_{\breve K}}\OO_C$.
For a topological $\OO_C$-algebra $A$,
let $t(\Spf A)$ denote the adic space $\Spa(A,A)$
and let $t(\Spf A)_\eta$ denote its generic fiber
$t(\Spf A)\setminus V(\varpi)$.
These constructions extend to formal schemes
by gluing.

We put
$M_{\infty,C}
:=t(\MM_{H_0,\infty,\OO_C})_\eta$,
which is a perfectoid space.
We also write
$\widetilde H_{\eta,C}
:=t(\widetilde H_{\OO_C})_\eta$.
The generic fiber of $\MM_{\wedge H_0,\infty,\OO_C}$
is identified with the space of nonzero elements
of $V(\wedge H)$.
Consequently, the Cartesian diagram above gives
the following description of the $C$-points of
$M_{\infty,C}$:
\begin{align*}
M_{\infty,C}(C)
&=\left\{
x=(x_1,\ldots,x_n)\in\widetilde H_{\eta,C}(C)^n
\,\middle|\,
\delta(x)\neq 0,\,
\qlog_{\wedge H}\delta(x)=0
\right\}\\
&=\left\{
x\in\widetilde H_{\eta,C}(C)^n
\,\middle|\,
\delta(x)\neq 0,\,
\det(\qlog_H(x_1),\ldots,\qlog_H(x_n))=0
\right\}.
\end{align*}
The formal scheme $\widetilde{H}^n_{\OO_C}$ is equipped with the action of $\GG=\GL_n(K)\times D^\times \times W_K$,
which stabilizes $\MM_{H_0,\infty, \OO_C}$.

We use subscripts to distinguish operations on different objects;
for the $\mathcal{O}_K$-action and the sum operation on $H$ we write $[a]_H$ and $x+_H y$ respectively.
We use the same notation for the action and operation on $\widetilde{H}$.
We write $[a]_{\flat}$ (resp.\ $[a]_{\flat,\wedge}$) for $a\in K$
and $\bm{x}+_{\flat} \bm{y}$ (resp.\ $\bm{x}+_{\flat,\wedge} \bm{y}$) for $\bm{x},\bm{y}\in \Nil^\flat(A)$ 
for the corresponding operations on $\Nil^\flat$ induced by $\lambda$ (resp.\ $\lambda'$).
By abuse of notation, we write $x+_\flat y$ for the first component of $\bm{x}+_\flat \bm{y}$
(similar for $x+_{\fw} y$, $[a]_\flat(x)$ and $[a]_{\fw}(x)$).

We write the representing algebra of $\Nil^{\flat,n}_{\mathcal{O}_C}$ as $R=\mathcal{O}_C[[X_1^{q^{-\infty}},\dots,X_n^{q^{-\infty}}]]$
and that of $\Nil^\flat_{\mathcal{O}_C}$ as $\mathcal{O}_C[[T^{q^{-\infty}}]]$.
For these coordinates, the action of $\GG$ on $\Nil_{\mathcal{O}_C}^{\flat,n}$ is described as follows:
\begin{itemize}
    \item For $g\in \GL_n(K)$, $g^*\Xb_i = (\flat)\sum_{j=1}^n [g_{ji}]_\flat(\Xb_j)$.
    Here, $(\flat)\sum$ means that the sum is computed using $+_{\flat}$.
    \item For $d\in D^\times$, $d^*\Xb_i = [d^{-1}]_\flat(\Xb_i)$.
    \item For $\sigma\in W_K$, $\sigma^*\Xb_i=\Xb_i^{q^{-m_\sigma}}$ and $\sigma^*a=\sigma(a)$ for $a\in \mathcal{O}_C$.
    Here, $m_\sigma=v_K(\Art_K^{-1}(\sigma))$.
\end{itemize}

\subsection{Affinoids}
We continue to take $H_0$ as the reduction of the standard formal $\OO_K$-module $H_{\st}$ of height $n$.
We identify $\OO_D$ with $\End_{\OO_K}H_0$ by $\varphi_D\mapsto \Pi$ and $\zeta\mapsto [\zeta]$ for $\zeta\in \mu_{q^n-1}(K_n)$.
Then the $K_f$-embedding $i^D\colon L\hookrightarrow D$ in Subsection~\ref{subsec:D}
allows us to view $H_0$ as a formal $\OO_L$-module of height $1$.
We take a lift $H$ of $H_0$ to $\OO_C$, a formal $\OO_L$-module of height $1$.
We write $T_LH(\OO_C)=\varprojlim_m H[\varphi_L^m]$ in order to distinguish from $TH(\OO_C)=\varprojlim_m H[\varpi^m]$, the Tate module when $H$ is viewed as a formal $\OO_K$-module of height $n$.
These are naturally isomorphic over $\OO_L$, via
\begin{align*}
    T_LH(\OO_C)\xrightarrow{\sim} TH(\OO_C);\, (\ldots,t_2,t_1)\mapsto (\ldots,t_{2e},t_e).
\end{align*}
Since $T_LH(\OO_C)$ is a rank one $\OO_L$-module, we can choose a basis $\widetilde{\xi}_{0,L}$.
Let $\widetilde{\xi}_0$ be the image in $TH(\OO_C)$ and $\xxi_0=\lambda(\widetilde{\xi}_0)\in \Nil^\flat(\OO_C)$.
It is easy to see from the construction that $|\xi_0|=|\varpi|^{1/eq^{f(e-1)}(q^f-1)}$.

We put 
\[A_0=\begin{pmatrix}
    \alpha_1 & \alpha_2 & \cdots & \alpha_f\\
    \Phi(\alpha_1) & \Phi(\alpha_2) & \cdots & \Phi(\alpha_f)\\
    \vdots & \vdots & \ddots & \vdots\\
    \Phi^{f-1}(\alpha_1) & \Phi^{f-1}(\alpha_2) & \cdots & \Phi^{f-1}(\alpha_f)
\end{pmatrix}\in M_f(K_f)\]
and $A=I_e\otimes A_0\in M_n(K_f)$.
We set
\[\vec{\xxi}'=(\varphi_L^{e-1},\varphi_L^{e-2},\dots ,1)\otimes (\xxi_0,0,\dots ,0),\quad \vec{\xxi}=\vec{\xxi}' A.\]
Then $\vec{\xxi}$ is an $\OO_K$-basis of $TH(\OO_C)$,
which gives an $\OO_C$-point $x_0$ of $\MM_{H_0,\infty,\OO_C}$ with CM by $\OO_L$
whose corresponding embeddings in the sense of \cite[Section 3.1]{MR3402698}
are $i^G$ and $i^D$ taken in Subsections~\ref{subsec:GL} and~\ref{subsec:D}.
The point $\tilde{t}:=\delta(x_0)\in \MM_{\wedge H_0,\infty}^{(0)}(\OO_C)$
gives an isomorphism $(\varinjlim_m \OO_{\breve K_m})^\wedge\xrightarrow{\sim} \OO_{\widehat{K}^{\ab}}$,
where $\widehat{K}^{\ab}\subset C$ is the completion of the maximal abelian extension of $K$.
Thus $A_\infty^{(0)}$ is endowed with the $\OO_{\widehat{K}^{\ab}}$-algebra structure,
and we use it to form $\MM_{H_0,\infty,t}^{(0)}:=\Spf (A_\infty^{(0)}\widehat{\otimes}_{\OO_{\widehat{K}^{\ab}}}\OO_C)$,
the fiber of $\tilde{t}$ along $\delta$.
We put $\bm{t}:=\lambda'(\tilde{t})\in \Nil^\flat(\OO_C)$.

We relabel the coordinates of $\Nil_{\OO_C}^{\flat,n}$ as
\[R=\OO_C[[X_{ij}^{q^{-\infty}}\mid 1\leq i\leq e, 1\leq j\leq f]].\]
We write $\vec{X}$ for the horizontal vector $(\vec{X}_1,\ldots, \vec{X}_f)$
where $\vec{X}_j=(X_{1j},\ldots, X_{ej})$ for $j\in \{1,\ldots,f\}$.
We take a new coordinate $\vec{Y}=(\vec{Y}_1,\ldots ,\vec{Y}_f)$ making the CM point $\vec{\xi}$ the origin,
together with the auxiliary coordinate $Z$ that will be used as a substitute for $Y_{11}$ to define an appropriate affinoid.
\begin{defn}\label{def:affinoid}
    Let $m\geq 1$ be an integer.
    We define $\vec{Y}=(Y_{ij})_{1\leq i\leq e,1\leq j\leq f}$ and $Z$ by
    \[\vec{X}=\vec{\xi}+_\flat \vec{Y}A=(\vec{\xi}'+_\flat \vec{Y})A,\quad Z=(\flat)\sum_{i=1}^e [\varphi_L^{i-1}]_\flat (Y_{i1}).\]
    Let $\widetilde{U}_m$ be the rational subset of $t(\Nil_{\OO_C}^{\flat,n})$ defined by
    \begin{align*}
        |Y_{ij}|\leq |\xi_{i1}'|^{c_{m,j}}\neq 0,\quad |Z|\leq |\xi_{11}'|^{q^{fm}}\neq 0,
    \end{align*}
    where
    \begin{itemize}
        \item if $m$ is even, $c_{m,j}=q^{fm/2}$ for any $j$, and
        \item if $m$ is odd,
        \begin{align*}
            c_{m,j}=\begin{dcases}
                q^{f(m-1)/2}\frac{1+q^f}{2} & \text{if $j=1$},\\
                q^{f(m-1)/2+j-1} & \text{if $j\geq 2$}.
            \end{dcases}
        \end{align*}
    \end{itemize}
    We define the rational subset $U_m$ of $t(\MM_{H_0,\infty,\OO_C})$
    as the inverse image of $\widetilde{U}_m$ under the map $t(\MM_{H_0,\infty,\OO_C})\to t(\Nil^{\flat,n}_{\OO_C})$.
    We also define $U_m^t=U_m\cap t(\MM_{H_0,\infty,t}^{(0)})$.
\end{defn}

\subsection{Comparison with \cite{MR4126896}}
We apply \cite{MR4126896} to the totally ramified extension $L/K_f$ of degree $e$.
Since $H_0$ can also be regarded as a formal $\mathcal O_{K_f}$-module
of height $e$, we may consider the corresponding Lubin-Tate space
at infinite level $\mathcal M'_{H_0,\infty}$.
For a deformation $H'$ of $H_0$ as a formal $\mathcal O_{K_f}$-module,
a $K_f$-basis $(\eta'_1,\ldots,\eta'_e)$ of its rational Tate
module $VH'$ gives the following $K$-basis:
\[
\bigl(
 [\alpha_1]_{H'}(\eta'_1),\ldots,[\alpha_1]_{H'}(\eta'_e),
 \ldots,
 [\alpha_f]_{H'}(\eta'_1),\ldots,[\alpha_f]_{H'}(\eta'_e)
\bigr).
\]
This construction gives a natural morphism
$\mathcal M'_{H_0,\infty}\to\mathcal M_{H_0,\infty}$,
and the following diagram is commutative:
\[
\xymatrix{
\MM_{H_0,\infty}' \ar[r] \ar[d]&
\MM_{H_0,\infty} \ar[r] \ar[d]&
\MM_{\wedge H_0,\infty} \ar[d]\\
\Nil^{\flat,e} \ar@{^{(}->}[r]&
\Nil^{\flat,n} \ar[r]^\Delta&
\Nil^{\flat}.
}
\]
Here the embedding
$\operatorname{Nil}^{\flat,e}\hookrightarrow
\operatorname{Nil}^{\flat,n}$ is given by
\[
(X_1,\ldots,X_e)
\longmapsto
(X_1,\ldots,X_e,0,\ldots,0)A.
\]
Under this embedding, a point with $Y$-coordinates
$(Y_1,\ldots,Y_e)$ in the notation of \cite{MR4126896} maps to a point of
$\operatorname{Nil}^{\flat,n}$ whose $Y$-coordinates in our notation are
\[(Y_1,\ldots,Y_e,0,\ldots,0).\]
Moreover, the CM point $\vec{\xi}\in\operatorname{Nil}^{\flat,n}$
lies in the image of $\operatorname{Nil}^{\flat,e}$.
The restriction to $t(\mathcal M'_{H_0,\infty,\mathcal O_C})$
of the affinoid in $t(\mathcal M_{H_0,\infty,\mathcal O_C})$
defined in Definition~\ref{def:affinoid} therefore coincides with the affinoid
defined in \cite[Section~3.2]{MR4126896}, applied to the totally ramified
extension $L/K_f$ of degree $e$.

\subsection{Reduction of the affinoids}\label{subsec:reduction}

The main result of this section is the following theorem.
\begin{thm}\label{thm:red}
    Let $m\geq 1$ be an integer coprime to $e$.
    \begin{enumerate}
        \item The subgroup $J_{m,s}^1$ stabilizes $U_m^t$ and $U_m$.
        \item There exists a formal model $\mathcal{X}_m$ (resp.\ $\mathcal{X}_m^t$) over $\mathcal{O}_C$ of $U_m$ (resp.\ $U_m^t$)
        equipped with a $J_{m,s}^1$-action.
        \item There exists a $J_{m,s}^1$-equivariant map $\overline{\mathcal{X}}_m\to X_m^{\perf}$,
        restricting to a $J_{m,s}^1$-equivariant isomorphism $\overline{\mathcal{X}}_m^t \xrightarrow{\sim} X_m^{0,\perf}$.
        Here, the bars represent the special fibers of the formal models.
    \end{enumerate}
\end{thm}

We explain how we are going to prove this theorem.
Throughout this section, we assume $\gcd(m,e)=1$.
We put $\mathcal{Y}_m=\Spf S_+$, where $S_+=\OO_C\langle z'^{q^{-\infty}}, y_{ij}'^{q^{-\infty}}\mid (i,j)\neq (1,1)\rangle$,
and regard it as a formal model of $\widetilde{U}_m$ by
\begin{align*}
    \bm y_{ij}'=\bm \xi_{i1}'^{-c_{m,j}}\bm Y_{ij},\quad \bm z'=\bm \xi_{11}'^{-q^{fm}}\bm Z.
\end{align*}
Let $\overline{m}$ denote the representative of $m$ in $\{1,\ldots,e-1\}$ modulo $e$ if $e>1$,
and otherwise set $\overline{m}=0$.
We define a rational subset $\wedge \widetilde{U}_m$ of $t(\Nil^\flat_{\OO_C})$ as 
\begin{align*}
    |T -_{\fw} t|\leq |\tau|^{q^r},
    \quad r=\lfloor m/e \rfloor,\quad \tau=\xi_0^{\nu_m}\det A
\end{align*}
where $\nu_m=\overline{m}q^n+(e(1+q+\cdots+q^{f-1})-\overline{m})q^{n-f}$.
We define $\wedge \mathcal{Y}_m=\Spf \OO_C\langle u^{q^{-\infty}}\rangle$
and regard it as a formal model of $\wedge \widetilde{U}_m$ by
\begin{align*}
    \bm u=\bm \tau^{-q^r}(\bm T-_{\fw} \bm t).
\end{align*}
We will prove that $\Delta-_{\fw}t=(-1)^{(n-1)r}(\tau g)^{q^r}$ for some $g\in S_+$.
This implies that $\Delta$ sends $\widetilde{U}_m$ into $\wedge \widetilde{U}_m$,
and that this morphism lifts to
the morphism of formal schemes $\mathcal Y_m\to \wedge \mathcal Y_m$
defined by
\begin{align*}
    \OO_C\langle u^{q^{-\infty}}\rangle \to S_+,
    \quad \bm u\mapsto \bm \tau^{-q^r}(\bm \Delta(\bm z',\bm y_{ij}')-_{\fw} \bm t).
\end{align*}
Let $\wedge U_m$ be the preimage of $\wedge \widetilde{U}_m$ in $t(\MM_{\wedge H_0,\infty,\OO_C})$.
\begin{lem}\label{lem:translate}
    The pullback $\wedge U_m$ of $\wedge \widetilde{U}_m$ in $t(\MM_{\wedge H_0,\infty,\OO_C})$
    is the $G_m$-translates of $\bm t$, where $G_m:=N_{L/K}(U_L^m)=U_K^{\lceil m/e\rceil}$.
    In addition, if $e>1$, then every element $x\in G_m$ satisfies
    \begin{align*}
        |[x]_{\fw}(t)-_{\fw} t|<|\tau|^{q^r}
    \end{align*}
    for $r=\lfloor m/e \rfloor$.
\end{lem}
\begin{prf}
    By \cite[Lemma~3.9]{MR4126896}, for $x\in K^\times$,
    we have $v([x]_{\fw}(t)-_{\fw} t)=v([x-1]_{\fw}(t))=q^{v(x-1)}v(t)$,
    where $v$ is the valuation normalized as $v(\varpi)=1$.
    Therefore $[x]_{\fw}(t)\in \wedge \widetilde{U}_m$ if and only if
    $q^{v(x-1)}v(t)\geq q^rv(\tau)$.
    Since $v(t)=1/(q-1)$ and $v(\tau)=1/(q-1)+\overline{m}/e$,
    this is equivalent to $q^{v(x-1)}\geq q^r(1+\frac{\overline{m}}{e}(q-1))$.
    We have
    \begin{align*}
        q^r\Big(1+\frac{\overline{m}}{e}(q-1)\Big)
        \begin{cases}
            \in (q^r,q^{r+1}) & \text{if $e>1$},\\
            =q^m & \text{if $e=1$},
        \end{cases}
    \end{align*}
    so the lemma follows.
\end{prf}
By the above lemma, $\wedge \mathcal{X}_m=\Spf \Cont(G_m,\OO_C)$
is a formal model of $\wedge U_m$,
and the morphism $\wedge U_m\to \wedge \widetilde{U}_m$ lifts to
$\wedge \mathcal{X}_m\to \wedge \mathcal{Y}_m$ defined by
\begin{align*}
    \OO_C\langle u^{q^{-\infty}}\rangle \to \Cont(G_m,\OO_C),\quad
    \bm u\mapsto \Big( x\mapsto \bm \tau^{-q^r}([x]_{\fw}(\bm t)-_{\fw}\bm t) \Big).
\end{align*} 
We also consider the morphism $\Spf \OO_C\to \wedge \mathcal{X}_m$
corresponding to $\bm t$, that is,
the evaluation at $1\in G_m$ for $\Cont(G_m,\OO_C)\to \OO_C$.
Then we define the formal schemes $\mathcal{X}_m$ and $\mathcal{X}_m^t$
by the following Cartesian diagram:
\[
\xymatrix{
\mathcal{X}_m^t \ar[r] \ar[d]&
\Spf \OO_C \ar[d]\\
\mathcal{X}_m \ar[r] \ar[d]&
\wedge \mathcal{X}_m \ar[d]\\
\mathcal{Y}_m \ar[r]^\Delta&
\wedge \mathcal{Y}_m.
}
\]
The special fibers give the following Cartesian diagram:
\[
\xymatrix{
\overline{\mathcal{X}}_m^t \ar[r] \ar[d]&
\{1\} \ar[d]\\
\overline{\mathcal{X}}_m \ar[r] \ar[d]&
G_m \ar[d]\\
\A^{n,\perf} \ar[r]^{\overline{g}}&
\A^{1,\perf},
}
\]
where $G_m$ denotes the constant scheme $\Spec \Cont(G_m,\overline{\F}_q)$
and $\{1\}$ denotes $\Spec \overline{\F}_q$.
Lemma~\ref{lem:translate} implies that the morphism $G_m\to \A^{1,\perf}$
factors through $G_m/G_{m+1}$, which is trivial unless $e=1$,
in which case it is isomorphic to $\F_q$.
On the other hand,
we will also see that $h'(\Delta)=\varpi^r\tau f$ for some $f\in S_+$
and $X_m^{\perf}$ is given by $V(\overline{f})^{\perf}$.
In addition, we will see that $\overline{f}=\overline{g}$ if $e>1$
and $\overline{f}=\overline{g}-\overline{g}^q$ if $e=1$,
so that $X_m^{\perf}$ is the pullback of $G_m/G_{m+1}$ by $\overline{g}\colon \A^{n,\perf}\to \A^{1,\perf}$.
This gives the morphism $\overline{\mathcal{X}}_m\to X_m^{\perf}$
and the identification $\overline{\mathcal{X}}_m^t\cong V(\overline{g})^{\perf}$.
The connectedness follows by Proposition~\ref{prop:connected}.

When $e=1$, our construction reduces to the unramified case
treated by Boyarchenko and Weinstein in \cite{MR3402698},
and the assertions of Theorem~\ref{thm:red} follow from
their results. Hence, for the remainder of the proof,
we assume $e>1$.

We begin the proof of Theorem~\ref{thm:red} with the approximation of $h'(\Delta)$.
Since $\qlog_{\wedge H}\delta(x)=\det(\qlog_H (x_i)_{1\leq i\leq n})$
and $\qlog_H$ is expressed as (\ref{eqn:qlog}),
we have
\begin{align*}
    h'(\Delta)=\qlog_{\wedge H}\delta|_{\Nil^{\flat,n}}
    =\begin{vmatrix}
        h(X_{11}) & \cdots & h(X_{ef})\\
        h(X_{11}^q) & \cdots & h(X_{ef}^q)\\
        \vdots & & \vdots\\
        h(X_{11}^{q^{n-1}}) & \cdots & h(X_{ef}^{q^{n-1}})
    \end{vmatrix}.
\end{align*}
\begin{lem}
    Let $R$ be a non-commutative ring and $M$ be a left $R$-module.
    For horizontal vectors $\vec{x}\in M^n$ and $\vec{r}\in R^n$, and an invertible matrix $P\in M_n(R)$,
    we have
    \begin{align*}
        {}^t\! \vec{r}(\vec{x}P)=
        \begin{pmatrix}
            \vec{x}\cdot{}^P\! r_1\\
            \vdots\\
            \vec{x}\cdot{}^P\! r_n
        \end{pmatrix}P.
    \end{align*}
    Here, for an element $r\in R$, ${}^P\! r$ denotes the matrix $P(rI_n)P^{-1}$.
\end{lem}
\begin{prf}
    This follows by direct calculation.
\end{prf}
For the following lemma, we recall that $i^G$ is defined by the tensor product of $i^G_1\colon E_e\hookrightarrow M_e(K)$ and $i^G_2\colon K_f\hookrightarrow M_f(K)$,
where each one is induced by the ordered basis $(\varphi_L^{e-1},\ldots, 1)$ and $(\alpha_1,\ldots ,\alpha_f)$ respectively,
and $s\in M_f(K)$ is the representing matrix of the Frobenius of $K_f/K$.
\begin{lem}\label{lem:conj}
    \begin{enumerate}
        \item For $x\in K_f$, we have
        \begin{align*}
            Ai^G(x)A^{-1}=I_e\otimes\begin{pmatrix}
                x & & & \\
                 & \Phi(x) & & \\
                 & & \ddots & \\
                 & & & \Phi^{f-1}(x)
            \end{pmatrix}.
        \end{align*}
        \item $A_0 s A_0^{-1}=\begin{pmatrix}
             & I_{f-1}\\
            1 & 
        \end{pmatrix}$.
        \item ${}^A\Pi=A(\Pi I_n)A^{-1}=\Big(I_e\otimes \begin{pmatrix}
             & I_{f-1}\\
            1 & 
        \end{pmatrix}\Big)\Pi$.
    \end{enumerate}
\end{lem}
\begin{prf}
    (i) Since $i^G(x)=1\otimes i^G_2(x)$, we only need to compute $A_0 i^G_2(x) A_0^{-1}$.
    By definition, we have $\begin{pmatrix}
        \alpha_1 & \cdots & \alpha_f
    \end{pmatrix} i^G_2(x)=\begin{pmatrix}
        x\alpha_1 & \cdots & x\alpha_f
    \end{pmatrix}$.
    We let $\Phi$ act on this equation to get $\begin{pmatrix}
        \Phi^i(\alpha_1) & \cdots & \Phi^i(\alpha_f)
    \end{pmatrix} i^G_2(x)=\begin{pmatrix}
        \Phi^i(x)\Phi^i(\alpha_1) & \cdots & \Phi^i(x)\Phi^i(\alpha_f)
    \end{pmatrix}$, as desired.
    The remaining statements follow from similar computations.
\end{prf}
By the above lemmas, we have
\begin{align*}
    {}^t\!(1,\Pi,\ldots, \Pi^{n-1})\vec{X}
    &={}^t\!(1,\Pi,\ldots, \Pi^{n-1}) ((\vec{\xi}'+_\flat \vec{Y})A)\\
    &=\begin{pmatrix}
        \vec{\xi}'+_\flat \vec{Y}\\
        (\vec{\xi}'+_\flat \vec{Y})\cdot {}^A \Pi\\
        \vdots\\
        (\vec{\xi}'+_\flat \vec{Y})\cdot {}^A \Pi^{n-1}
    \end{pmatrix}A\\
    &=\begin{pmatrix}
        \vec{\xi}_1'+_\flat \vec{Y}_1 & \vec{Y}_2 & \cdots & \vec{Y}_f\\
        \Pi\vec{Y}_f & \Pi(\vec{\xi}_1'+_\flat \vec{Y}_1) & \cdots & \Pi\vec{Y}_{f-1}\\
        \vdots & \vdots & \ddots & \vdots \\
        \Pi^{n-1}\vec{Y}_2 & \Pi^{n-1}\vec{Y}_3 & \cdots & \Pi^{n-1}(\vec{\xi}_1'+_\flat \vec{Y}_1)
    \end{pmatrix}A.
\end{align*}
For two matrices $P$ and $Q$, we write $P\sim Q$ when $Q$ is obtained by
permuting rows and columns of $P$.
Then we have
\begin{align*}
    {}^t\!(1,\Pi,\ldots, \Pi^{n-1})\vec{X}
    \sim \begin{pmatrix}
        \Yf_1 & \Yf_2 & \cdots & \Yf_f\\
        \Pi\Yf_f & \Pi\Yf_1 & \cdots & \Pi\Yf_{f-1}\\
        \vdots & \vdots & \ddots & \vdots \\
        \Pi^{f-1}\Yf_2 & \Pi^{f-1}\Yf_3 & \cdots & \Pi^{f-1}\Yf_1
    \end{pmatrix}A',
\end{align*}
where $\Yf_j={}^t\!(1,\Pi^f,\ldots ,\Pi^{(e-1)f})(\vec{\xi}_j'+_\flat \vec{Y}_j)$
and $A\sim A'$ (note that $\vec{\xi}_j'=0$ for $j\geq 2$).
By definition, we have $h(X+_\flat Y)=h(X)+h(Y)$ and $h([a]_\flat(X))=ah(X)$ for $a\in \OO_K$,
so the above computation leads
\begin{align*}
    \pm h'(\Delta)=\det h\begin{pmatrix}
        \Yf_1 & \Yf_2 & \cdots & \Yf_f\\
        \Pi\Yf_f & \Pi\Yf_1 & \cdots & \Pi\Yf_{f-1}\\
        \vdots & \vdots & \ddots & \vdots \\
        \Pi^{f-1}\Yf_2 & \Pi^{f-1}\Yf_3 & \cdots & \Pi^{f-1}\Yf_1
    \end{pmatrix}\cdot \det A,
\end{align*}
where we let $h(P)$ denote the matrix $(h(p_{ij}))$ for a matrix $P=(p_{ij})$.
We write $\det h(\Xi)$ for the first determinant and view it as a polynomial on $h(\xi_0^{q^l})$, $h(Y_{ij}^{q^l})$ and $h(Z^{q^l})$
for $l\in \Z_{\geq 0}$ and $(i,j)\in (\{1,\ldots, e\}\times \{1,\ldots, f\})\setminus {(1,1)}$.
We assign degree $0$ to $h(\xi_0^{q^l})$ and degree $1$ to $h(Y_{ij}^{q^l})$ and $h(Z^{q^l})$,
and regard the resulting expression as a graded polynomial.
We approximate this separately for the grading.
Before we start approximation, we prepare the following lemma.
\begin{lem}\label{lem:special-bound}
Suppose $e>1$.
Put
\begin{align*}
 c=v(\xi_0)=\frac{1}{e q^{n-f}(q^f-1)},\qquad
 b_l=cq^{n-f+l}-\frac{n-f+l}{n}\quad(0\le l<f).
\end{align*}
Let $a\in\mathbb Z$ and let $T$ satisfy $v(T)\ge cq^a$. Write $a\equiv l\pmod f$ with $0\le l<f$. Then
\begin{align*}
 v(h(T))\ge \frac an+b_l.
\end{align*}
We call $a/n+b_l$ the \emph{special lower bound}. Write $a=nr+w$, where $r\in\mathbb Z$ and $1\le w\le n$. We call $w$ \emph{special} if
\begin{align*}
 \begin{cases}
    w=n-f\ \text{or}\ n &\text{if }l=0,\\
    w=n-f+l &\text{if }1\le l<f.
 \end{cases}
\end{align*}
The special lower bound is attained if and only if $v(T)=cq^a$ and $w$ is special.
In this case the unique summand of lowest valuation in $h(T)$ is indexed by $i=-r$, and
$h(T)\equiv\varpi^r T^{q^{-nr}}$.
\end{lem}
\begin{prf}
Recall that $h(T)=\sum_{i\in\mathbb Z}\varpi^{-i}T^{q^{in}}$. For integers $t$ put
\begin{align*}
 F(t)&=cq^t-\frac tn.
\end{align*}
We have
\begin{align*}
 F(t+f)-F(t)
 &=cq^t(q^f-1)-\frac fn
 =\frac1e\bigl(q^{t-(n-f)}-1\bigr).
\end{align*}
Thus $F(t+f)<F(t)$ for $t<n-f$, equality holds for $t=n-f$, and $F(t+f)>F(t)$ for $t>n-f$. It follows that for $t\equiv l\pmod f$,
\begin{align*}
 F(t)\ge b_l,
\end{align*}
with equality exactly when $t=n-f$ or $n$ for $l=0$, and when $t=n-f+l$ for $1\le l<f$.
For each $i\in\mathbb Z$,
\begin{align*}
 v\bigl(\varpi^{-i}T^{q^{in}}\bigr)
 &\ge-i+cq^{a+in}
 =\frac an+F(a+in)
 \ge\frac an+b_l.
\end{align*}
This proves the lower bound. Equality requires $v(T)=cq^a$ and one of the above minimizing values of $a+in$. Since $a=nr+w$, this is equivalent to $w$ being special; the corresponding index is $i=-r$. The two minimizing integers $n-f$ and $n$ for $l=0$ differ by $f<n$, so at most one of them belongs to the sequence $a+n\mathbb Z$. Thus the smallest-valuation summand is unique, as asserted.
\end{prf}

Let $\Gamma=(\gamma_{ij})\in M_n(S_+)$ and consider the leading terms of $\det h(\Gamma)$. Suppose that $v(\gamma_{ij})=cq^{m_{ij}}$ for $m_{ij}\in\mathbb Z$ and put
\begin{align*}
 M=(m_{ij})=nR+W,
 \qquad R=(r_{ij})\in M_n(\mathbb Z),\quad
 W=(w_{ij})\in M_n(\{1,\ldots,n\}).
\end{align*}
For a fixed subset $S\subset \mathfrak{S}_n$ we study the leading terms in
\begin{align*}
 \sum_{\sigma\in S}\sgn(\sigma)
 \prod_{i=1}^n h(\gamma_{i\sigma(i)}).
\end{align*}
In each exponent matrix below, the entries in the $j$-th block column (or the $j$-th block row in (III)(b)) are congruent to $j-1$ modulo $f$.
Each such block contains $e$ columns (respectively, rows); as a permutation uses each column and each row exactly once, every residue class $l=0,\ldots,f-1$ occurs exactly $e$ times among $m_{i\sigma(i)}\bmod f$.
Lemma~\ref{lem:special-bound} therefore yields
\begin{align*}
 v\left(\prod_{i=1}^n h(\gamma_{i\sigma(i)})\right)
 &\ge\frac1n\sum_{i=1}^n m_{i\sigma(i)}
       +e\sum_{l=0}^{f-1}b_l.
\end{align*}
Since
\begin{align*}
 e\sum_{l=0}^{f-1}b_l
 &=\frac1{q-1}-n+\frac{f+1}{2},
\end{align*}
putting $\mathfrak M=n(n-(f+1)/2)$ gives
\begin{align*}
 v\left(\prod_{i=1}^n h(\gamma_{i\sigma(i)})\right)
 &\ge\frac1{q-1}
 +\frac{\sum_i m_{i\sigma(i)}-\mathfrak M}{n}.
\end{align*}
We call this the \emph{special lower bound for the determinant term}. The special lower bound is attained if and only if every selected $w_{i\sigma(i)}$ is special.

We use the following strategy.

\noindent\textbf{Step (i).} We compute $\sum_i m_{i\sigma(i)}$ and, in the cases below, show that it is independent of $\sigma\in S$. This determines a common special lower bound for the determinant terms.

\noindent\textbf{Step (ii).} We determine the permutations for which every $w_{i\sigma(i)}$ is special. These are precisely the permutations attaining the special lower bound. Their distribution of residue values is also determined as follows. Let
\begin{align*}
 d=\#\{i\mid w_{i\sigma(i)}=n\}.
\end{align*}
Then $n$ occurs $d$ times, $n-f$ occurs $e-d$ times, and $n-f+l$ occurs $e$ times for each $1\le l<f$. Consequently,
\begin{align*}
 \sum_i w_{i\sigma(i)}
 &=dn+(e-d)(n-f)+e\sum_{l=1}^{f-1}(n-f+l)
 =\mathfrak M+df.
\end{align*}
If Step (i) gives $\sum_i m_{i\sigma(i)}=\mathfrak M+mf$, then
\begin{align*}
 \mathfrak M+mf
 &=n\sum_i r_{i\sigma(i)}+\mathfrak M+df,
 \qquad d\equiv m\pmod e.
\end{align*}
Since $\gcd(m,e)=1$, writing $\overline m\in\{1,\ldots,e-1\}$ for the residue of $m$ modulo $e$, we obtain
\begin{align*}
 d&=\overline m,\qquad
 \sum_i r_{i\sigma(i)}=\left\lfloor\frac me\right\rfloor,\\
 \sum_i q^{w_{i\sigma(i)}}
 &=\overline m q^n+
 \left(e-\overline m+e\sum_{l=1}^{f-1}q^l\right)q^{n-f}
 =:\nu_m.
\end{align*}
Thus all terms attaining the special lower bound have the same scalar factor
\begin{align*}
 \varpi^{\lfloor m/e\rfloor}\xi_0^{\nu_m},\qquad
 v\left(\varpi^{\lfloor m/e\rfloor}\xi_0^{\nu_m}\right)
 =\frac1{q-1}+\frac me.
\end{align*}

\noindent\textbf{Step (iii).} For the permutations found in Step (ii), we use the leading summand $\varpi^{r_{ij}}\gamma_{ij}^{q^{-nr_{ij}}}$ from Lemma~\ref{lem:special-bound} and compute the signed sum. If permutations attaining the special lower bound exist, every other term has strictly higher valuation.

Once Step (i) is established, Step (ii) concerns only $W$.
The following combinatorial lemma will be used to find the requisite permutations.
Indices in its statement are understood modulo the size of the matrix.

\begin{lem}\label{lem:W-shift}
Let $W=(w_{ij})\in M_n(\{1,\ldots,n\})$, where $w_{ij}\equiv c_0-j+i\pmod n$ for some integer $c_0$.
All indices below are read modulo $n$ and represented in $\{1,\ldots,n\}$.

\noindent\textup{(i)} For $a,d\in\{1,\ldots,n\}$, put $W(a;d)=(w^{(a;d)}_{ij})$, where
\begin{align*}
 w^{(a;d)}_{ij}&=
 \begin{cases}
 w_{a,j-d},&i=a,\\
 w_{ij},&i\ne a.
 \end{cases}
\end{align*}
If $1\le d<n$, there are exactly two permutations $\sigma\in \mathfrak{S}_n$ such that all $w^{(a;d)}_{i\sigma(i)}$ belong to $\{n-1,n\}$.
For both permutations, $n$ occurs $d$ times and $n-1$ occurs $n-d$ times.
If $d=n$, there are exactly two such permutations: one selects only $n$, and the other selects only $n-1$.

\noindent\textup{(ii)} For $a,d\in\{1,\ldots,n\}$, put $W'(a;d)=(w'^{(a;d)}_{ij})$, where
\begin{align*}
 w'^{(a;d)}_{ij}&=
 \begin{cases}
 w_{i+d,a},&j=a,\\
 w_{ij},&j\ne a.
 \end{cases}
\end{align*}
If $1\le d<n$, there are exactly two permutations selecting only $n-1$ and $n$. For both, $n$ occurs $d$ times and $n-1$ occurs $n-d$ times.
If $d=n$, there are exactly two such permutations: one selects only $n$ and the other selects only $n-1$.
In either case, a permutation selecting only $n-1$ exists if and only if $d=n$, and is then unique.
\end{lem}
\begin{prf}
For the original matrix $W$, the condition $w_{ij}\in\{n-1,n\}$ allows exactly the two columns
\begin{align*}
 j\equiv c_0+i\pmod n\quad(w_{ij}=n),\qquad
 j\equiv c_0+i+1\pmod n\quad(w_{ij}=n-1)
\end{align*}
for each row $i$.
In $W(a;d)$ the corresponding choices for the $a$-th row are shifted by $d$.
Once one of these two columns is chosen in row $a$, the requirement that the chosen columns be distinct forces the selections in the successive rows around the cyclic order.
For $1\le d<n$, either initial choice gives a unique permutation, with $d$ entries equal to $n$ and $n-d$ equal to $n-1$.
For $d=n$, the modified row is unchanged, and the only two possibilities are to select $n$ in every row or $n-1$ in every row.
This proves (i).

For (ii), apply the same argument to columns.
In the original matrix, a column $j$ can select $n$ in row $i\equiv j-c_0$ or $n-1$ in row $i\equiv j-c_0-1$ modulo $n$.
The two choices in the $a$-th column are shifted by $-d$.
Injectivity of the selected rows forces the remaining choices around the cyclic order.
For $1\le d<n$ this gives exactly two permutations, each selecting $n$ exactly $d$ times.
For $d=n$ we recover the two constant choices.
The last assertion follows at once.
\end{prf}

\subsubsection*{(I) Constant terms}

This is the value of $h'(\Delta)$ at $\widetilde Y=0$, which equals $\operatorname{qlog}_{\wedge H}\delta(x_0)=0$. For later use we compute its exponent matrix. Put $J=(e-1-j+i)_{1\le i,j\le e}$. We have
\begin{align*}
 M&=\begin{pmatrix}
 fJ+[0]&&&\\
 &fJ+[1]&&\\
 &&\ddots&\\
 &&&fJ+[f-1]
 \end{pmatrix},\\
 \sum_i m_{i\sigma(i)}
 &=f(e-1)ef+e(1+\cdots+(f-1))
 =\mathfrak M
\end{align*}
for all permutations respecting these blocks. Here $[j]$ denotes the $e\times e$ matrix with every entry equal to $j$.

Let $\overline J\in M_e(\{1,\ldots,e\})$ and $\widetilde J\in M_e(\{0,\ldots,e-1\})$ be the matrices congruent to $J$ modulo $e$. Then the residue matrix is
\begin{align*}
 W&=\begin{pmatrix}
 f\overline J&&&\\
 &f\widetilde J+[1]&&\\
 &&\ddots&\\
 &&&f\widetilde J+[f-1]
 \end{pmatrix}.
\end{align*}
The special lower bound is $1/(q-1)$. In the second through the $f$-th blocks, the residue values have the form $f\widetilde J+[l]$ with $1\le l<f$. Hence every selected value must be $n-f+l$, and there is a unique special permutation in each of these blocks. In the first block the two special permutations select respectively only $n$ and only $n-f$.
The degree-0 part itself vanishes, as noted above, irrespective of cancellation among the terms attaining the bound.

\subsubsection*{(II) Contributions from diagonal blocks}

Let $\det h(\Xi)_{\mathrm{diag}}$ denote this part of the determinant. We have
\begin{align*}
 \det h(\Xi)_{\mathrm{diag},1}
 &=\sum_{j=1}^f\left(
 \det h(\Pi^{j-1}\mathfrak{Y}_1)_1
 \prod_{k\ne j}\det h(\Pi^{k-1}\mathfrak{Y}_1)_0\right).
\end{align*}
The calculation of $\det h(\mathfrak{Y}_1)$ is the totally ramified calculation for $L/K_f$ in \cite{MR4126896}. Indeed, $H_0$ is also the standard formal $\mathcal O_{K_f}$-module of height $e$. Let $\wedge' H_0$ be its top exterior power, a formal $\mathcal O_{K_f}$-module of height one, and let $\Delta_{\mathrm{tot}}:\operatorname{Nil}^{\flat,e}\to\operatorname{Nil}^{\flat}$ denote the corresponding determinant morphism. Writing $\wedge'H$ for the associated lift, put
\begin{align*}
 h'_{\mathrm{tot}}(T)&=\sum_{i\in\mathbb Z}
 \frac{(-1)^{(e-1)i}T^{q^{fi}}}{\varpi^i}.
\end{align*}
Then
\begin{align*}
 \det h(\mathfrak{Y}_1)
 &=h'_{\mathrm{tot}}(\Delta_{\mathrm{tot}})
 =h'_{\mathrm{tot}}\bigl(\Delta_{\mathrm{tot}}-_{\flat,\wedge'}\Delta_{\mathrm{tot}}(x_0)\bigr).
\end{align*}
In particular $\det h(\mathfrak{Y}_1)_0=0$, so the summands in $\det h(\Xi)_{\mathrm{diag},1}$ vanish unless $j=1$. Similarly,
\begin{align*}
 \det h(\Xi)_{\mathrm{diag},2}
 &=\det h(\mathfrak{Y}_1)_2\prod_{j>1}\det h(\Pi^{j-1}\mathfrak{Y}_1)_0\\
 &\quad+\sum_{j>1}\det h(\mathfrak{Y}_1)_1\det h(\Pi^{j-1}\mathfrak{Y}_1)_1
 \prod_{k\ne1,j}\det h(\Pi^{k-1}\mathfrak{Y}_1)_0.
\end{align*}
For $j>1$, the degree-zero part $\det h(\Pi^{j-1}\mathfrak{Y}_1)_0$ has a unique leading permutation by (I), hence no cancellation at its leading valuation. Definition 8.2 implies that replacing a degree-zero entry by a positive-degree entry strictly decreases its absolute value. Comparing determinant monomials therefore gives
\begin{align*}
 \left|\det h(\Pi^{j-1}\mathfrak{Y}_1)_1\right|
 &<\left|\det h(\Pi^{j-1}\mathfrak{Y}_1)_0\right|\qquad(j>1).
\end{align*}
Consequently, the second summand of $\det h(\Xi)_{\mathrm{diag},2}$ is strictly smaller than $\det h(\Xi)_{\mathrm{diag},1}$. The same argument excludes all terms with positive degree in a block $j>1$; in the first block, the terms of degree at least three are negligible by \cite{MR4126896}. Thus
\begin{align*}
 \det h(\Xi)_{\mathrm{diag}}
 &\equiv \det h(\mathfrak{Y}_1)_{1,2}
 \prod_{j>1}\det h(\Pi^{j-1}\mathfrak{Y}_1)_0
\end{align*}
modulo terms of strictly higher valuation.

We compute the valuation of the degree-one main terms using the new strategy. The exponent matrix is obtained from the matrix in (I) by replacing its $(1,1)$-block with
\begin{align*}
 M_0&=fJ+f\begin{pmatrix}
 m&0&\cdots&0\\
 \vdots&\vdots&&\vdots\\
 m&0&\cdots&0
 \end{pmatrix}.
\end{align*}
\noindent\textbf{Step (i).} For every relevant permutation,
\begin{align*}
 \sum_i m_{i\sigma(i)}&=\mathfrak M+mf,
\end{align*}
so the special lower bound is $1/(q-1)+m/e$.

\noindent\textbf{Step (ii).} The residue part of $M_0$ is
\begin{align*}
 W_0&=f\overline J'(1;\overline m).
\end{align*}
By Lemma~\ref{lem:W-shift}, exactly two permutations in the first block select only the special values $n$ and $n-f$, using $n$ exactly $\overline m$ times. The remaining blocks uniquely select $n-f+l$ for each $1\le l<f$. Thus the only permutations attaining the special lower bound have the distribution from Step (ii) of the general strategy, and their common scalar factor is
\begin{align*}
 \varpi^{\lfloor m/e\rfloor}\xi_0^{\nu_m},\qquad
 \nu_m=\overline m q^n+
 \left(e-\overline m+e(q+\cdots+q^{f-1})\right)q^{n-f}.
\end{align*}
\noindent\textbf{Step (iii).}
To compute the leading coefficient, we use the totally ramified approximation from \cite{MR4126896}.
Put $r=\lfloor m/e\rfloor$ and
\begin{align*}
 M'_3(m)&=\left(1-\frac{\overline m}{e}\right)q^{fr}
 +\frac{\overline m}{e}q^{f(r+1)}.
\end{align*}
The calculation in \cite{MR4126896} gives
\begin{align*}
 \Delta_{\mathrm{tot}}-_{\flat,\wedge'}\Delta_{\mathrm{tot}}(x_0)
 &\equiv\xi_0^{e q^{f(e-1)}M'_3(m)}f_{\mathrm{tot}}.
\end{align*}
The factor $e$ converts the valuation normalization in \cite{MR4126896} to our choice of $\xi_0$.
The polynomial $f_{\mathrm{tot}}\in\mathbb Z[z'^{q^{-\infty}},y'_{21}{}^{q^{-\infty}},\ldots,y'_{e1}{}^{q^{-\infty}}]$ is the polynomial computed there.
For any integer $a\geq 0$, we write
$f_{\mathrm{tot}}^{(q^{-a})}$ for the polynomial obtained
from $f_{\mathrm{tot}}$ by raising each variable to the
$q^{-a}$-th power, leaving the coefficients unchanged.
Then we have $\overline{f_{\mathrm{tot}}^{(q^{-a})}}=\overline{f_{\mathrm{tot}}}^{q^{-a}}$,
and taking the leading summand of $h'_{\mathrm{tot}}$ gives
\begin{align*}
 h'_{\mathrm{tot}}\bigl(\Delta_{\mathrm{tot}}-_{\flat,\wedge'}\Delta_{\mathrm{tot}}(x_0)\bigr)
 &\equiv (-1)^{(e-1)r}\varpi^r\bigl(\Delta_{\mathrm{tot}}-_{\flat,\wedge'}\Delta_{\mathrm{tot}}(x_0)\bigr)^{q^{-fr}}\\
 &\equiv (-1)^{(e-1)r}\varpi^r\xi_0^{(e-\overline m)q^{n-f}+\overline m q^n}f_{\mathrm{tot}}^{(q^{-fr})}.
\end{align*}
By the special-lower-bound computation in Step (ii), multiplying by the degree-zero contributions from the remaining blocks gives
\begin{align*}
 \det h(\Xi)_{\mathrm{diag}}
 &\equiv (-1)^{(e-1)r}\varpi^r\xi_0^{\nu_m}
 f_{\mathrm{tot}}^{(q^{-fr})}.
\end{align*}
By the computation in
\cite[Theorem~3.23]{MR4126896}, after the corresponding
change of coordinates, we have
\begin{align*}
 f_{\mathrm{tot}}^{(q^{-fr})}
 \equiv (-1)^{m-r+1}(z^{q^f}-z-f_m).
\end{align*}
Combining this with the factor $(-1)^{(e-1)r}$ above,
we obtain
\begin{align*}
 (-1)^{(e-1)r}(-1)^{m-r+1}
 =(-1)^{\overline m+1},
\end{align*}
where $m=er+\overline m$.

\subsubsection*{(III) Contributions from off-diagonal blocks}

\smallskip\noindent\textbf{(a) Case of $m\ge3$ odd.} Put $s=(m-1)/2$. The exponent matrix is
\begin{align*}
 M
 &=f[J]+f\begin{pmatrix}
 0&[s]&\cdots&[s]\\
 [s]&0&\cdots&[s]\\
 \vdots&\vdots&\ddots&\vdots\\
 [s]&[s]&\cdots&0
 \end{pmatrix}
 +\bigl([\overline{j-i}+i-1]\bigr)_{1\le i,j\le f},
\end{align*}
where $[J]\in M_n(\mathbb Z)$ is the block matrix all of whose blocks are $J$, and $\overline{j-i}\in\{0,\ldots,f-1\}$ denotes the residue modulo $f$.

Suppose that a determinant monomial selects $t>0$ off-diagonal $Y$-entries and degree-zero entries elsewhere. Block-row and block-column balance gives $t\ge2$. Let $N_-$ count the selected entries from blocks $(i,j)$ with $j<i$. Since the contributions of $f[J]$ and $[i-1]$ are independent of the permutation and $\sum_{\mathrm{selected}}\overline{j-i}=fN_-$, we have
\begin{align*}
 \sum_i m_{i\sigma(i)}=\mathfrak M+f(st+N_-),\qquad N_-\ge1.
\end{align*}
For $t\ge3$, this exceeds $\mathfrak M+fm$, since $s\ge1$ and $m=2s+1$. For $t=2$, the two blocks must be symmetric and $N_-=1$. By Definition 8.2, any additional diagonal $Y$- or $Z$-entry also has strictly higher valuation. Thus only two off-diagonal $Y$-entries, with degree-zero entries elsewhere, can contribute at the desired valuation. Let $S\subset\mathfrak{S}_n$ be the set of these permutations.

\noindent\textbf{Step (i).} For every $\sigma\in S$, we have $t=2$ and $N_-=1$, so
\begin{align*}
 \sum_i m_{i\sigma(i)}=\mathfrak M+f(2s+1)=\mathfrak M+mf.
\end{align*}
The common special lower bound is $1/(q-1)+m/e$.

\noindent\textbf{Step (ii).} Suppose the two $Y$-entries are selected from block rows $k,l$ with $1\le k<l\le f$. The affected blocks are
\begin{align*}
 \begin{pmatrix}
 fJ+[k-1]&fJ+[fs+l-1]\\
 fJ+[fs+f+k-1]&fJ+[l-1]
 \end{pmatrix}.
\end{align*}
Suppose the $Y$-entry in the upper-right block is chosen from row $a$, and the other from row $b$ of the lower-left block. Put
\begin{align*}
 d_1=\overline{s+1+b-a},\qquad
 d_2=\overline{s+a-b},
\end{align*}
where the bars in this display are representatives in $\{1,\ldots,e\}$. Replacing the $a$-th row of the upper-left block by the $b$-th row of the lower-left block, and the $b$-th row of the lower-right block by the $a$-th row of the upper-right block, the residue matrix for the remaining entries is
\begin{align*}
 \begin{cases}
 \begin{pmatrix}
 f\overline J(a;d_1)&0\\
 0&f\widetilde J(b;d_2)+[l-1]
 \end{pmatrix},&k=1,\\[7pt]
 \begin{pmatrix}
 f\widetilde J(a;d_1)+[k-1]&0\\
 0&f\widetilde J(b;d_2)+[l-1]
 \end{pmatrix},&k\ge2.
 \end{cases}
\end{align*}
These two row replacements change the sign by $-1$. The remaining blocks have the unique special permutations described in (I).

If $k\ge2$, attaining the special lower bound requires every selected entry in each affected block to equal $n-f+k-1$ or $n-f+l-1$, respectively. The observation following Lemma~\ref{lem:W-shift} then forces $d_1=d_2=e$, which contradicts
\begin{align*}
 d_1+d_2\equiv m\not\equiv0\pmod e.
\end{align*}
If $k=1$, the second affected block requires $d_2=e$, and hence $d_1=\overline m$. Thus the special lower bound is attained precisely when
\begin{align*}
 k=1,\qquad (d_1,d_2)=(\overline m,e),\qquad
 b-a\equiv s\pmod e.
\end{align*}
For each $l=2,\ldots,f$, there are $e$ possible pairs $(a,b)$; Lemma~\ref{lem:W-shift} supplies two special permutations from $f\overline J(a;\overline m)$ and one from $f\widetilde J(b;e)+[l-1]$. Hence exactly $2e(f-1)$ terms attain the special lower bound. They all have the distribution established in the general strategy and therefore have common scalar factor $\varpi^{\lfloor m/e\rfloor}\xi_0^{\nu_m}$.

\noindent\textbf{Step (iii). Approximation.} Decompose
\begin{align*}
 J&=J_1+J_2,\quad J_1=(i-1)_{ij},\quad J_2=(e-j)_{ij},\\
 M&=M_1+M_2,\\
 M_1&=f[J_1]+\bigl([i-1]\bigr)_{1\le i,j\le f},\\
 M_2&=f[J_2]+f\begin{pmatrix}
 0&[s]&\cdots&[s]\\
 [s]&0&\cdots&[s]\\
 \vdots&\vdots&\ddots&\vdots\\
 [s]&[s]&\cdots&0
 \end{pmatrix}
 +\bigl([\overline{j-i}]\bigr)_{1\le i,j\le f}.
\end{align*}
Let $m_{ij}=m^{(1)}_{ij}+m^{(2)}_{ij}=nr_{ij}+w_{ij}$ denote a selected exponent corresponding to an off-diagonal $Y$-entry. The corresponding entry of $h(\Xi)$ has the form $h(\Pi^{m^{(1)}_{ij}}Y)$, with $v(Y)=cq^{m^{(2)}_{ij}}$. For a term attaining the special lower bound, $w_{ij}$ is special, so Lemma~\ref{lem:special-bound} yields
\begin{align*}
 h(\Pi^{m^{(1)}_{ij}}Y)
 &\equiv\varpi^{r_{ij}}\bigl(Y^{q^{m^{(1)}_{ij}}}\bigr)^{q^{-nr_{ij}}}
 =\varpi^{r_{ij}}Y^{q^{w_{ij}-m^{(2)}_{ij}}}.
\end{align*}
Thus, in the primed coordinates, its contribution apart from the common scalar factor is of the form $(y'_{ij})^{q^{w_{ij}-m^{(2)}_{ij}}}$.

For the $(1,l)$-block, the condition for the single surviving $Y$-entry in row $a$ and column $t$ is
\begin{align*}
 e-1-t+a+s\equiv e-1\pmod e,
 \qquad t\equiv a+s\pmod e.
\end{align*}
Writing $b=\overline{a+s}\in\{1,\ldots,e\}$, this term contributes $(y'_{b,l})^{q^{f(b-s-1)}}$.
For the $(l,1)$-block, the two surviving entries from row $b$ occur in columns $\overline{b+s}$ and $\overline{b+s+1}$, with opposite signs.
Substituting the special exponents $w-m^{(2)}$ yields the polynomial $g_m'$, namely
\[
(-1)^{\overline{m}-1}\sum_{b=1}^e
\begingroup
\renewcommand{\arraystretch}{1.4}
    \begin{vmatrix}
        1 & y_{b,2}'^{q^{f(b-s-1)}} & y_{b,3}'^{q^{f(b-s-1)}} & \cdots & y_{b,f}'^{q^{f(b-s-1)}}\\
        y_{b+s,f}'^{q^{f(\overline{b+s}-s-1)+1}}-y_{b+s+1,f}'^{q^{f(\overline{b+s+1}-s-2)+1}} & 1 & 0 & \cdots & 0\\
        y_{b+s,f-1}'^{q^{f(\overline{b+s}-s-1)+2}}-y_{b+s+1,f-1}'^{q^{f(\overline{b+s+1}-s-2)+2}} & 0 & 1 & \cdots & 0\\
        \vdots & \vdots & \vdots & \ddots & \vdots \\
        y_{b+s,2}'^{q^{f(\overline{b+s}-s-1)+f-1}}-y_{b+s+1,2}'^{q^{f(\overline{b+s+1}-s-2)+f-1}} & 0 & 0 & \cdots & 1
    \end{vmatrix}^*
\endgroup
\]
where $|\cdot|^{*}$ denotes the sum of the non-identity permutation terms of the determinant, as in Definition~\ref{defn:var}.
The corresponding part of $\det h(\Xi)$ therefore has leading term $\varpi^{\lfloor m/e\rfloor}\xi_0^{\nu_m}g'_m$.
The change of variables
\begin{align*}
 y_{b,l}'^{q^{f(b-s-1)+1-l}}=:y_{b,l}
\end{align*}
identifies $g'_m$ with $(-1)^{\overline{m}}g_m$ from Definition~\ref{defn:var}.

\smallskip\noindent\textbf{(b) Case of $m\ge2$ even.} Put $s=m/2$. The exponent matrix is
\begin{align*}
 M=f[J]+f\begin{pmatrix}
 0&[s]&\cdots&[s]\\
 [s]&0&\cdots&[s]\\
 \vdots&\vdots&\ddots&\vdots\\
 [s]&[s]&\cdots&0
 \end{pmatrix}+\bigl([i-1]\bigr)_{1\le i,j\le f}.
\end{align*}
If $t>0$ off-diagonal $Y$-entries are selected, block balance gives $t\ge2$, and the exponent sum, for degree-zero choices elsewhere, is
\begin{align*}
 \sum_i m_{i\sigma(i)}=\mathfrak M+fst.
\end{align*}
Since $s\ge1$, the special lower bound exceeds $1/(q-1)+m/e$ when $t\ge3$. An additional diagonal $Y$- or $Z$-entry also gives strictly higher valuation by Definition 8.2. Thus we need only two off-diagonal $Y$-entries from symmetric blocks. Let $S\subset\mathfrak S_n$ be the corresponding permutations.

\noindent\textbf{Step (i).} For each $\sigma\in S$, $t=2$, and hence
\begin{align*}
 \sum_i m_{i\sigma(i)}=\mathfrak M+mf.
\end{align*}
The common special lower bound is $1/(q-1)+m/e$.

\noindent\textbf{Step (ii).} Select the two $Y$-entries from block rows $k,l$, with $1\le k<l\le f$. The affected blocks are
\begin{align*}
 \begin{pmatrix}
 fJ+[k-1]&fJ+[fs+k-1]\\
 fJ+[fs+l-1]&fJ+[l-1]
 \end{pmatrix}.
\end{align*}
Suppose that the entry in the upper-right block is chosen from column $a$, and that in the lower-left block from column $b$. Put
\begin{align*}
 d'_1=\overline{s+b-a},\qquad d'_2=\overline{s+a-b},
\end{align*}
where the bars mean representatives in $\{1,\ldots,e\}$. Interchanging the affected columns changes the sign by $-1$, and the residue matrices for the remaining choices become
\begin{align*}
 \begin{cases}
 \begin{pmatrix}
 f\overline J'(b;d'_1)&0\\
 0&f\widetilde J'(a;d'_2)+[l-1]
 \end{pmatrix},&k=1,\\[7pt]
 \begin{pmatrix}
 f\widetilde J'(b;d'_1)+[k-1]&0\\
 0&f\widetilde J'(a;d'_2)+[l-1]
 \end{pmatrix},&k\ge2.
 \end{cases}
\end{align*}
If $k\ge2$, selecting only special values forces $d'_1=d'_2=e$ by Lemma~\ref{lem:W-shift}, contradicting $d'_1+d'_2\equiv m\not\equiv0\pmod e$. If $k=1$, it forces $d'_2=e$, so $d'_1=\overline m$. Consequently the special lower bound is attained exactly when
\begin{align*}
 k=1,\qquad (d'_1,d'_2)=(\overline m,e),\qquad b-a\equiv s\pmod e.
\end{align*}
For each $l=2,\ldots,f$ there are $e$ such pairs $(a,b)$, and Lemma~\ref{lem:W-shift} gives two special permutations in the first affected block and one in the second. Thus there are $2e(f-1)$ contributing terms, all with the common scalar factor $\varpi^{\lfloor m/e\rfloor}\xi_0^{\nu_m}$.

\noindent\textbf{Step (iii). Approximation.} Write $J=J_1+J_2$ as in (a), and decompose
\begin{align*}
 M&=M_1+M_2,\\
 M_1&=f[J_1]+\bigl([i-1]\bigr)_{1\le i,j\le f},\\
 M_2&=f[J_2]+f\begin{pmatrix}
 0&[s]&\cdots&[s]\\
 [s]&0&\cdots&[s]\\
 \vdots&\vdots&\ddots&\vdots\\
 [s]&[s]&\cdots&0
 \end{pmatrix}.
\end{align*}
By Lemma~\ref{lem:special-bound} the same entrywise approximation as in (a) applies. In the $(l,1)$-block, the unique surviving $Y$-entry lies in internal row $a=\overline{b-s}$ and contributes
\begin{align*}
 (y'_{\overline{a+s},l^*})^{q^{f(\overline{a+s}-s-1)+l-1}}.
\end{align*}
In the $(1,l)$-block, the two surviving choices in column $a$ occur in internal rows $\overline{a-s+1}$ and $\overline{a-s}$, with opposite signs. Their contributions are
\begin{align*}
 (y'_{a,l})^{q^{f(a-s)}}\qquad\text{and}\qquad
 (y'_{a,l})^{q^{f(a-s-1)}},
\end{align*}
respectively. Thus the leading coefficient of the off-diagonal part is
\begin{align*}
    (-1)^{\overline{m}-1}\sum_{a=1}^e
    \begingroup
    \renewcommand{\arraystretch}{1.4}
        \begin{vmatrix}
            1 & y_{a+s,f}'^{q^{f(\overline{a+s}-s-1)+1}} & y_{a+s,f-1}'^{q^{f(\overline{a+s}-s-1)+2}} & \cdots & y_{a+s,2}'^{q^{f(\overline{a+s}-s-1)+f-1}}\\
            y_{a,2}'^{q^{f(a-s)}}-y_{a,2}'^{q^{f(a-s-1)}} & 1 & 0 & \cdots & 0\\
            y_{a,3}'^{q^{f(a-s)}}-y_{a,3}'^{q^{f(a-s-1)}} & 0 & 1 & \cdots & 0\\
            \vdots & \vdots & \vdots & \ddots & \vdots \\
            y_{a,f}'^{q^{f(a-s)}}-y_{a,f}'^{q^{f(a-s-1)}} & 0 & 0 & \cdots & 1
        \end{vmatrix}^*
    \endgroup
\end{align*}
which we call $g_m'$.
The leading term is therefore $\varpi^{\lfloor m/e\rfloor}\xi_0^{\nu_m}g'_m$. The change of variables
\begin{align*}
 y_{a,l}'^{q^{f(a-s-1)}}=:y_{a,l}
\end{align*}
identifies $g'_m$ with $(-1)^{\overline{m}}g_m$ in Definition~\ref{defn:var}.

\smallskip\noindent\textbf{(c) Case of $m=1$.} Here more than two off-diagonal $Y$-entries may contribute. The exponent matrix is
\begin{align*}
 M=f[J]+\bigl([\overline{j-i}+i-1]\bigr)_{1\le i,j\le f}.
\end{align*}
Index rows and columns by $(r,a)$ and $(c,b)$, respectively, where $1\le r,c\le f$ and $1\le a,b\le e$. Then
\begin{align*}
 m_{(r,a),(c,b)}=f(e-1-b+a)+(c-1)+f\,\boldsymbol{1}_{c<r}.
\end{align*}

\noindent\textbf{Step (i).} Let $k$ be the number of selected entries from blocks below the block diagonal. As in (a), $\sum_{\mathrm{selected}}(c-r)=0$, so
\begin{align*}
 \sum_{r,a}m_{(r,a),\sigma(r,a)}=\mathfrak M+fk.
\end{align*}
For a term containing an off-diagonal block we have $k\ge1$. Its special lower bound is therefore $1/(q-1)+k/e$, so all terms with $k\ge2$ have strictly higher valuation than $1/(q-1)+1/e$. Any additional diagonal $Y$- or $Z$-entry is also negligible by Definition 8.2. We henceforth take $k=1$.

\noindent\textbf{Step (ii).} By the general strategy, a term attains this bound if and only if it selects $n$ once and $n-f$ exactly $e-1$ times in the first block column, and $n-f+c-1$ exactly $e$ times in every block column $c\ge2$.
The sole selected entry below the block diagonal lies in the first block column. Otherwise that entire column would be selected from the $(1,1)$-block; its special entries form the two cyclic matchings selecting either all $n$ or all $n-f$, neither of which has exactly one $n$ since $e>1$.

Suppose the exceptional entry belongs to the $(l,1)$-block, $2\le l\le f$, and has internal indices $(i,j)$. There are precisely two cases (indices are read modulo $e$):
\begin{itemize}
 \item \textbf{Case A:} $(i,j)=(i,i)$, with residue $n$. All other first-column choices are $(a,a)$ in the $(1,1)$-block for $a\ne i$.
 \item \textbf{Case B:} $(i,j)=(i,\overline{i+1})$, with residue $n-f$. The $(1,1)$-block contributes $(\overline{i+1},i)$, with residue $n$, and $(a,a)$ for $a\ne i,\overline{i+1}$.
\end{itemize}
In each block column $c\ge2$ the selected special entries on or above the block diagonal have equal internal row and column indices. Thus the other off-diagonal choices form an increasing sequence of block indices
\begin{align*}
 1=r_0<r_1<\cdots<r_t=l,
\end{align*}
with one selected entry in each $(r_j,r_{j+1})$-block; the other entries are diagonal. In Case B the two first-column choices are obtained from Case A by interchanging their columns. Hence these and only these permutations attain the special lower bound.

\noindent\textbf{Step (iii). Approximation.} Decompose $J=J_1+J_2$ and $M=M_1+M_2$ as in (a), with
\begin{align*}
 M_1=f[J_1]+\bigl([i-1]\bigr)_{1\le i,j\le f},\qquad
 M_2=f[J_2]+\bigl([\overline{j-i}]\bigr)_{1\le i,j\le f}.
\end{align*}
By Lemma~\ref{lem:special-bound}, the entrywise leading terms are obtained
as in (a). For an upper-diagonal block $(r,c)$ with
$r<c$, the entry with internal index $i$ contributes
$(y'_{i,c-r+1})^{q^{f(i-1)+r-1}}$.
The two choices in the $(l,1)$-block contribute
\begin{align*}
y_{i,f+2-l}'^{q^{f(i-1)+l-1}},
\qquad
y_{\overline{i+1},f+2-l}'^{q^{f(\overline{i+1}-2)+l-1}},
\end{align*}
with opposite signs.
Summing over the increasing block sequences and the
two cases, we obtain $g_1'$ as
\begin{align*}
\sum_{i=1}^e
\begin{vmatrix}
1
&y_{i,2}'^{q^{f(i-1)}}
&y_{i,3}'^{q^{f(i-1)}}
&\cdots
&y_{i,f}'^{q^{f(i-1)}}\\
y_{i,f}'^{q^{f(i-1)+1}}
-y_{\overline{i+1},f}'^{q^{f(\overline{i+1}-2)+1}}
&1
&y_{i,2}'^{q^{f(i-1)+1}}
&\cdots
&y_{i,f-1}'^{q^{f(i-1)+1}}\\
y_{i,f-1}'^{q^{f(i-1)+2}}
-y_{\overline{i+1},f-1}'^{q^{f(\overline{i+1}-2)+2}}
&0
&1
&\cdots
&y_{i,f-2}'^{q^{f(i-1)+2}}\\
\vdots
&\vdots
&\vdots
&\ddots
&\vdots\\
y_{i,2}'^{q^{f(i-1)+f-1}}
-y_{\overline{i+1},2}'^{q^{f(\overline{i+1}-2)+f-1}}
&0
&0
&\cdots
&1
\end{vmatrix}^{*}.
\end{align*}
The corresponding leading part of $\det h(\Xi)$
is $\xi_0^{\nu_1}g'_1$.
The change of variables
\begin{align*}
y_{i,j}'^{q^{f(i-1)+1-j}}=:y_{i,j}
\end{align*}
identifies $g'_1$ with $-g_1$ in Definition~\ref{defn:var}.

\subsubsection*{(IV) Combining the leading terms}
Let $\epsilon\in\{\pm1\}$ be the sign determined by
the row and column permutations above, so that
\[
 \epsilon h'(\Delta)=\det h(\Xi)\det A.
\]
Define
\[
 f'_{\mathrm{tot}}
 :=(-1)^{(e-1)r}f_{\mathrm{tot}}^{(q^{-fr})}.
\]
Combining the computations in (II) and (III), we obtain
\[
 \epsilon h'(\Delta)
 \equiv \varpi^r\tau(f'_{\mathrm{tot}}+g'_m)
\]
modulo terms of strictly higher valuation.
The preceding estimates show that
\[
 f:=(\varpi^r\tau)^{-1} h'(\Delta)\in S_+,
 \quad
 \overline f
 =\epsilon\,
 \overline{f'_{\mathrm{tot}}+g'_m}.
\]
Suppose that $e>1$.
By the calculation in (II) and (III),
up to change of variables, we have
\[
 \overline{f'_{\mathrm{tot}}}
 =(-1)^{\overline m+1}(z^{q^f}-z-f_m),\quad
 \overline{g'_m}=(-1)^{\overline m}g_m.
\]
Consequently,
\begin{align*}
 \overline f
 &=\epsilon(-1)^{\overline m+1}
   (z^{q^f}-z-f_m-g_m),
\end{align*}
which implies
\[
 V(\overline f)^{\perf}
 \cong X_m^{\perf}.
\]

\subsubsection*{(V) Connected component}
In this part, we also consider the case $e=1$ in order to clarify the difference between the cases $e=1$ and $e>1$.
Our aim is to prove that there exists $g\in S_+$ such that
\begin{align*}
    \Delta-_{\flat,\wedge}t
    =[\varpi^r]_{\flat,\wedge}(\tau g),\quad
    f\equiv\begin{cases}
        g & \text{if $e>1$},\\
        g-g^q & \text{if $e=1$}.
    \end{cases}
\end{align*}
Then the argument at the beginning of this subsection
proves Theorem~\ref{thm:red}.
Since
\[
 h'(\Delta)=\log_{\wedge H}\delta_0
     =\varpi^r\tau f,
 \qquad
 \log_{\wedge H}\delta_0
     =\varpi^r\log_{\wedge H}\delta_r,
\]
we have $\log_{\wedge H}\delta_r=\tau f$.
Let $\exp_{\wedge H}(T)\in C[[T]]$ denote the formal inverse of
$\log_{\wedge H}(T)$.

\begin{lem}
We have
\[
 \exp_{\wedge H}(T)\in t\mathcal O_C[[T/t]],
 \qquad
 \exp_{\wedge H}(T)\notin t\mathcal O_C\langle T/t\rangle.
\]
In particular, for any nonzero $a\in\mathcal O_C$ with $|a|<|t|$,
we have
\[
 \exp_{\wedge H}(T)\in a\mathcal O_C\langle T/a\rangle.
\]
\end{lem}

\begin{prf}
Put
\[
 L_t(T):=\frac{1}{t}\log_{\wedge H}(tT),
 \qquad
 E_t(T):=\frac{1}{t}\exp_{\wedge H}(tT).
\]
These power series are formal inverses of each other. We have
\[
 L_t(T)=\sum_{i\geq0}
 \frac{(-1)^{(n-1)i}t^{q^i-1}}{\varpi^i}T^{q^i}
 \in\mathcal O_C[[T]].
\]
Consequently, $E_t(T)\in\mathcal O_C[[T]]$, and hence
$\exp_{\wedge H}(T)\in t\mathcal O_C[[T/t]]$.
Writing $u=t^{q-1}/\varpi\in\mathcal O_C^\times$, we obtain
\[
 \overline{L_t}(T)
   =T+(-1)^{n-1}\overline u\,T^q
 \quad\text{in }\overline{\mathbb F}_q[T].
\]
Since the coefficient of $T^q$ is nonzero, the inverse power series
$\overline{E_t}(T)$ cannot be a polynomial: otherwise the degree of
$\overline{L_t}(\overline{E_t}(T))=T$ would be at least $q$.
Thus $E_t(T)\notin\mathcal O_C\langle T\rangle$, which proves the
second assertion.

Write $E_t(T)=\sum_{i\geq1}a_iT^i$ with $a_i\in\mathcal O_C$
and $a_1=1$. For any nonzero $c\in\mathfrak m_C$,
\[
 \frac{1}{ct}\exp_{\wedge H}(ctT)
  =\frac{1}{c}E_t(cT)
  =\sum_{i\geq1}a_ic^{i-1}T^i
  \in\mathcal O_C\langle T\rangle.
\]
Taking $c=a/t$ proves the final assertion.
\end{prf}

Suppose first that $e>1$. Then $|\tau|<|t|$, so the lemma shows
that $\exp_{\wedge H}(T)\in\tau\mathcal O_C\langle T/\tau\rangle$.
Hence the series
\[
 \widetilde g:=\tau^{-1}\exp_{\wedge H}(\tau f)
\]
converges to an element of $S_+$.
As $\log_{\wedge H}\delta_r=\tau f$, the difference
$\delta_r-_{\wedge H}\tau\widetilde g$ is killed by
$[\varpi^i]_{\wedge H}$ for some $i$.
Since $S_+$ is an integral domain and the $\varpi$-power torsion
points of $\wedge H$ over $C$ form a discrete set, this difference
is constant. Evaluating at the CM point gives
$\delta_r-_{\wedge H}\tau\widetilde g=\widetilde t_r$.
It follows that
$\delta_r-_{\wedge H}\widetilde t_r\in\tau S_+$.

\begin{rem}
When $e=1$, the proof of \cite[Theorem~3.8.1]{MR3402698}
uses the same argument. In this case, however, $|t|=|\tau|$, and
one has only
\[
 \exp_{\wedge H}(T)\in\tau\mathcal O_C[[T/\tau]].
\]
Thus the assertion in the proof of \cite[Theorem~3.8.1]{MR3402698}
that the exponential belongs to $\tau\OO_C\langle T/\tau\rangle$ does not hold in this case,
and the convergence of $\exp_{\wedge H}(\tau f)$ in $S_+$
is not guaranteed by this argument.
Nevertheless, with a slight modification of the argument,
we can still obtain the desired inclusion
$\delta_r-_{\wedge H}\widetilde t_r\in\tau S_+$
as follows.

Since $|\varpi\tau|<|t|=|\tau|$, the series
$\exp_{\wedge H}(\varpi\tau f)$ converges, and we may write
\[
 \exp_{\wedge H}(\varpi\tau f)
   =\varpi\tau\widetilde g'
\]
for some $\widetilde g'\in S_+$.
As
$\varpi\tau f=\log_{\wedge H}\delta_{r-1}$,
the preceding argument gives
\[
 \delta_{r-1}-_{\wedge H}\varpi\tau\widetilde g'
   =\widetilde t_{r-1}.
\]
Consequently,
\[
 [\varpi]_{\wedge H}
 \bigl(\delta_r-_{\wedge H}\widetilde t_r\bigr)
 =\delta_{r-1}-_{\wedge H}\widetilde t_{r-1}
 \in\varpi\tau S_+.
\]
To conclude, note that $|\tau|^{q-1}=|\varpi|$.
For a topologically nilpotent element $x\in S_+$, the valuation of
$[\varpi]_{\wedge H}(x)$ equals $qv(x)$ whenever
$0<v(x)<1/(q-1)$, by the Newton polygon of $[\varpi]_{\wedge H}$.
Since $v(\tau)=1/(q-1)$, this applies when $x\not\in \tau S_+$.
Therefore
$[\varpi]_{\wedge H}(x)\in\varpi\tau S_+$
implies $x\in\tau S_+$.
Applying this to $x=\delta_r-_{\wedge H}\widetilde t_r$
proves the desired inclusion.
\end{rem}

Now let
\[
 (\delta-_{\wedge H}\widetilde t)^{(r)}
 :=(\delta_{i+r}-_{\wedge H}\widetilde t_{i+r})_{i\geq0}
 \in\widetilde{\wedge H},
\]
and write
\[
 (\alpha,\alpha^{1/q},\ldots)\in \Nil^\flat
\]
for the element corresponding to it.
We claim that $\delta_r-_{\wedge H}\widetilde{t}_r\in \tau S_+$ implies $\alpha\in\tau S_+$.
To prove this, set
\[
 z_i:=\delta_{r+i}-_{\wedge H}\widetilde t_{r+i}
 \quad(i\geq0).
\]
Then $[\varpi]_{\wedge H}(z_{i+1})=z_i$ and $z_0\in\tau S_+$.
Let $\|\cdot\|$ denote the Gauss norm on
$S_+$, and put $v(x)=\log\|x\|/\log|\varpi|$ for $x\neq0$,
with $v(0)=+\infty$.
By the definition of $\tau$, we have
\[
 v(\tau)=\frac{1}{q-1}+\frac{\overline m}{e}
   <\frac{q}{q-1}.
\]
As noted above, we have
$v([\varpi]_{\wedge H}(x))=qv(x)$
if $0<v(x)<\frac{1}{q-1}$.
We deduce by induction that
\[
 v(z_i)\geq\frac{v(\tau)}{q^i}\quad(i\geq0).
\]
Indeed, the case $i=0$ is already known. If
$v(z_{i+1})<v(\tau)/q^{i+1}$, then
$v(z_{i+1})<1/(q-1)$, because $v(\tau)<q/(q-1)$.
Hence
\[
 v(z_i)=qv(z_{i+1})<\frac{v(\tau)}{q^i},
\]
contrary to the induction hypothesis.

As described in \cite[equation~(2.4)]{MR4126896},
we have
\[
 \alpha=\lim_{i\to\infty}(-1)^{q(n-1)i} z_i^{q^i}.
\]
Since
$v(z_i^{q^i})=q^iv(z_i)\geq v(\tau)$,
each term lies in $\tau S_+$,
so that $\alpha\in\tau S_+$.
We may therefore write
$\alpha=\tau g$ with $g\in S_+$.
Since
\[
 [\varpi^r]_{\wedge H}
     \bigl((\delta-_{\wedge H}\widetilde t)^{(r)}\bigr)
   =\delta-_{\wedge H}\widetilde t,
\]
we obtain
\[
 \Delta-_{\flat,\wedge}t
   =[\varpi^r]_{\flat,\wedge}(\tau g),
\]
as required.

The identity $h'(\Delta)=\varpi^r\tau f$ equals
$h'(\Delta-_{\fw}t)
=\varpi^r h'(\tau g)$,
so that $\tau f=h'(\tau g)$.
If $e>1$, we have $h'(\tau g)\equiv\tau g$, hence $\overline f=\overline g$.
If $e=1$, as in \cite[proof of Theorem~3.8.1]{MR3402698},
we have $h'(\tau g)\equiv \tau(g-g^q)$, hence $\overline{f}=\overline{g}-\overline{g}^q$.
This completes the proof of Theorem~\ref{thm:red}.

\subsection{Group action on the reduction}
We compute the action of $J_{m,s}^1$ on $\Nil^{\flat,n}$
and prove that $J_{m,s}^1$ stabilizes the affinoid domain $\widetilde{U}_m$,
hence $U_m$ and $U_m^t$.
Also we obtain the explicit formula of the action on $X_m$,
which turns out to be identical with that described in Definition~\ref{defn:action-on-reduction}.
This completes the proof of Theorem~\ref{thm:red}.

We note that any element of $\mathcal{G}^1$ is of the form $(g,d\varphi_D^{-m_\sigma},\sigma)$.
(Recall that $m_\sigma=v_K(\Art_K^{-1}(\sigma))$.)
\begin{lem}\label{lem:action}
    For $(g,d\varphi_D^{-m_\sigma},\sigma)\in \mathcal{G}^1$, we have
    \begin{align*}
        (g,d\varphi_D^{-m_\sigma},\sigma)^*\vec{Y}
        &=\vec{Y}A\begin{pmatrix}
            d^{-1} & & \\
             & \ddots & \\
             & & d^{-1}
        \end{pmatrix}
        g\sigma(A)^{-1}\\
        &\phantom{=}+_\flat \vec{\xi}'A\begin{pmatrix}
            d^{-1} & & \\
             & \ddots & \\
             & & d^{-1}
        \end{pmatrix}
        g\sigma(A)^{-1}
        -_\flat \sigma(\vec{\xi}').
    \end{align*}
    If $(g,d\varphi_D^{-m_\sigma},\sigma)$ stabilizes $\vec{\xi}$,
    then the second line of the right hand side vanishes.
\end{lem}
\begin{prf}
    By the definition of the group action, we have
    \begin{align*}
        (g,d\varphi_D^{-m_\sigma},\sigma)^*\vec{X}
        =(\vec{\xi}'+_\flat \vec{Y})A\begin{pmatrix}
            d^{-1} & & \\
             & \ddots & \\
             & & d^{-1}
        \end{pmatrix}g.
    \end{align*}
    The left hand side equals
    \begin{align*}
        (g,d\varphi_D^{-m_\sigma},\sigma)^*(\vec{\xi}+_\flat \vec{Y}A)
        =\sigma(\vec{\xi})+_\flat ((g,d\varphi_D^{-m_\sigma},\sigma)^*\vec{Y})\sigma(A),
    \end{align*}
    from which the first assertion follows.
    The second assertion is obtained by substituting $\vec{Y}=0$.
\end{prf}

In fact, we can prove a little stronger result saying that
$\J_{\GL}^m$, $\J_D^m$, $(a,a,1)$ for $a\in L^\times$ and $(1,a_\sigma^{-1},\sigma)$ for $\sigma\in W_L$
and $a_\sigma=\Art_L^{-1}(\sigma)$
stabilize $\widetilde{U}_m$,
and we prove them respectively.
We introduce the auxiliary symbol to simplify the formula:
\begin{align*}
    \pi_{ik}=
    \begin{cases}
        \varpi &\text{ if }1\leq i\leq k,\\
        1 &\text{ if }k<i\leq e.
    \end{cases}
\end{align*}
Recall from Subsections~\ref{subsec:GL}, \ref{subsec:D} that,
for $g\in \GL_n(K)$ and $d\in D^\times$, we can write
\begin{align*}
    g=1+\sum_{i=1}^e\sum_{j=0}^{f-1}\gamma_{ij}s_i^j,\quad 
    d^{-1}=1+\sum_{j=0}^{n-1}\beta_j\varphi_D^j,
\end{align*}
where $\gamma_{ij}\in L$ and $\beta_j\in K_n$.
We further put
\begin{align*}
    \gamma_{ij}=\sum_{k=0}^{e-1} x_{ijk}\varphi_L^k\, (x_{ijk}\in K_f),\quad \mu_j=\sum_{i=0}^{e-1}\beta_{fi+j}\varphi_L^i,
\end{align*}
so that $d^{-1}=1+\sum_{j=0}^{f-1}\mu_j\varphi_D^j$.
Under these notations, we have the following corollary.
\begin{cor}\label{cor:action}
    \begin{enumerate}
        \item For $g\in \GL_n(K)$, we have
        \begin{align*}
            &g^*Y_{ij}
            =Y_{ij}+_\flat
            (\flat)\sum_{b=0}^{f-1}\sum_{k=0}^{e-1}[\pi_{ik}\Phi^{j-b-1}(x_{ibk})]_\flat(Y_{i-k,j-b})+_\flat
            [\gamma_{i,j-1}\varphi_L^{e-i}]_\flat (\xi_0),\\
            &g^*Z
            =Z+_\flat (\flat)\sum_{i=1}^e\sum_{b=0}^{f-1}\sum_{k=0}^{e-1}[\varphi_L^{i-1}\pi_{ik}\Phi^{-b}(x_{ibk})]_\flat(Y_{i-k,1-b})+_\flat
            [e\pr_1(g-1)\varphi_L^{e-1}]_\flat (\xi_0).
        \end{align*}
        \item For $d\in D^\times$, we have
        \begin{align*}
            &d^*Y_{ij}
            =Y_{ij}+_\flat (\flat)\sum_{b=0}^{f-1} [\mu_b \varphi_D^b]_\flat (Y_{i,j-b})+_\flat [\mu_{j-1}\varphi_D^{j-1}\varphi_L^{e-i}]_\flat(\xi_0),\\
            &d^*Z
            =Z+_\flat (\flat)\sum_{i=1}^e\sum_{b=0}^{f-1}[\varphi_L^{i-1}\mu_b\varphi_D^b]_\flat(Y_{i,1-b})
            +_\flat [e\pr_2(d^{-1}-1)\varphi_L^{e-1}]_\flat(\xi_0).
        \end{align*}
        \item We have
        \begin{align*}
            &(\varphi_L,\varphi_L,1)^*Y_{ij}=[\varphi_L^{-1}\pi_{i1}]_\flat(Y_{i-1,j}),\quad (\varphi_L,\varphi_L,1)^*Z=Z.
        \end{align*}
        \item For $u\in U_L^1$, we have
        \begin{align*}
            (u,u,1)^*Y_{ij}
            =[u^{-1}]_\flat \Big(Y_{ij}+_\flat \sum_{k=0}^{e-1}[\pi_{ik}\Phi^{j-1}(x_{i0k})]_\flat (Y_{i-k,j})\Big),\quad
            (u,u,1)^*Z=Z.
        \end{align*}
        \item For $u\in U_{K_f}$, we have
        \begin{align*}
            &(u,u,1)^*Y_{ij}=[u^{-1}\Phi^{j-1}(u)]_\flat(Y_{ij}),\quad (u,u,1)^*Z=Z.
        \end{align*}
        \item For $\sigma\in W_L$, we have
        \begin{align*}
            (1,a_\sigma^{-1},\sigma)^*Y_{ij}=[u_\sigma]_\flat(Y_{ij}),\quad
            (1,a_\sigma^{-1},\sigma)^*Z=[u_\sigma]_\flat(Z),
        \end{align*}
        where $u_\sigma=a_\sigma\varphi_L^{-n_\sigma}$.
    \end{enumerate}
\end{cor}
\begin{prf}
(i) By Lemma~\ref{lem:action}, we have
\begin{align*}
    g^*\vec{Y}
        =\vec{Y}A gA^{-1}
        +_\flat \vec{\xi}'(A gA^{-1}-1)
        =\vec{Y}+_\flat (\vec{Y}+_\flat \vec{\xi}')A(g-1)A^{-1}.
\end{align*}
Lemma~\ref{lem:conj} gives
\begin{align*}
    A(g-1)A^{-1}
    &=\sum_{a=1}^e\sum_{b=0}^{f-1} \sum_{k=0}^{e-1} \begin{pmatrix}
         & I_{e-1}\\
        \varpi & 
    \end{pmatrix}^k t_a
    \otimes
    \begin{pmatrix}
        x_{abk} & & \\
         &\ddots & \\
         & & \Phi^{f-1}(x_{abk})
    \end{pmatrix}
    \begin{pmatrix}
         & I_{f-1}\\
        1 &
    \end{pmatrix}^b.
\end{align*}
We substitute this in the above formula of $g^*\vec{Y}$ to get
\begin{align*}
    g^*Y_{ij}
    =Y_{ij}+_\flat
    (\flat)\sum_{b=0}^{f-1}\sum_{k=0}^{e-1}[\pi_{ik}\Phi^{j-b-1}(x_{ibk})]_\flat(Y_{i-k,j-b}+_\flat \xi_{i-k,j-b}').
\end{align*}
We note that $\xi_{ij}'$ is zero unless $j=1$ and $\xi_{i1}'=[\varphi_L^{e-i}]_\flat(\xi_0)$,
so we have
\begin{align*}
    [\pi_{ik}]_\flat(\xi_{i-k,1}')=[\varphi_L^{e-i+k}]_\flat (\xi_0)
\end{align*}
in either case, $1\leq i\leq k$ or $k<i\leq e$.
Therefore we obtain
\begin{align*}
    g^*Y_{ij}
    &=Y_{ij}+_\flat
    (\flat)\sum_{b=0}^{f-1}\sum_{k=0}^{e-1}[\pi_{ik}\Phi^{j-b-1}(x_{ibk})]_\flat(Y_{i-k,j-b})+_\flat
    (\flat)\sum_{k=0}^{e-1} [x_{i,j-1,k}\varphi_L^{e-i+k}]_\flat (\xi_0)\\
    &=Y_{ij}+_\flat
    (\flat)\sum_{b=0}^{f-1}\sum_{k=0}^{e-1}[\pi_{ik}\Phi^{j-b-1}(x_{ibk})]_\flat(Y_{i-k,j-b})+_\flat
    [\gamma_{i,j-1}\varphi_L^{e-i}]_\flat (\xi_0).
\end{align*}
By the definition of $Z$ and the fact that $e\pr_1(g-1)=\sum_{i=1}^e \gamma_{i0}$,
we also obtain the result for $g^*Z$.

(ii) By Lemma~\ref{lem:action}, we have
\begin{align*}
    d^*\vec{Y}
    =\vec{Y}+_\flat (\vec{Y}+_\flat \vec{\xi}')A
    \begin{pmatrix}
        d^{-1}-1 & & \\
         & \ddots & \\
         & & d^{-1}-1
    \end{pmatrix}
    A^{-1}.
\end{align*}
Lemma~\ref{lem:conj} gives
\begin{align*}
    A\begin{pmatrix}
        d^{-1}-1 & & \\
         & \ddots & \\
         & & d^{-1}-1
    \end{pmatrix}
    A^{-1}
    &=\sum_{b=0}^{f-1}
    \mu_b\varphi_D^b I_e\otimes \begin{pmatrix}
         & I_{f-1}\\
        1 & 
    \end{pmatrix}^b,
\end{align*}
from which we deduce the result in (ii).

The remaining cases are also verified in a similar way,
noting that in these cases the actions stabilize the CM point $\vec{\xi}$.
\end{prf}

Put $c=v(\xi_0)$ and suppose that
$m=2s$ is even.
Since $g-1\in\PA^s$, we have
$\gamma_{ib}\in\p_L^s$ for all $i,b$.
Writing
\[
\gamma_{ib}
=\sum_{k=0}^{e-1}x_{ibk}\varphi_L^k,
\]
we obtain
\[
k+e\,v_{K_f}(x_{ibk})\geq s
\]
for every $i,b,k$.
By the definition of $\widetilde{U}_m$,
we have
\[
v(Y_{ij})\geq q^{f(e-i+s)}c,
\qquad
v(Z)\geq q^{f(e-1+m)}c.
\]
First, we consider $g^*Y_{ij}$.
Since
\[
v_L(\pi_{ik})
+e-\overline{i-k}=e-i+k,
\]
we have
\begin{align*}
v\Big(
[\pi_{ik}\Phi^{j-b-1}(x_{ibk})]_\flat
(Y_{i-k,j-b})
\Big)
\geq
q^{f(e-i+s+k+e v_{K_f}(x_{ibk}))}c
>q^{f(e-i+s)}c.
\end{align*}
Moreover, since
$\gamma_{i,j-1}\in\p_L^s$, we have
\[
|[\gamma_{i,j-1}\varphi_L^{e-i}]_\flat
(\xi_0)|
\leq
|[\varphi_L^s]_\flat(\xi_{i1}')|.
\]
Consequently,
\[
|g^*Y_{ij}|
\leq|\xi_{i1}'|^{q^{fs}}.
\]

Next, we consider $g^*Z$.
For every $i,b,k$, we obtain
\begin{align*}
v\Big(
[\varphi_L^{i-1}\pi_{ik}
\Phi^{-b}(x_{ibk})]_\flat
(Y_{\overline{i-k},\overline{1-b}})
\Big)
\geq
q^{f(e-1+s+k+e v_{K_f}(x_{ibk}))}c
\geq q^{f(e-1+m)}c.
\end{align*}
In addition, since
$e\pr_1(g-1)
\in\p_L^m$,
we have
\[
|[e\pr_1(g-1)\varphi_L^{e-1}]_\flat (\xi_0)|
\leq
|[\varphi_L^{m+e-1}]_\flat(\xi_0)|.
\]
Therefore,
\[
|g^*Z|
\leq|\xi_{11}'|^{q^{fm}}.
\]
This proves that $g$ preserves
$\widetilde{U}_m$.
Also the formula of the action of $S_{\GL}^m$ on $X_m$ is obtained by these calculations.

For elements $d\in \J_D^m$, $(a,a,1)$ for $a\in \{\varphi_L\}\cup U_L^1 \cup U_{K_f}$
and $(1,a_\sigma^{-1},\sigma)$ for $\sigma\in W_L$,
we can also check the stabilization of $\widetilde{U}_m$ and the formula of actions on $X_m$
by the similar computations and Corollary~\ref{cor:action},
which we omit.
Since $\varphi_L^\Z U_L^1 U_{K_f}=L^\times$,
we have proved Theorem~\ref{thm:red} for $m\geq 2$ even.
The other cases are treated in the same way.

\begin{prop}\label{prop:disjoint}
    For $g\in J_m^1\setminus J_{m,s}^1$, we have $gU_m\cap U_m=\emptyset$.
\end{prop}
\begin{prf}
    We see that it suffices to prove that every element $x\in L^\times\setminus U_L^m$ satisfies $i^G(x)\widetilde{U}_m\cap \widetilde{U}_m=\emptyset$.
    Let $g\in J_m^1\setminus J_{m,s}^1$ satisfy $gU_m\cap U_m\neq \emptyset$.
    Then, since $J_m^1=(L^\times \J_{\GL}^m\times \{1\}\times \{1\})^1J_{m,s}^1$
    and $J_{m,s}^1$ stabilizes $\widetilde{U}_m$,
    there exists an element $(a,1,1)\in (L^\times \J_{\GL}^m\times \{1\}\times \{1\})^1\setminus J_{m,s}^1$
    such that $(a,1,1)\widetilde{U}_m\cap \widetilde{U}_m\neq \emptyset$.
    We write $a=i^G(x)h$ where $x\in L^\times$ and $h\in \J_{\GL}^m$.
    Since $\J_{\GL}^m$ stabilizes $\widetilde{U}_m$, $i^G(x)\widetilde{U}_m\cap \widetilde{U}_m\neq \emptyset$.
    However, such $x$ should belong to $U_L^m\subset \J_{\GL}^m$,
    which contradicts the assumption $(a,1,1)\notin J_{m,s}^1$.
    
    Let $x\in L^\times$ and $u\in \widetilde{U}_m$ and we suppose that $i^G(x)(u)\in \widetilde{U}_m$.
    We write the corresponding valuation as $|\cdot|_u$.
    We put $v_L(x-1)=l$.
    Since we have
    \begin{align*}
        i^G(x)^*Z=Z+_\flat [x-1]_\flat(Z)+_\flat [e(x-1)\varphi_L^{e-1}]_\flat(\xi_0),
    \end{align*}
    the inequalities $|Z|_u\leq |\xi_{11}'|^{q^{fm}}$ and $|i^G(x)^*Z|_u\leq |\xi_{11}'|^{q^{fm}}$
    implies that 
    \begin{align*}
        |[x-1]_\flat(Z+_\flat [e]_\flat(\xi_{11}'))|_u
        =|Z+_\flat [e]_\flat(\xi_{11}')|_u^{q^{fl}}
        \leq |\xi_{11}'|^{q^{fm}}.
    \end{align*}
    As $|Z|_u\leq |\xi_{11}'|^{q^{fm}}<|[e]_\flat(\xi_{11}')|$,
    we deduce
    \begin{align*}
        |Z+_\flat [e]_\flat(\xi_{11}')|_u^{q^{fl}}
        =|\xi_{11}'|^{q^{fl}}
        \leq |\xi_{11}'|^{q^{fm}},
    \end{align*}
    which implies $l\geq m$ as desired.
\end{prf}

\section{Proof of the main theorem}\label{sec:prf}
The following theorem is essentially contained in \cite{mieda2016geometricapproachexplicitlocal}.
Its formulation in the present form was communicated to us by Yoichi Mieda,
and we include a proof here with his permission.
We fix an element $\gamma\in K^\times$ such that $v_K(\gamma)\geq 1$,
and we consider a discrete subgroup $\Gamma=\gamma^\Z\subset K^\times$.
\begin{thm}\label{thm:Mieda}
Let $J_1$ (resp.\ $J_2$) be an open compact-mod-center subgroup of $\GL_n(K)$
(resp.\ $D^\times$), and $J_3$ a finite index subgroup of $W_K$.
Let $\overline{\pi}$ (resp.\ $\overline{\rho}$) be a finite-dimensional irreducible smooth
representation of $J_1$ (resp.\ $J_2$) with central character trivial on $\Gamma$
such that
\[
\pi = \cInd_{J_1}^{\GL_n(K)} \overline{\pi}
\quad
(\text{resp.\ } \rho = \cInd_{J_2}^{D^\times} \overline{\rho})
\]
is irreducible.
Let $\overline{\tau}$ be a finite-dimensional representation of $J_3$, and put
\[
\tau = \Ind_{J_3}^{W_K} \overline{\tau}.
\]

We put $J = J_1 \times J_2 \times J_3$ and $J^1 = J \cap \mathcal{G}^1$.
Suppose that there exist
\begin{itemize}
  \item a rational subset $U$ of $t(\mathcal{M}_{H_0,\infty,t}^{(0)})$
        contained in $M^{(0)}_{H_0,\infty,t}$,
  \item a flat formal model $\mathcal{X}$ over $\mathcal{O}_C$ of $U$,
  \item a finite index open subgroup $J_s^1$ of $J^1$, and
  \item a smooth equidimensional affine scheme $Y$ over $\overline{\mathbb{F}}_q$
        equipped with a $J_s^1$-action
\end{itemize}
such that
\begin{itemize}
  \item $J_s^1$ stabilizes $U$, and the action of $J_s^1$ on $U$ extends to that on
        $\mathcal{X}$,
  \item for $g \in J^1 \setminus J_s^1$, $gU \cap U = \emptyset$,
  \item $\overline{\mathcal{X}} = \mathcal{X} \otimes_{\mathcal{O}_C}
        \overline{\mathbb{F}}_q$ is $J_s^1$-equivariantly isomorphic to the
        perfection of $Y$, and
  \item the natural map
        \[
        H_c^{n-1}(Y,\overline{\mathbb{Q}}_l)
        \to
        H^{n-1}(Y,\overline{\mathbb{Q}}_l)
        \]
        is an isomorphism.
\end{itemize}

Further, assume one of the following conditions:
\begin{enumerate}
  \item
  The composite $J_s^1 \hookrightarrow J \twoheadrightarrow J_1 \times J_3$ is surjective and
  $\Hom_{J_s^1 \cap J_2}(\overline{\rho}, H_c^{n-1}(Y,\overline{\mathbb{Q}}_l))$
  is isomorphic to $\overline{\pi} \boxtimes \overline{\tau}$ as a representation of
  $J_1 \times J_3$.
  \item
  The composite $J_s^1 \hookrightarrow J \twoheadrightarrow J_2 \times J_3$ is surjective and
  $\Hom_{J_s^1 \cap J_1}(\overline{\pi}, H_c^{n-1}(Y,\overline{\mathbb{Q}}_l))$
  is isomorphic to $\overline{\rho} \boxtimes \overline{\tau}$ as a representation of
  $J_2 \times J_3$.
\end{enumerate}
Then, there exists a non-zero $\GL_n(K)\times D^\times \times W_K$-equivariant map
\[
\pi \boxtimes \rho \boxtimes \tau \to H^{n-1}_{\mathrm{LT},\Gamma}.
\]
\end{thm}
\begin{prf}
    Every cohomology is taken over $\overline{\Q}_l$.
    We set $H_{\LT}':=\varinjlim_m H^{n-1}_c(M_{H_0,m}^{(0)}\widehat{\otimes}_{\breve K_m}C)$.
    We freely use the notations in \cite{mieda2016geometricapproachexplicitlocal}.
    We put \[
    \mathcal G_\Gamma
    =
    \GL_n(K)\times(D^\times/\Gamma)\times W_K,
    \quad
    J_\Gamma=J_1\times J_2/\Gamma\times J_3.
    \]
    Since $\Gamma\cap\mathcal G^1=\{1\}$,
    the composite $\mathcal G^1\hookrightarrow \mathcal G\to \mathcal G_\Gamma$ is injective,
    by which we identify $\mathcal{G}^1$ with its image in $\mathcal{G}_\Gamma$.
    By \cite[Corollary~4.6 (i)]{mieda2016geometricapproachexplicitlocal},
    we have a $J_s^1$-equivariant homomorphism
    \[\mathrm{sp}^*\colon H_c^{n-1}(Y)\to H_{\LT}'.\]
    By \cite[Corollary~4.6 (iii)]{mieda2016geometricapproachexplicitlocal},
    the induced map
    \[\Ind_{J_s^1}^{J^1}H_c^{n-1}(Y)\to H_{\LT}'\]
    is a $J^1$-equivariant injection.

    We assume that we have the condition (i) in the theorem.
    Applying \cite[Lemma~4.8 (iii)]{mieda2016geometricapproachexplicitlocal}
    to $G=J_1\times J_3$, $H=J_2/\Gamma$, $\Delta=J_s^1$ and $\Delta'=J^1$,
    we obtain the $J_1\times J_3$-equivariant isomorphism
    \[\Hom_{J_s^1\cap (J_2/\Gamma)}(\overline{\rho}, H_c^{n-1}(Y))^{J_1\times J_3\text{-sm}}
    \cong \Hom_{J^1\cap (J_2/\Gamma)}(\overline{\rho}, \Ind_{J_s^1}^{J^1} H_c^{n-1}(Y))^{J_1\times J_3\text{-sm}}.\]
    This injects into $\Hom_{J^1\cap (J_2/\Gamma)}(\overline{\rho}, H_{\LT}')^{J_1\times J_3\text{-sm}}$.
    By \cite[Lemma~4.8 (ii)]{mieda2016geometricapproachexplicitlocal},
    this space is $J_1\times J_3$-equivariantly isomorphic to 
    $\Hom_{J_2/\Gamma} (\overline{\rho}, \Ind_{J^1}^{J_\Gamma}H_{\LT}')$.
    By \cite[Lemma~4.7 (ii)]{mieda2016geometricapproachexplicitlocal},
    since $J_\Gamma$ is open in $\mathcal G_\Gamma$ and $\mathcal G^1$ is closed in $\mathcal G_\Gamma$,
    we have a $J_\Gamma$-equivariant isomorphism $\Ind_{J^1}^{J_\Gamma} H_{\LT}'\cong \Ind_{\mathcal{G}^1}^{\mathcal{G}^1J_\Gamma} H_{\LT}'$.
    Further, by (i) of the same lemma, this injects into $\Ind_{\mathcal{G}^1}^{\mathcal{G}_\Gamma} H_{\LT}'$,
    which is isomorphic to $H_{\LT,\Gamma}=H_{\LT,\Gamma}^{n-1}$.
    Therefore, the condition (i) gives a $J_1\times J_3$-equivariant injection
    $\overline{\pi}\boxtimes \overline{\tau}\hookrightarrow \Hom_{J_2/\Gamma}(\overline{\rho},H_{\LT,\Gamma})$,
    which gives a nonzero $J_\Gamma$-equivariant homomorphism
    \[\overline{\pi}\boxtimes \overline{\rho}\boxtimes \overline{\tau}\to H_{\LT,\Gamma}.\]
    Then the desired homomorphism is obtained by the Frobenius reciprocity.
\end{prf}

\begin{prop}\label{prop:LT-twist}
    For an unramified character $\chi$ of $K^\times$, we have
    \[H^{n-1}_{\mathrm{LT}}\otimes ((\chi\circ \det)\boxtimes (\chi^{-1}\circ \Nrd)\boxtimes (\chi^{-1}\circ \Art_K^{-1}))\cong H^{n-1}_{\mathrm{LT}}\]
    as a representation of $\GL_n(K)\times D^\times \times W_K$.
\end{prop}
\begin{prf}
    Let $\mathcal{G}^0$ be the kernel of the homomorphism 
    \[\mathcal G=\GL_n(K)\times D^\times \times W_K\to \Z; (g,d,\sigma)\mapsto v_K(\det(g)^{-1}\Nrd(d)\Art_K^{-1}(\sigma)).\]
    Then, $\mathcal{G}^0$ equals the stabilizer of $M_{H_0,m}^{(0)}\subset M_{H_0,m}$ for $m\geq 0$,
    and we have $H_{\mathrm{LT}}^{n-1}\cong \cInd_{\mathcal{G}^0}^G (\varinjlim_m H_c^{n-1}(M_{H_0,m}^{(0)}\otimes_{\breve{K}}C, \overline{\Q}_l))$.
    Since $\chi$ is unramified, $(\chi\circ \det)\boxtimes (\chi^{-1}\circ \Nrd)\boxtimes (\chi^{-1}\circ \Art_K^{-1})$
    is trivial on $\mathcal{G}^0$.
    Therefore the statement follows from the Frobenius reciprocity.
\end{prf}

\begin{prf}[of Theorem~\ref{thm:main} (ii)]
We apply Theorem 9.1 with
\[
J_1=L^\times \J_{\GL}^m,\quad
J_2=L^\times \J_D^m,\quad
J_3=W_L,
\]
and
\[
\overline{\pi}=\overline{\pi}_{\GL,\theta}^{\vee},\quad
\overline{\rho}=\overline{\pi}_{D,\theta},\quad
\overline{\tau}=\theta,\quad
U=U_m^t,\quad
\mathcal{X}=\mathcal{X}_m^t.
\]
The flatness of $\mathcal{X}_m^t$ follows by \cite[Corollary~A.6]{mieda2016geometricapproachexplicitlocal}.
Here and below, we regard $\theta$ as a character of $W_L$ via
$\operatorname{Art}_L^{-1}$. Put
$H=H_c^{n-1}(X_m^0,\overline{\mathbb Q}_l)$.

We first assume that the central characters of $\overline{\pi}$ and $\overline{\rho}$
are trivial on $\Gamma$.
By Theorem~\ref{thm:main} (i), we have an isomorphism
\begin{equation}\label{eq:hom-ulm1}
\Hom_{U_L^{m,1}}(\theta,H)
\cong \overline{\pi}\boxtimes \overline{\rho}
\end{equation}
of $J_1\times J_2$-representations.
In order to apply condition (i) of Theorem~\ref{thm:Mieda}, it suffices
to show that
\begin{equation}\label{eq:hom-j2}
\Hom_{J_{m,s}^1\cap J_2}(\overline\rho,H)
\cong \overline{\pi}\boxtimes\overline{\tau}
\end{equation}
as $J_1\times J_3$-representations, since the projection
$J_{m,s}^1\to J_1\times J_3$ is surjective by Lemma~\ref{lem:surj}.

We write
\[
\J_D^{m,1}=\J_D^m\cap \Ker(\Nrd)=J_{m,s}^1\cap J_2.
\]
By the equality
$\Nrd(\J_D^m)=N_{L/K}(U_L^m)$,
we have
$\J_D^m=\J_D^{m,1}U_L^m$.
By Theorem~\ref{thm:characterize} (ii), the restriction $\overline \rho|_{\J_D^m}$ is
irreducible and $U_L^m$ acts on $\overline \rho$ by the scalar
character $\theta$. It follows that $\overline \rho|_{\J_D^{m,1}}$
is irreducible.

For each $u\in U_L^{m,1}$, the actions
of $(1,1,\operatorname{Art}_L(u))$ and $(1,u,1)$
on $H$ coincide.
Let $H[\theta]$ denote the $\theta$-isotypic component
of $H$ for this action.
Then
\begin{align*}
    \Hom_{\J_D^{m,1}}(\overline{\rho},H)
    =\Hom_{\J_D^{m,1}}(\overline{\rho},H[\theta]).
\end{align*}

The canonical evaluation map
\[
\Hom_{U_L^{m,1}}(\theta,H)\otimes\theta
\longrightarrow H[\theta]
\]
is a $J_{m,s}^1$-equivariant isomorphism, as used
in the proof of Theorem~\ref{thm:characterize}.
Together with \eqref{eq:hom-ulm1}, this gives
\[
H[\theta]\cong
(\overline{\pi}\boxtimes\overline{\rho}\boxtimes\overline{\tau})
|_{J_{m,s}^1}.
\]
Thus, by the irreducibility of $\overline{\rho}|_{\J_D^{m,1}}$ we obtain
\begin{align*}
\Hom_{\J_D^{m,1}}(\overline{\rho},H)
&\cong
\Hom_{\J_D^{m,1}}
   (\overline{\rho},\overline{\pi}\boxtimes\overline{\rho}\boxtimes\overline{\tau})\\
&\cong
(\overline{\pi}\boxtimes\overline{\tau})\otimes
\operatorname{End}_{\J_D^{m,1}}(\overline{\rho})\\
&\cong \overline{\pi}\boxtimes\overline{\tau}.
\end{align*}
These isomorphisms are compatible with the action of
$J_{m,s}^1$, and hence are $J_1\times J_3$-equivariant
by Lemma~\ref{lem:surj}.
This proves \eqref{eq:hom-j2}.

The remaining hypotheses of Theorem 9.1 are also
satisfied. The required irreducibility statements
follow from Proposition~\ref{prop:comparison},
the affinoid and its
equivariant formal model are provided by Theorem~\ref{thm:red},
and the disjointness condition follows from
Proposition~\ref{prop:disjoint}.
The morphism
\[
H_c^{n-1}(X_m^0,\overline{\Q}_l)
\longrightarrow
H^{n-1}(X_m^0,\overline{\Q}_l)
\]
is an isomorphism by Proposition~\ref{prop:forget-supp}.
Thus Theorem~\ref{thm:Mieda} gives a nonzero
$\GL_n(K)\times D^\times\times W_K$-equivariant map
\begin{align*}
\pi_{\GL,\theta}^{\vee}
\boxtimes \pi_{D,\theta}
\boxtimes \Ind_{W_L}^{W_K}\theta
\longrightarrow H_{\mathrm{LT},\Gamma}^{n-1}.
\end{align*}
By Theorem~\ref{thm:NALT} (ii), together with the irreducibility
of the three factors, we deduce
\[
\JL(\pi_{\GL,\theta})
\cong\pi_{D,\theta},\quad
\rec_K(\pi_{\GL,\theta})
\cong
\Ind_{W_L}^{W_K}\theta
\Big(\frac{n-1}{2}\Big).
\]

We now treat the general case, without the assumption
on the central characters.
First observe that the central characters of
$\overline{\pi}_{\GL,\theta}$ and
$\overline{\pi}_{D,\theta}$ coincide on $K^\times$.
Indeed, the diagonal subgroup
$\{(a,a,1)\mid a\in K^\times\}$
acts trivially on $X_m^0$ by Definition~\ref{defn:action-on-reduction}.
The assertion follows from the
$J_{m,s}^1$-equivariant isomorphism
\[
H[\theta]\cong
\left.
\left(
\overline{\pi}_{\GL,\theta}^{\vee}
\boxtimes
\overline{\pi}_{D,\theta}
\boxtimes\theta
\right)
\right|_{J_{m,s}^1}.
\]
We denote their common central character by $\omega$.

We choose an unramified character $\chi$ of $K^\times$
such that 
$\chi(\gamma)^n=\omega(\gamma)^{-1}$.
We write
\begin{align*}
\overline{\pi}\chi^{-1}&=\overline{\pi}\otimes(\chi^{-1}\circ\det),\quad
\overline{\rho}\chi=\overline{\rho}\otimes(\chi\circ\operatorname{Nrd}),\quad
\overline{\tau}\chi=\overline{\tau}\otimes(\chi\circ\Art_K^{-1})|_{W_L}.
\end{align*}
Then the central characters of $\overline{\pi}\chi^{-1}$ and $\overline{\rho}\chi$
are trivial on $\Gamma$.
Furthermore, the character
\[
(\chi^{-1}\circ\det)\boxtimes
(\chi\circ\Nrd)\boxtimes
(\chi\circ\Art_K^{-1})
\]
is trivial on $\mathcal{G}^1$.
Therefore the argument above gives
\[
\Hom_{\J_D^{m,1}}(\overline{\rho}\chi,H)
\cong \overline{\pi}\chi^{-1}\boxtimes\overline{\tau}\chi
\]
as $J_1\times J_3$-representations.
Applying Theorem~\ref{thm:Mieda} to $(\overline{\pi}\chi^{-1},\overline{\rho}\chi,\overline{\tau}\chi)$
and then Theorem~\ref{thm:NALT} (ii), we obtain
\begin{align*}
   \JL(\pi^\vee \chi) \cong \rho\chi,\quad
   \rec_K
    (\pi^\vee \chi)
    \cong
    \tau \chi
    \Big(\frac{n-1}{2}\Big).
\end{align*}

Finally, by Proposition~\ref{prop:LT-twist},
we have
\begin{align*}
    &\Hom_{\GL_n(K)} (H_{\LT}^{n-1},\pi^\vee)\\
    &=\Hom_{\GL_n(K)} (H_{\LT}^{n-1},\pi^\vee\chi\otimes (\chi^{-1}\circ \det))\\
    &\cong \Hom_{\GL_n(K)} (H_{\LT}^{n-1},\pi^\vee\chi)\otimes ((\chi^{-1}\circ\Nrd)\boxtimes (\chi^{-1}\circ \Art_K^{-1})),
\end{align*}
and by Theorem~\ref{thm:NALT} (i) this is isomorphic to 
\begin{align*}
    (\rho\chi\boxtimes \tau\chi(n-1))\otimes((\chi^{-1}\circ\Nrd)\boxtimes (\chi^{-1}\circ \Art_K^{-1}))
    \cong \rho\boxtimes \tau(n-1).
\end{align*}
Hence the theorem follows by Theorem~\ref{thm:NALT} (i).
\end{prf}

\bibliographystyle{my_amsalpha}
\bibliography{refs} 

\end{document}